\documentclass[fleqn]{LMCS}

\def\dOi{11(1:10)2015}
\lmcsheading%
{\dOi}
{1--23}
{}
{}
{Feb.~15, 2014}
{Mar.~25, 2015}
{}

\keywords{operad, non-symmetric operad, Cat-operad, weak
Cat-operad, 2-category, bicategory, multicategory, nominal arity,
coherence, symmetric groups, hemiassociahedron, associahedron,
permutohedron}

\ACMCCS{[{\bf Theory of computation}]: Logic---Proof theory; [{\bf
Mathematics of
  computing}]:  Discrete mathematics---Combinatorics---Combinatorics on words;  Discrete mathematics---Graph theory---Hypergraphs}
\amsclass{03B22, 03G30, 08A55, 18A15, 18D05, 18D20, 18D50, 52B10,
55P48}

\usepackage{enumitem,hyperref}

\def\O{$\mathcal{O}$}
\def\Oe{$\mathcal{O}_e$}
\def\Ou{$\mathcal{O}_u$}
\def\WOe{$\mathcal{WO}_e$}
\def\WOu{$\mathcal{WO}_u$}

\def\Om{$\mathcal{O}^-$}
\def\Oem{$\mathcal{O}_e^-$}
\def\Oum{$\mathcal{O}_u^-$}
\def\WOem{$\mathcal{WO}_e^-$}
\def\WOum{$\mathcal{WO}_u^-$}
\def\WOuth{$\mathcal{WO}_u^\theta$}
\def\WOemX{$\mathcal{WO}_e^-(X)$}
\def\WOub{$\mathcal{WO}_u(b)$}
\def\WOubp{$\mathcal{WO}_u({b\:\bullet})$}

\def\Cg{$\mathcal{C}_\Gamma$}
\def\BCg{$\mathcal{BC}_\Gamma$}

\def\Np{${\mathbf N}^+$}

\def\ins{\triangleleft}

\def\cirk{\,{\raisebox{.3ex}{\tiny $\circ$}}\,}
\def\mj{{\mathbf 1}}
\def\pl{\!+\!}
\def\mn{\!-\!}
\def\str{\rightarrow}

\def\Iota{{\mathbf I}}
\def\Monu{${\mbox{\it Mon}}_u$}

\def\bsl{\backslash}

\begin{document}

\title[Weak Cat-Operads]{Weak Cat-Operads}

\author[K.~Do\v sen]{Kosta Do\v sen}  
\address{Mathematical Institute, SANU, Knez Mihailova 36, p.f.\ 367, 11001 Belgrade,
Serbia} 
\email{\{kosta,zpetric\}@mi.sanu.ac.rs}  

\author[Z.~Petri\' c]{Zoran Petri\' c}    
\address{\vspace{-18 pt}} 

\begin{abstract}
\noindent An operad (this paper deals with non-symmetric operads)
may be conceived as a partial algebra with a family of insertion
operations, which correspond to substitution of an operation
within an operation. These insertion operations are Gerstenhaber's
circle-i products, and they satisfy two kinds of associativity,
one of them involving commutativity. A Cat-operad is an operad
enriched over the category Cat of small categories, as a
2-category with small hom-categories is a category enriched over
Cat. This means that the operadic operations of the same arity in
a Cat-operad do not make just a set, but they are the objects of a
small category. The notion of weak Cat-operad is to the notion of
Cat-operad what the notion of bicategory is to the notion of
2-category. This means that the equations of operads like
associativity of insertions are replaced by isomorphisms in a
category. The goal of this paper is to formulate conditions
concerning these isomorphisms that ensure coherence, in the sense
that all diagrams of canonical arrows commute. This is the sense
in which the notions of monoidal category and bicategory are
coherent. (The coherence of monoidal categories, which is due to
Mac Lane, is the best known coherence result.) The coherence proof
in the paper is much simplified by indexing the insertion
operations in a context-independent way, and not in the usual
manner. This proof, which is in the style of term rewriting,
involves an argument with normal forms that generalizes what is
established with the completeness proof for the standard
presentation of symmetric groups. This generalization may be of an
independent interest, and related to other matters than those
studied in this paper. Some of the coherence conditions for weak
Cat-operads lead to the hemiassociahedron, which is a polyhedron
related to, but different from, the three-dimensional
associahedron and permutohedron.
\end{abstract}

\maketitle\vfill

\vspace{2ex}

\begin{tabbing}
{\bf Contents}
\\*[1.5ex]
1. \hspace{.3em} \= {\it Introduction}
\\*[1ex]
{\small \sc Part I}
\\*[.5ex]
2. \> {\it The operad \O}
\\*
3. \> {\it The structure \Oe}
\\
4. \> {\it \O\ and \Oe}
\\
5. \> {\it The structure \Ou}
\\
6. \> {\it \Oe\ and \Ou}
\\[1ex]
{\small \sc Part II}
\\[.5ex]
7. \> {\it The category \WOum}
\\
8. \> {\it The category \WOu}
\\
9. \> {\it The category \WOem}
\\
10. \> {\it The category \WOe}
\\
11. \> {\it \WOe\ and \WOu}
\\
12. \> {\it Operads, Cat-operads and weak Cat-operads}
\\
13. \> {\it \WOem\ and hemiassociahedra}
\\[1ex]
{\small \sc Part III}
\\[.5ex]
14. \> {\it Coherence of \Monu}
\\
15. \> {\it The category \WOuth}
\\
16. \> {\it \Cg\ and \BCg}
\\
17. \> {\it Coherence of \Cg}
\\
18. \> {\it \Cg\ and \WOuth\ --- Coherence of \WOe}
\\[1ex]
{\it References}
\end{tabbing}

\section{Introduction}\label{sec1}
An operad may be conceived as a partial algebra whose elements,
called \emph{operadic operations}, are of various arities; to
these elements as arguments one applies partial binary operations
--- these are operations applied to operadic operations, partiality
being induced by arity --- which we call \emph{insertions}. We
have this new name to distinguish insertions from the related
partial operation of composition in categories, which will appear
together with insertions later in this paper. Insertions
correspond to Gerstenhaber's ``$\cirk_i$-products'' (see
\cite{MSS}, Sections I.1.3 and II.1.3, and the ``circle-$i$'' of
\cite{L04}, Section 2.3) or to Gentzen's cut (see \cite{G35}).

For insertions one assumes two kinds of associativity, one of them
involving commutativity up to a certain point (see the equations
(\emph{assoc}~1) and (\emph{assoc}~2) in Section~2 below, and
related equations given later). One assumes also a unit operadic
operation and appropriate equations tying it to insertions. When
this unit is missing we have a \emph{non-unitary} operad;
otherwise, the operad is \emph{unitary} (for this terminology see
\cite{MSS}, Section II.1.3). This notion of operad, with which we
deal in this paper, is not the original symmetric notion, but the
non-symmetric (non-$\Sigma$) notion (see \cite{MSS}, Section
I.1.3, \cite{M72} and \cite{L04}).

A Cat-operad is an operad enriched in the category Cat of all
small categories, whose arrows are functors, as a 2-category is a
category enriched over Cat (provided the hom-categories of the
2-category are small; see \cite{ML98}, Section XII.3). The
operadic operations of the same arity in a Cat-operad do not make
just a set, but they are the objects of a small category, and the
structure of the operad involving insertions is related by some
assumptions to the categorial structure (see the precise
definition in Section 12).

The notion of weak Cat-operad will be to the notion of Cat-operad
what the notion of bicategory is to the notion of 2-category (see
\cite{ML98}, Section XII.6). The equations of operads like
associativity of insertion are replaced by isomorphisms in the
categorial structure, and one has to make assumptions concerning
these isomorphisms to ensure coherence. Coherence means here, as
in Mac Lane's original use of the term, that ``all diagrams
commute'', i.e.\ all diagrams of canonical arrows do so. Coherence
for weak Cat-operads is like Mac Lane's coherence for monoidal
categories (see \cite{ML63} and \cite{ML98}, Section VII.2), and like coherence
for bicategories of \cite{MLP85}.

Besides this motivation from the theory of operads, this paper may
be taken as being motivated by the theory of multicategories.
Multicategories, as conceived by Lambek in \cite{L69} and
\cite{L89}, are a generalization of the notion of operad with
insertions primitive, where not only the arity of the operadic
operations counts. Instead of one-sorted operadic operations we
have many-sorted operations (see also Section 12). This notion is
the categorial reconstruction of Gentzen's sequents of \cite{G35}
(singular sequents, with a single conclusion), and it is
interesting for proof theory.

For multicategories we have, as for operads with insertions
primitive, two kinds of associativity (the corresponding equations
analogous to (\emph{assoc}~1) and (\emph{assoc}~2) are in Lambek's
papers), and the mathematics involved in finding a weak notion of
multicategory analogous to our weak Cat-operads would not differ
essentially from what we have in this paper. That would only
involve additions not influencing in a significant way the
mathematical core. Our category \WOu\ in Part II is not far from
this notion, but in this paper we will not go further into this
matter.

Yet another motivation for the present paper would come from
matters investigated in \cite{DP06}. There one finds insertions as
operations used to present in a non-standard manner an algebra
with a free binary operation (a groupoid in the sense of universal
algebra), with or without unit. These insertions satisfy again the
two kinds of associativity. As, by replacing the equations by
isomorphisms, the two kinds of associativity are weakened in the
notion of weak Cat-operad, so this can be done in the context of
\cite{DP06}. In this weakened context we would be interested in
the \emph{constructing} of an element of our algebra rather than
in this element itself, which is related to matters treated in
\cite{DP10} and in Section 13 of this paper.

A motivation for this paper is also in the theory of polyhedra
related to permutohedra and associahedra (see below and Section
13). This matter is also interesting because, due to the convexity
of the polyhedra in question (which follows from the realizations
presented in \cite{DP11}), it would yield an alternative proof of
coherence in the style of Stasheff, as a proof of Mac Lane's
coherence for monoidal categories may be based on \cite{S63} (see
\cite{K93} and \cite{DP07}).

Up to now we have been concerned with the motivation for our
paper. We will now survey its content. Our goal is to formulate a
notion of weak Cat-operad, spell out the coherence conditions, and
demonstrate coherence. The paper is organized in three parts.

In the first part (Sections 2-6) we introduce the free unitary
operad \O\ obtained from a generating set $G$ of free operadic
operations. The general notion of operad with insertions primitive
is based on \O\ (see Section 12). Next we formulate a partial
algebraic structure we call \Oe, which is essentially a notational
variant of \O. The structure of \Oe\ is less handy for basing the
general notion of operad on it; for that \O\ is better. The
structure \Oe\ is however handier to work with than \O\ when it
comes to indexing the insertions.

Instead of indexing insertions by natural numbers, which stand for
the number of the place where the insertion is made counting from
the left, we introduce a \emph{name} for that place, which is a
word in the alphabet of natural numbers, a finite sequence of
natural numbers. While the number of the place is
context-dependent --- it may increase when further insertions are
made on its left --- the name we introduce is context-independent;
it is invariant. The arity of an operadic operation, the ordinary
\emph{numerical} arity, which is just a natural number in \O\ and
in the general notion of operad, becomes a set of our names. We
call such a set a \emph{nominal arity}. The nominal arities of
this paper may be conceived as made of leaves of trees (see the
examples in Section~13).

Invariance makes \Oe\ handier in the following sense. At the cost
of making more involved the indexing of insertions, we have
simplified the equations expressing the two kinds of
associativity.

This will become very important when these equations are replaced
by isomorphisms in the second part of the paper (Sections 7-13),
and in the remainder. Expressing the coherence conditions for the
isomorphisms using numerical arities, instead of our nominal
arities, is possible, but it would be extremely and unnecessarily
complicated (for an example see Section~12).

We introduce also in the first part a structure more general than
\Oe, which we call \Ou. The structure \Oe\ corresponds, roughly,
to just one level of \Ou, which may be understood as a
multicategory freely generated by some particular generators (see
the end of Section~5).

In \Ou\ we have instead of a family of insertions, indexed by our
names, just one partial operation of insertion. This is achieved
at the cost of, so to speak, moving the indices of insertions into
the structure in order to distinguish various occurrences of the
same element of the structure, various occurrences of the same
operadic operation.

Because of this \emph{diversification} the structure may look more
complicated. As a matter of fact, it is still handier to work
with. The advantage of having just one insertion operation proves
very big when it comes to axiomatizing the coherence conditions,
in the second part of the paper.

With \Ou\ the two kinds of associativity become, on the one hand,
just plain associativity and, on the other hand, associativity
mixed with commutativity, of a \emph{single} insertion operation.
The coherence conditions for the corresponding isomorphisms
reproduce now to a great extent the coherence conditions for
associativity and commutativity, well known from Mac Lane's
coherence results for monoidal and symmetric monoidal categories
(see \cite{ML63} and \cite{ML98}, Chapter VII), and from the standard presentation
of symmetric groups. There are however new coherence conditions
mixing the two kinds of associativity isomorphisms; let us call
them \emph{mixed} coherence conditions.

Proving that the previously known coherence conditions together
with the new mixed coherence conditions are sufficient for
coherence is possible to do with a single insertion operation, and
we do that in the third, and last, part of the paper (Sections
14-18). We suppose it would be possible to prove that also
directly in the \Oe\ version, with a family of insertions, but it
would be much more complicated, and moreover it would be
unnecessarily so. The complications would not alter the underlying
combinatorial core --- they would just obscure~it.

In the second part of the paper we axiomatize our coherence
conditions in the \Ou\ and \Oe\ versions, which results in the
categories \WOu\ and \WOe. We establish that \WOu\ is the disjoint
union of isomorphic copies of \WOe, and we use \WOe\ to give our
notion of weak Cat-operad in Section 12.

The last section of that part, Section 13, is optional, and a
reader not interested in the matters treated there may skip it,
and continue reading the remainder. Through this section our paper
is however connected with important and interesting matters of algebraic
topology and combinatorics. It provides an insight into the regularity
underlying the equations with which we define weak Cat-operads. In Section 13 we consider
how the mixed coherence conditions engender a new kind of
polyhedron, related to, but different from, the three-dimensional
associahedron and permutohedron. The family to which this
polyhedron, called hemiassociahedron, belongs (the family includes
also the three-dimensional cyclohedron) is investigated in
\cite{DP10} and \cite{DP11} (which should be consulted for further
references). The mixed conditions may also however lead in some
cases to the three-dimensional associahedron and permutohedron,
which provides another perspective on these well-known polyhedra
tied to associativity and commutativity. The hemiassociahedron
arises in the non-unitary structure of our weak Cat-operads, as
the associahedron and permutohedron arise in the non-unitary
structure of monoidal categories and symmetric groups.

In the third part of the paper we prove the sufficiency of our
conditions for coherence (these conditions are of course also
necessary). We do that with \WOu, which, as we said, is easier to
work with. We \emph{strictify} first the monoidal structure of
\WOu\ (all the arrows of that structure become identity arrows),
which leaves a category \WOuth, equivalent with \WOu, which is
similar to a symmetric strictly monoidal category.

To prove that \WOuth\ is a preorder (i.e.\ that there is not more
than one arrow in it with a given source and target), which is
what coherence here amounts to, we reconsider the standard
presentation of symmetric groups, and a proof that this
presentation is complete. This proof, which involves a reduction
to a normal form, implicit already in Moore's paper \cite{M96}
(the first paper to deal with the matter), yields, as a matter of
fact, something more, pertaining to a wider family of structures.

In this family we have categories we call \Cg, which generalize
symmetric strictly monoidal categories; symmetric groups belong to
the family as one-object categories. All one has in the categories
\Cg\ comes from the symmetric structure, but this structure may be
incomplete; it may, roughly speaking, have gaps. It need not even
give a groupoid in the categorial sense (its arrows need not be
invertible, they need not be isomorphisms). The category \WOuth\
is just a particular \Cg, and that it is a preorder follows from a
coherence result for an arbitrary \Cg. This coherence involves
graphs corresponding to permutations, but in connection with
\WOuth, because of diversification (see above), the graphs need
not be mentioned. It follows that \WOu\ is a preorder, which
implies that \WOe\ is a preorder, and with that our notion of weak
Cat-operad is justified. (For other examples that fall within the
range of this general coherence result one may consult
\cite{DP06a} and \cite{DP12a}.)

Throughout the paper we distinguish matters pertaining to unitary
operads from those pertaining to non-unitary operads. The
non-unitary versions of our operads, categories and related
structures have a superscript $^-$ in their names. Whatever we
established for the unitary versions can be established for the
non-unitary ones, which make the non-unitary core of our notions.
As a matter of fact, the more interesting mathematics pertains to
the non-unitary notions (as Section 13 illustrates).

We do not claim that our notion of weak Cat-operad is the right
notion of weak Cat-operad in general, but that it is the right
notion if we are motivated by coherence involving insertions. In
our definition, this coherence has been combinatorially analyzed
by our axiomatic equations. We prove that no equation is missing.

The notion of monoidal category was introduced in a non-axiomatic
way via coherence by B\' enabou in \cite{Ben63}, and in the
equational axiomatic way, such as we favour, by Mac Lane in
\cite{ML63}. Mac Lane's definition is the standard one, while B\'
enabou's is rarely mentioned. For B\' enabou, coherence is built
into the definition, and for Mac Lane it is a theorem. One could
analogously define the theorems of classical propositional logic
as being the tautologies (this is done, for example, in
\cite{CK73}, Sections 1.2-3), in which case completeness would not
be a theorem, but would be built into the definition. Analogously
to what was done for monoidal categories in \cite{ML63}, and later
for bicategories in \cite{MLP85},  we are not only proposing a
definition, but we are proving a coherence theorem concerning~it.

An anonymous referee of our text asks whether our notion of weak
Cat-operad is equivalent with the notion of pseudo-operad of
\cite{DS01}. Another anonymous referee is convinced that the right
notion of weak Cat-operad is given by, as suggested to him by Mark
Weber, ``pseudo-algebras of a monad for which the strict algebras
are the existing notion of Cat-operad'' (without going into
details), and that referee wonders  whether this notion is the
``same'' as a ours.

To start answering these questions one would have to consider
first the difference in language, since the other notions are not
based like ours on insertions. If the other notions could be based
on insertions (in the non-unitary case defining insertions would
pose a problem), the question would reduce to the question whether
these other notions of weak Cat-operad require coherence (in our and Mac Lane's sense) involving
insertions, which seems possible. If they require this
coherence, then our paper provides a positive answer to the
question. It would show that notions, which are not like Mac
Lane's, are equivalent with a notion in the style of Mac Lane.

The study of this matter --- whether the alternative notions could
be based on insertions, and whether they require coherence
involving insertions
--- would be a study of these alternative notions, conducted independently
of the content of this paper. This is a different
topic, best left for a separate treatment, because, as far as we can see,
to be presented with sufficient detail this would require a lot of
space.

\vspace{4ex}

{\samepage \noindent{\large {\sc Part I}}

\section{The operad \O}}\label{sec2}
Let $G$ be a set, for whose members we use $x,y,z,\ldots,$ perhaps
with indices, and let ${\alpha_G\!:G\str{\mathbf N}}$ be a
function; $\alpha_G(x)$ is intuitively the arity of $x$. The
elements of $G$ are the free, generating, \emph{primitive operadic
operations} of the operad \O, which is the unitary operad freely
generated out of $G$ and $\alpha_G$ in the following manner. (As a
matter of fact, mentioning ${\alpha_G\!:G\str{\mathbf N}}$ is
enough; it carries the information about $G$.)

We define first inductively the set of \emph{terms} of \O; these
terms will stand for the \emph{operadic operations} of \O.
Together with the terms of \O\ we define simultaneously a function
$\alpha$ from the terms of \O\ to $\mathbf N$; the values of this
function are intuitively the arity. Here are the three clauses of
these two simultaneous inductive definitions:
\begin{tabbing}
\hspace{0em}\=(0) \hspace{.5em} \=if $x\in G$, then $x$ is a
term;
\hspace{.2em} $\alpha(x)=\alpha_G(x)$;\\[.5ex]
\>(1)\>$\Iota$ is a term; \hspace{.2em} $\alpha(\Iota)=1$;\\[.5ex]
\>(2)\>if $\varphi$ and $\gamma$ are terms and $1\leq n\leq
\alpha(\gamma)$, then $\gamma\ins_n \varphi$ is a term;\\*
\>\>$\alpha(\gamma\ins_n \varphi)=\alpha(\gamma)\mn 1\pl
\alpha(\varphi)$.
\end{tabbing}
Officially, in (2) we should have $(\gamma\ins_n \varphi)$ instead
of $\gamma\ins_n \varphi$, but, as usual, the outermost
parentheses of these and other terms later will be taken for
granted, and omitted. We use $\varphi,\gamma,\chi,\ldots,$ perhaps
with indices, for the terms of \O.

The term $\Iota$ stands for the \emph{unit} operadic operation,
and $\ins_k$ stands for a partial operation of \emph{insertion}.
The expression $\gamma\ins_n \varphi$ does not become a term for
every substitution for $n$, $\varphi$ and $\gamma$; substitutions
for $n$ and $\gamma$ are restricted. We express this by saying
that $\gamma\ins_n \varphi$ is \emph{legitimate} when $1\leq n\leq
\alpha(\gamma)$.

Analogously, substitutions will be restricted for the equations
between terms of the operad \O. As for terms, we express this by
saying that $\varphi=\gamma$ is legitimate when
$\alpha(\varphi)=\alpha(\gamma)$. An equation cannot hold between
terms of different arity.

The equations of \O\ between terms are given through an axiomatic
equational system, and the operadic operations of \O\ will be
formally equivalence classes of terms of \O\ such that these
equations are satisfied. Besides $\varphi=\varphi$, the
\emph{axiomatic equations} of \O\ are the following:
\begin{tabbing}
\hspace{0em}\=(\emph{unit})\hspace{2em}\=$\varphi\ins_n
\Iota=\varphi$, \hspace{1em} $\Iota\ins_1\varphi=\varphi$,
\\[1ex]
\>(\emph{assoc}~1)\>$(\chi\ins_n\gamma)\ins_m\varphi=\chi\ins_n(\gamma\ins_{m-n+1}\varphi)$,
\hspace{.5em} provided $n\leq m<n\pl\alpha(\gamma)$,
\\[1ex]
\>(\emph{assoc}~2)\>$(\chi\ins_n\gamma)\ins_m\varphi=
(\chi\ins_{m-\alpha(\gamma)+1}\varphi)\ins_n\gamma$, \hspace{.5em}
provided $n\pl\alpha(\gamma)\leq m$.
\end{tabbing}

The proviso for (\emph{assoc}~1) may be derived from the
legitimacy of $\gamma\ins_{m-n+1}\varphi$, for which we must have
$1\leq m\mn n\pl 1\leq\alpha(\gamma)$. So, as a matter of fact,
this proviso need not have been mentioned, and in the future we
will not always mention such provisos, which may be inferred from
the legitimacy of an equation, or of our notation for terms. We
assume that all the expressions for terms that occur in an
equation are legitimate.

The proviso for (\emph{assoc}~2) is not derivable in this manner.
This equation could be replaced by
\begin{tabbing}
\hspace{4em}$(\chi\ins_n\gamma)\ins_m\varphi=
(\chi\ins_m\varphi)\ins_{n+\alpha(\varphi)-1}\gamma$,
\hspace{.5em} provided $m<n$.
\end{tabbing}

The remaining equations of \O\ are derived with the help of the
rules of symmetry and transitivity of $=$, and of the rule of
$\ins_n$-\emph{congruence}:
\begin{tabbing}
\hspace{4em}from $\varphi_1=\varphi_2$ and $\gamma_1=\gamma_2$
derive $\gamma_1\ins_n\varphi_1=\gamma_2\ins_n\varphi_2$,
\end{tabbing}
provided both sides of the last equation are legitimate. As a
matter of fact, it is superfluous to state this proviso; we
understand rules like $\ins_n$-congruence always with such
provisos.

Once we have defined this axiomatic equational system, it can be
verified by an easy induction on the length of derivation that for
every equation $\varphi=\gamma$ of \O\ we have
$\alpha(\varphi)=\alpha(\gamma)$. So these equations are indeed
legitimate.

This concludes our definition of the operad \O. This operad may be
conceived as a partial algebra, i.e.\ algebra with partial
operations, $\langle C(\mathcal{O}),\{\ins_n\mid {n\in{\mathbf
N}^+}\},\Iota\rangle$, with $C$(\O), the \emph{carrier} of \O,
being the set of operadic operations of \O.

The non-unitary operad \Om\ freely generated out of $G$ and
$\alpha_G$ is defined like \O\ save that we omit clause (1) from
the definition of terms, and we omit (\emph{unit}) from the
axiomatic equations.

\section{The structure \Oe}\label{sec3}
Before we introduce \Oe\ we deal with preliminary matters
concerning nominal arities.

An \Np-\emph{word} is a finite (possibly empty) sequence of
natural numbers greater than $0$. We use $a,b,c,\ldots,$ perhaps
with indices, for \Np-words; we reserve $e$ to denote the empty
\Np-word.

An \Np-\emph{language} is a set of \Np-words. An \Np-language is
called a \emph{nominal arity} when there are no two distinct \Np-words of
the form $a$ and $ab$ in it; i.e.\ such that one is a proper initial
segment of the other. This definition
allows also infinite nominal arities, but in this paper we have
use only for finite ones. The empty \Np-language and every
singleton \Np-language are nominal arities. We use $X,Y,Z,\ldots,$
perhaps with indices, for nominal arities.

We say that the \Np-word $a$ is a \emph{prefix} of the nominal
arity $X$ when for every $c$ in $X$ we have that $c$ is of the
form $ab$; i.e., $a$ is an initial segment, not necessarily
proper, of every member of $X$. Note that $a$ may be $e$, which is
a prefix of every nominal arity. Note also that $X$ may have more
than one prefix, of which $e$ is always one. It is trivially
satisfied that every $a$ is a prefix of the empty nominal arity.
The set of prefixes of $X$ is denoted by $P_X$.

For every \Np-language $M$, and, in particular, for $M$ a nominal
arity, let
\begin{tabbing}
\hspace{4em}$a\cdot M=_{df}\{ab\mid b\in M\}$.
\end{tabbing}
We have, of course, $e\cdot M=M$.

Every nominal arity is linearly ordered by the lexicographical
order $\prec$, whose definition for nominal arities is simpler,
and is given by the following:
\begin{tabbing}
\hspace{4em}$a_1\prec a_2$ \hspace{.1em} iff \hspace{.1em}
$\exists a,b,c(\exists n,m\in {\mathbf N}^+)(a_1=anb \;\;\&\;\;
a_2=amc \;\;\&\;\; n<m)$.
\end{tabbing}

If $a\in P_X\cap Y$, then we define the result of \emph{inserting}
$X$ in $Y$ at $a$:
\begin{tabbing}
\hspace{4em}$Y\ins_a X=_{df}(Y\mn\{a\})\cup X$.
\end{tabbing}
Note that this union is disjoint. Otherwise, for some $b$ in $Y$
different from $a$ we would have that it is in $X$, and hence that
it has $a$ as an initial segment; this contradicts the assumption
that $Y$ is a nominal arity. We conclude in a similar manner, by
going through all possible cases, that $Y\ins_a X$ is a nominal
arity. The expression $Y\ins_a X$ is legitimate when $a\in P_X\cap
Y$.

For $|M|$ being the cardinality of the set $M$, we have the
following.

\begin{rem}\label{rem3.1}
$|a\cdot X|=|X|$, \quad $|Y\ins_a X|=|Y|\mn 1\pl |X|$.
\end{rem}

\noindent For the second equation we rely on the disjointness
mentioned after the definition of $\ins_a$. We also have the
following.

\begin{rem}\label{rem3.2}
$a\cdot P_X\subseteq P_{a\cdot X}$,
\quad $P_Y\subseteq P_{Y\ins_a X}$.
\end{rem}

\noindent The inclusion converse to the first one holds only if
$a$ is $e$; if $a$ is not $e$, then $e\in P_{a\cdot X}$ and
$e\notin a\cdot P_X$. The inclusion converse to the second one
does not hold for $Y=\{a\}$ and $X=\emptyset$. We also have the
following two remarks.

\begin{rem}\label{rem3.3}
$a\cdot(b\cdot X)=ab\cdot X$,\quad $a\cdot(Y\ins_b X)=a\cdot
Y\ins_{ab} a\cdot X$.
\end{rem}
\begin{rem}\label{rem3.4}
For $a\in Y$ and $b\in Z$, $(Z\ins_b Y)\ins_a X=Z\ins_b(Y\ins_a
X)$;\\*[.5ex] for $a,b\in Z$, $(Z\ins_b Y)\ins_a X=(Z\ins_a
X)\ins_b Y.$
\end{rem}
The condition $b\in Z$ is implied by the legitimacy of $Z\ins_b Y$
on the left-hand sides of both equations of the last remark. The
condition $a\in Y$ is implied by the legitimacy of $Y\ins_a X$ on
the right-hand side of the first equation, and the condition $a\in
Z$ by the legitimacy of $Z\ins_a X$ on the right-hand side of the
second equation. In the first equation we must also have that $a$
is of the form $bc$, since $b\in P_Y$. In the second equation this
is excluded, since $Z$ is a nominal arity.

Let $\bar{n}=_{df}\{1,\dots,n\}$ for $n\geq 1$, and let
$\bar{0}=\emptyset$. It is clear that $\bar{n}$ is a nominal arity
for every $n\geq 0$, and $P_{\bar{n}}=\{e\}$ for $n>0$.

We pass now to the definition of the structure \Oe. Let $G$ and
$\alpha_G$ be as for \O\ (see the beginning of the preceding
section). We define first inductively with three clauses the set
of \emph{terms} of \Oe\ together with a function $s$ from these
terms to nominal arities:
\begin{tabbing}
\hspace{0em}\=(0$_e$) \hspace{.5em} \=if $x\in G$, then $x$ is a
term; \hspace{.2em} $s(x)= \left\{
\begin{array}{l}
\overline{\alpha_G(x)}, \hspace{1em} {\mbox{\rm if }}
\alpha_G(x)\neq 1,
\\[.5ex]
\{e\}, \hspace{2.2em} {\mbox{\rm if }} \alpha_G(x)=1;
\end{array}
\right .$\\
\>(1$_e$)\>$\Iota$ is a term; \hspace{.2em} $s(\Iota)=\{e\}$;\\[1.7ex]
\>(2$_e$)\>if $f$ and $g$ are terms and $a\in s(g)$, then $g\ins_a
f$ is a term;\\* \>\>$s(g\ins_a f)=s(g)\ins_a a\cdot s(f)$.
\end{tabbing}
(Why we do not have $s(x)= \overline{\alpha_G(x)}= \{1\}$ if
$\alpha_G(x)=1$ is explained after the definition of the terms of
\Ou\ in Section~5.) We use $f,g,h,\ldots,$ perhaps with indices,
for the terms of \Oe.

Here, as for \O, the term $\Iota$ stands for the unit operadic
operation, while $\ins_a$ stands for a partial operation of
insertion. The expression $g\ins_a f$ is legitimate when $a\in
s(g)$, and $f=g$ is legitimate when $s(f)=s(g)$.

The equations of \Oe\ between terms are given as for \O\ through
an axiomatic equational system, whose \emph{axiomatic equations}
besides $f=f$ are the following equations:
\begin{tabbing}
\hspace{0em}\=(\emph{unit}$_e$)\hspace{2em}\=$f\ins_a \Iota=f$,
\hspace{1em} $\Iota\ins_e f=f$,
\\[1ex]
\>(\emph{assoc}~1$_e$)\>$(h\ins_b g)\ins_{ba} f= h\ins_b(g\ins_a
f)$,
\\[1ex]
\>(\emph{assoc}~2$_e$)\>$(h\ins_b g)\ins_a f= (h\ins_a f)\ins_b
g$.
\end{tabbing}
As rules we have symmetry and transitivity of $=$, as for \O, and
$\ins_a$-\emph{congruence}:
\begin{tabbing}
\hspace{4em}from $f_1=f_2$ and $g_1=g_2$ derive $g_1\ins_a
f_1=g_2\ins_a f_2$.
\end{tabbing}

This concludes our definition of the equations of \Oe. As for \O,
the operadic operations of \Oe\ are equivalence classes of terms
such that the equations of \Oe\ are satisfied. The structure \Oe\
should be conceived as a partial algebra $\langle
C(\mathcal{O}_e),\{\ins_a\mid a\; \mbox{\rm is an
\Np-word}\},\Iota\rangle$, with the carrier $C$(\Oe) of \Oe\ being
the set of operadic operations of \Oe.

The legitimacy of $h\ins_b g$ on the left-hand sides of
(\emph{assoc}~1$_e$) and (\emph{assoc}~2$_e$) implies $b\in s(h)$;
the legitimacy of $g\ins_a f$ on the right-hand side of
(\emph{assoc}~1$_e$) implies $a\in s(g)$, and the legitimacy of
$h\ins_a f$ on the right-hand side of (\emph{assoc}~2$_e$) implies
$a\in s(h)$.

When the indices of (\emph{assoc}~1$_e$) are compared with those
of the first equation of Remark~\ref{rem3.4}, one should bear in
mind that we have the following:
\begin{tabbing}
\hspace{4em}\=$s((h\ins_b g)\ins_{ba} f)\;$\=$= (s(h)\ins_b
b\cdot s(g))\ins_{ba} ba\cdot s(f)$, \hspace{.5em}\=by definition,
\\[.5ex]
\>\>$= s(h)\ins_b (b\cdot s(g)\ins_{ba} ba\cdot s(f))$,\>by
Remark~\ref{rem3.4},
\\[.5ex]
\>\>$=s(h\ins_b(g\ins_a f))$,\>by definition.
\end{tabbing}
Hence we have verified that for $f_1=f_2$ being an instance of
(\emph{assoc}~1$_e$) we have $s(f_1)=s(f_2)$. This is verified
analogously for the other axiomatic equations of \Oe, and this
makes the basis of the induction on the length of derivation that
shows that for every equation $f=g$ of \Oe\ we have $s(f)=s(g)$.

We define the non-unitary structure \Oem\ like \Oe\ save that we
omit clause (1$_e$) from the definition of terms, and we omit
(\emph{unit}$_e$) from the axiomatic equations.

\section{\O\ and \Oe}\label{sec4}
In this section we show that \O\ and \Oe\ can be mapped one to the
other in such a manner that \Oe\ may be considered just an
alternative notation for \O.

For $X$ a finite nominal arity, let the bijection ${K_X\!:X\str
\overline{|X|}}$ be defined by
\begin{tabbing}
\hspace{4em}$K_X(a)=|\{b\in X\mid b\prec a\}|\pl 1.$
\end{tabbing}
Note that $|\bar{n}|=n$, so that $K_{\bar{n}}$ is a function from
$\bar{n}$ to $\bar{n}$.

It is easy to see that $a_1\prec a_2$ implies $K_X(a_1)<K_X(a_2)$,
provided $a_1,a_2\in X$. It is also easy to verify the following:
\begin{tabbing}
\hspace{0em}\=($K1$)\hspace{1.5em}\=$K_{\bar{n}}$ is the
identity function from $\bar{n}$ to $\bar{n}$,
\hspace{.2em} $K_{\{e\}}(e)=1$,\\[1.2ex]
\>($K2$)\>$K_{b\cdot X}(ba)=K_X(a)$,\\[1.5ex]
\>($K3.1$)\>$K_{Y\ins_b X}(a)=K_Y(b)\mn 1\pl
K_X(a)$, \hspace{.2em} if $a\in X$,\\[1ex]
\>($K3.2$)\>$K_{Y\ins_b X}(a)=K_Y(a)$, \hspace{.2em} if $a\in Y$
and $a\prec b$,\\[1ex]
\>($K3.3$)\>$K_{Y\ins_b X}(a)=K_Y(a)\mn 1\pl |X|$, \hspace{.2em}
if $a\in Y$ and $b\prec a$.
\end{tabbing}

Let ${K^{-1}_X\!:\overline{|X|}\str X}$ be the bijection inverse
to $K_X$. It is easy to verify for $K^{-1}_X$ the following
equations, which are interdeducible with those we have just given
for $K_X$:
\begin{tabbing}
\hspace{0em}\=($K^{-1}1$)\hspace{1.5em}\=$K^{-1}_{\bar{n}}=K_{\bar{n}}$,
\hspace{.2em} $K^{-1}_{\{e\}}(1)=e$,\\[1.2ex]
\>($K^{-1}2$)\>$K^{-1}_{b\cdot X}(n)=bK^{-1}_X(n)$, \hspace{.2em}
for $n\in \overline{|X|}$,\\[1.5ex]
\>($K^{-1}3.1$)\>$K^{-1}_{Y\ins_b X}(m)=
K^{-1}_X(m\mn n\pl 1)$, if $K_Y(b)\!=n$ and $n\leq m<n\pl |X|$,\\[1ex]
\>($K^{-1}3.2$)\>$K^{-1}_{Y\ins_b X}(m)=K^{-1}_Y(m)$,
\hspace{.2em} if $m<K_Y(b)$,\\[1ex]
\>($K^{-1}3.3$)\>$K^{-1}_{Y\ins_b X}(m)=K^{-1}_Y(m\mn |X|\pl 1)$,
\hspace{.2em} if $K_Y(b)\pl |X|\leq m$.
\end{tabbing}

These equations and those given above for $K_X$ serve to explain
how we pass from the indices in the axiomatic equations of \O\ to
those in the axiomatic equations of \Oe, and vice versa. In
particular the condition of ($K^{-1}3.1$) is transferred to
(\emph{assoc}~1), while the conditions of ($K^{-1}3.2$) and
($K^{-1}3.3$) are transferred to (\emph{assoc}~2).

Next we define inductively a map $\varepsilon$ from the terms of
\O\ to the terms of \Oe, for which we will show below (in
Propositions \ref{prop4.1} and \ref{prop4.2}) that it is a
bijection:
\begin{tabbing}
\hspace{4em}\=$\varepsilon(x)=x$, \hspace{.5em}
$\varepsilon(\Iota)=\Iota$,\\[.5ex]
\>$\varepsilon(\gamma\ins_n\varphi)= \varepsilon(\gamma)\ins_a
\varepsilon(\phi)$, \hspace{.2em} for
$K^{-1}_{s(\varepsilon(\gamma))}(n)=a$.
\end{tabbing}
We define inductively the map $\tau$ from the terms of \Oe\ to the
terms of \O, for which we will show that it is the inverse of
$\varepsilon$:
\begin{tabbing}
\hspace{4em}\=$\tau(x)=x$, \hspace{.5em}
$\tau(\Iota)=\Iota$,\\[.5ex]
\>$\tau(g\ins_a f)= \tau(g)\ins_n \tau(f)$, \hspace{.2em} for
$K_{s(g)}(a)=n$.
\end{tabbing}

Then we can establish the following propositions by
straightforward inductions on the complexity of $\varphi$ and~$f$.

\begin{prop}\label{prop4.1}
For every term $\varphi$ of \O\ we have that
$\tau(\varepsilon(\varphi))$ is~$\varphi$.
\end{prop}

\begin{prop}\label{prop4.2}
For every term $f$ of \Oe\ we have that
$\varepsilon(\tau(f))$ is~$f$.
\end{prop}

By inductions on the complexity of $\varphi$ and $f$ we can also
straightforwardly establish the following lemmata.

\begin{lem}\label{lem4.3}
For every term $\varphi$ of \O\ we have
$|s(\varepsilon(\varphi))|=\alpha(\varphi)$.
\end{lem}

\begin{lem}\label{lem4.4}
For every term $f$ of \Oe\ we have
$\alpha(\tau(f))=|s(f)|$.
\end{lem}

\noindent These lemmata are used in the proof of the following two
propositions, which consist in inductions on the length of
derivation. The main part of these inductions is however in the
basis, when we deal with axiomatic equations.

\begin{prop}\label{prop4.5}
If $\varphi=\gamma$ in \O, then
$\varepsilon(\varphi)=\varepsilon(\gamma)$ in \Oe.
\end{prop}

\begin{prop}\label{prop4.6}
If $f=g$ in \Oe, then $\tau(f)=\tau(g)$ in
\O.
\end{prop}

\noindent So $\varepsilon$ and $\tau$ induce bijections inverse to
each other between the operadic operations of \O\ and \Oe. This
shows that \Oe\ is just a notational variant of \O.

The bijections between the terms of \Om\ and \Oem\ are obtained by
just restricting $\varepsilon$ and $\tau$, and then for these
bijections we can establish as well as Propositions \ref{prop4.1}
and \ref{prop4.2} the analogues of Propositions \ref{prop4.5} and
\ref{prop4.6}, where \O\ and \Oe\ are replaced respectively by
\Om\ and \Oem. So the relationship between \Om\ and \Oem\ is
exactly analogous to that between \O\ and~\Oe.

\section{The structure \Ou}\label{sec5}
We introduce now the structure \Ou, which generalizes \Oe.

Let $G$ and $\alpha_G$ be as for \O\ and \Oe\ (see the beginning
of Section~2). We define first inductively with three clauses the
set of \emph{terms} of \Ou\ together with a function $s$ from
these terms to nominal arities (this function is related to the
function $s$ from the terms of \Oe, and this is why it bears the
same name) and a function $t$ from these terms to \Np-words:
\begin{tabbing}
\hspace{0em}\=(0$_u$) \hspace{.5em} \=if $x\in G$ and $a$ is an
\Np-word, then $a\cdot x$ is a term;\\*[.5ex]
\>\>$s(a\cdot x)=
\left\{
\begin{array}{l}
a\cdot \overline{\alpha_G(x)}, \hspace{.5em} {\mbox{\rm if }}
\alpha_G(x)\neq 1,
\\[.5ex]
\{a\}, \hspace{2.9em} {\mbox{\rm if }} \alpha_G(x)=1,
\end{array}
\right .$ \hspace{.5em} $t(a\cdot x)=a$;\\[1.5ex]
\>(1$_u$)\>if $a$ is an \Np-word, then $a\cdot \Iota$ is a term;
\hspace{.2em} $s(a\cdot \Iota)=\{a\}$,
\hspace{.2em} $t(a\cdot \Iota)=a$;\\[1.5ex]
\>(2$_u$)\>if $f$ and $g$ are terms and $t(f)\in s(g)$, then
$g\ins f$ is a term;\\* \>\>$s(g\ins f)=s(g)\ins_{t(f)} s(f)$,
\hspace{.2em} $t(g\ins f)=t(g)$.
\end{tabbing}
If $\alpha_G(x)=1$, then we can envisage having $s(a\cdot x)=
a\cdot \overline{\alpha_G(x)}= \{a1\}$ (which would entail $s(x)=
\overline{\alpha_G(x)}= \{1\}$ for \Oe). This way diversification
(see Section~1) would apply also to the unary members of $\Gamma$,
but a difference would arise with the treatment of $\Iota$, which
must have the clause above. We have preferred however not to
distinguish these two unary cases, because this is not essential.
The diversification we achieve is sufficient for our purposes (in
particular for Section~18).

In order to verify for (2$_u$) that $s(g)\ins_{t(f)} s(f)$ is
legitimate if $t(f)\in s(g)$, we have first that $t(a\cdot x)\in
P_{s(a\cdot x)}$ and $t(a\cdot \Iota)\in P_{s(a\cdot \Iota)}$, and
we have the following two remarks.

\begin{rem}\label{rem5.1}
If $t(f)\in P_{s(f)}$ and $t(f)\in s(g)$, then $s(g)\ins_{t(f)}
s(f)$ is legitimate.
\end{rem}

\begin{rem}\label{rem5.2}
If $t(g)\in P_{s(g)}$, then $t(g\ins f)\in
P_{s(g\ins f)}$.
\end{rem}

\noindent To justify the last remark we have
\begin{tabbing}
\hspace{4em}$t(g\ins f)=t(g)\in P_{s(g)}\subseteq
P_{s(g)\ins_{t(f)} s(f)}$
\end{tabbing}
by Remark~\ref{rem3.2}. These two remarks, together with what we
said before them for the basis of the induction, yield that
$t(f)\in P_{s(f)}$ holds for every term $f$ of \Ou.

This concludes our definition of the terms of \Ou. We use now
$f,g,h,\ldots,$ perhaps with indices, for the terms of \Ou. The
expression $g\ins f$ is legitimate when $t(f)\in s(g)$, and $f=g$
is legitimate when $s(f)=s(g)$ and $t(f)=t(g)$.

The equations of \Ou\ between terms are given as before through an
axiomatic equational system, whose \emph{axiomatic equations}
besides $f=f$ are the following equations:
\begin{tabbing}
\hspace{0em}\=(\emph{unit}$_u$)\hspace{2em}\=$f\ins a\cdot
\Iota=f$, \hspace{1em} $t(f)\cdot \Iota\ins f=f$,
\\[1ex]
\>(\emph{assoc}~1$_u$)\>$(h\ins g)\ins f= h\ins(g\ins f)$,
\\[1ex]
\>(\emph{assoc}~2$_u$)\>$(h\ins g)\ins f= (h\ins f)\ins g$.
\end{tabbing}
As rules we have symmetry and transitivity of $=$, and
$\ins$-congruence, which is like $\ins_a$-congruence of \Oe\ with
the subscript $a$ omitted.

This concludes our definition of the equations of \Ou, and of the
operadic operations of \Ou\ (which, as before, are equivalence
classes of terms). The structure \Ou\ should be conceived as a
partial algebra $\langle C(\mathcal{O}_u),\ins,\{a\cdot\Iota\mid
a\; \mbox{\rm is an \Np-word}\}\rangle$, with the carrier $C$(\Ou)
of \Ou\ being the set of operadic operations of \Ou.

The legitimacy of $h\ins g$ on the left-hand sides of
(\emph{assoc}~1$_u$) and (\emph{assoc}~2$_u$) implies $t(g)\in
s(h)$, while the legitimacy of $g\ins f$ on the right-hand side of
(\emph{assoc}~1$_u$) implies $t(f)\in s(g)$, and the legitimacy of
$h\ins f$ on the right-hand side of (\emph{assoc}~2$_u$) implies
$t(f)\in s(h)$. By induction on the length of derivation we
establish that for every equation $f=g$ of \Ou\ we have
$s(f)=s(g)$ and $t(f)=t(g)$.

We define the non-unitary structure \Oum\ like \Ou\ save that we
omit clause (1$_u$) from the definition of terms, and we omit
(\emph{unit}$_u$) from the axiomatic equations.

The structure \Ou\ amounts to the free multicategory generated by
the multigraph made of the objects in $\{a\mid a\; \mbox{\rm is an
\Np-word}\}$ and the multiarrows in $\{a\cdot x\mid x\in G$ and
$a$ is an \Np-word$\}$, with the source and target functions given
by the functions $s$ and $t$ (see \cite{L89}, Section~3).

\section{\Oe\ and \Ou}\label{sec6}
In this section we establish the correspondences that exist
between \Oe\ and some structures derived form \Ou.

For every \Np-word $a$ we define inductively a map $a\cdot$ from
the terms of \Ou\ to the terms of \Ou:
\begin{tabbing}
\hspace{4em}\=$a\cdot(b\cdot x)=ab\cdot x$, \hspace{.5em}
$a\cdot(b\cdot\Iota)=ab\cdot\Iota$,\\[.5ex]
\>$a\cdot(g\ins f)=a\cdot g\ins a\cdot f$;
\end{tabbing}
$a\cdot f$ stands for $a\cdot(f)$, and we read $a\cdot g\ins
a\cdot f$ as $(a\cdot g)\ins (a\cdot f)$. In order to verify that
$a\cdot g\ins a\cdot f$ is legitimate, we need to show that
$t(a\cdot f)\in s(a\cdot g)$. For that we establish first that we
have
\begin{tabbing}
\hspace{4em}($st\;a\cdot$)\hspace{1.5em}$s(a\cdot h)=a\cdot
s(h)$, \hspace{.2em} $t(a\cdot h)=at(h)$
\end{tabbing}
for $h$ being $b\cdot x$ and $b\cdot\Iota$ (we use here the first
equation of Remark~\ref{rem3.3}). Next we have the following two
remarks.

\begin{rem}\label{rem6.1}
If $t(a\cdot f)=at(f)$, $s(a\cdot g)=a\cdot s(g)$ and $t(f)\in
s(g)$, then $t(a\cdot f)\in s(a\cdot g)$.
\end{rem}

\begin{rem}\label{rem6.2}
If {\rm ($st\;a\cdot$)} holds for $h$ being $f$
and $g$, then it holds for $h$ being $g\ins f$.
\end{rem}

\noindent The last remark (for which we use the second equation of
Remark~\ref{rem3.3}) is the induction step, which together with
what we said above for the basis of the induction, yields that
($st\;a\cdot$) holds for every term $h$ of \Ou.

It is easy to infer that for every term $f$ of \Ou\ we have
$e\cdot f=f$ and $a\cdot(b\cdot f)=ab\cdot f$. We can also easily
establish the following by induction on the length of derivation.

\begin{rem}\label{rem6.3}
If $f=g$ in \Ou, then $a\cdot f=a\cdot g$ in
\Ou.
\end{rem}

\noindent So $a\cdot$ induces a map from the operadic operations
of \Ou\ to the operadic operations of \Ou.

We need $a\cdot$ to define inductively a map $U$ from the terms of
\Oe\ to the terms of \Ou:
\begin{tabbing}
\hspace{4em}\=$U(x)=e\cdot x$, \hspace{.5em}
$U(\Iota)=e\cdot\Iota$,\\[.5ex]
\>$U(g\ins_a f)=U(g)\ins a\cdot U(f)$.
\end{tabbing}
In order to verify that $U(g)\ins a\cdot U(f)$ is legitimate if
$g\ins_a f$ is legitimate, we need to show that $t(a\cdot U(f))\in
s(U(g))$ follows from $a\in s(g)$. For that we rely on
($st\;a\cdot$) and on
\begin{tabbing}
\hspace{4em}($st\;U$)\hspace{1.5em}$s(U(h))=s(h)$, \hspace{.2em}
$t(U(h))=e$,
\end{tabbing}
where $h$ is $g$ or $f$, as needed. It is clear that ($st\;U$)
holds for $h$ being $x$ and $\Iota$, and then we may establish the
induction step, which yields that ($st\;U$) holds for every term
$h$ of~\Oe.

If $c$ is the \Np-word $ab$, then $a\bsl c$ is defined, and is
$b$; i.e.\ $a\bsl ab=b$. Since $c$ is $ec$, we have that $e\bsl c$
is always defined, and is~$c$.

Then we define inductively a map $E$ from the terms of \Ou\ to the
terms of~\Oe:
\begin{tabbing}
\hspace{4em}\=$E(a\cdot x)=x$, \hspace{.5em}
$E(a\cdot\Iota)=\Iota$,\\[.5ex]
\>$E(g\ins f)=E(g)\ins_{t(g)\bsl t(f)} E(f)$.
\end{tabbing}
We have that $t(f)$ in the last line is of the form $t(g)b$
because $t(g)\in P_{s(g)}$ and $t(f)\in s(g)$. In order to verify
that $t(g)\bsl t(f)\in s(E(g))$, which we need for the legitimacy
of $E(g)\ins_{t(g)\bsl t(f)} E(f)$, we rely on
\begin{tabbing}
\hspace{4em}($st\;E$)\hspace{1.5em}$t(h)\cdot s(E(h))=s(h)$,
\end{tabbing}
where $h$ is $g$. It is clear that ($st\;E$) holds for $h$ being
$a\cdot x$ and $a\cdot\Iota$, and then we may establish the
induction step, which yields that ($st\;E$) holds for every term
$h$ of \Ou.

With the help of ($st\;a\cdot$), it is easy to establish the
following by induction on the complexity of $f$.

\begin{lem}\label{lem6.4}
For every term $f$ of \Ou\ we have that $E(a\cdot f)$ is $E(f)$.
\end{lem}

\noindent This lemma, together with ($st\;a\cdot$) and ($st\;U$),
serves for the first of the following two propositions, which are
proved by inductions on the complexity of $f$.

\begin{prop}\label{prop6.5}
For every term $f$ of \Oe\ we have that $E(U(f))$ is $f$.
\end{prop}

\begin{prop}\label{prop6.6}
For every term $f$ of \Ou\ we have that
$t(f)\cdot U(E(f))$ is $f$.
\end{prop}

Next we establish the following two propositions by inductions on
the length of derivation.

\begin{prop}\label{prop6.7}
If $f=g$ in \Oe, then $U(f)=U(g)$ in
\Ou.
\end{prop}

\begin{prop}\label{prop6.8}
If $f=g$ in \Ou, then $E(f)=E(g)$ in \Oe.
\end{prop}

\noindent The only case that is perhaps not quite straightforward
is with (\emph{assoc}~1$_u$) in the basis of the induction in the
proof of Proposition~\ref{prop6.8}. Here is how we proceed in that
case:
\begin{tabbing}
\hspace{4em}\=$E((h\ins g)\ins f)\;$\=$=(E(h)\ins_{t(h)\bsl
t(g)}
E(g))\ins_{t(h)\bsl t(f)} E(f)$,\\[1ex]
\>$E(h\ins (g\ins f))$\>$=E(h)\ins_{t(h)\bsl t(g)}
(E(g)\ins_{t(g)\bsl t(f)} E(f)).$
\end{tabbing}
Since we have that $t(f)$ is of the form $t(h)ba$, where $t(g)$ is
$t(h)b$, we may apply (\emph{assoc}~1$_e$).

As corollaries of the four propositions just established, we have
the following.

\begin{prop}\label{prop6.9}
$f=g$ in \Oe\ iff
$U(f)=U(g)$ in \Ou.
\end{prop}

\begin{prop}\label{prop6.10}
If $\;t(f)=t(g)$, then $f=g$ in \Ou\ iff
$E(f)=E(g)$ in \Oe.
\end{prop}

Consider the following set of operadic operations of \Ou\ for a
given \Np-word~$b$:
\begin{tabbing}
\hspace{4em}$C(\mathcal{O}_u(b))=_{df}\{f\in
C(\mathcal{O}_u)\mid t(f)=b\}$,
\end{tabbing}
and consider the structure
\begin{tabbing}
\hspace{4em}$\mathcal{O}_u(b)=_{df}\langle
C(\mathcal{O}_u(b)),\{\ins a\cdot\mid a\; \mbox{\rm is an
\Np-word}\},b\cdot\Iota\rangle$.
\end{tabbing}
The operation $\ins a\cdot$ is the partial operation on operadic
operations defined in \Ou, which applied to $f$ and $g$ yields
$g\ins a\cdot f$, provided $at(f)\in s(g)$. We can prove the
following.

\begin{prop}\label{prop6.11}
The structures \Oe\ and
$\mathcal{O}_u(e)$ are isomorphic.
\end{prop}

\begin{proof}
Propositions \ref{prop6.5}-\ref{prop6.8} show that, on the one
hand, the map $U$ and, on the other hand, the map $E$ restricted
to the terms of $\mathcal{O}_u(e)$ (i.e.\ those terms of \Ou\ that
stand for the members of $C(\mathcal{O}_u(e))$) induce bijections
inverse to each other between $C(\mathcal{O}_e)$ and
$C(\mathcal{O}_u(e))$. It remains only to consider the definition
of $U$ to establish that $U$, and hence $E$ too, are
homomorphisms.
\end{proof}

For a given \Np-word $b$ let
\begin{tabbing}
\hspace{4em}\=$C(\mathcal{O}_u(b\cdot))\;$\=$=_{df}\{f\in C(\mathcal{O}_u)
\mid t(f)=ba\; \mbox{\rm for some \Np-word}\; a\}$\\[1ex]
\>\>$=\cup\{C(\mathcal{O}_u(ba))\mid a\; \mbox{\rm is an
\Np-word}\}$,
\end{tabbing}
and consider the structure
\begin{tabbing}
\hspace{4em}$\mathcal{O}_u(b\cdot)=_{df}\langle
C(\mathcal{O}_u(b\cdot)),\ins,\{ba\cdot\Iota\mid a\; \mbox{\rm is
an \Np-word}\}\rangle$.
\end{tabbing}
Note that $\mathcal{O}_u(e\cdot)$ is \Ou.

We define inductively a map $b\bsl$ from
$C(\mathcal{O}_u(b\cdot))$ to $C(\mathcal{O}_u)$ by
\begin{tabbing}
\hspace{4em}\=$b\bsl(ba\cdot x)=a\cdot x$, \hspace{.5em}
$b\bsl(ba\cdot\Iota)=a\cdot\Iota$,\\[.5ex]
\>$b\bsl(g\ins f)=b\bsl g\ins b\bsl f$;
\end{tabbing}
$b\bsl f$ stands for $b\bsl(f)$, and we read $b\bsl g\ins b\bsl f$
as $(b\bsl g)\ins (b\bsl f)$. In the last clause, we have that
$b\bsl f$ is defined because $t(f)\in s(g)$ and $t(g)\in
P_{s(g)}$. We can prove the following.

\begin{prop}\label{prop6.12}
For every \Np-word $b$ the structures
\Ou\ and $\mathcal{O}_u(b\cdot)$ are isomorphic.
\end{prop}

\begin{proof}
From \Ou\ to $\mathcal{O}_u(b\cdot)$ we have the map
$b\cdot$, and it is easy to show that $b\bsl$ is its inverse.
\end{proof}

Let $\mathcal{O}_u^\ins(b)$ be the structure $\langle
C(\mathcal{O}_u(b)),\ins,b\cdot\Iota\rangle$, which is a
substructure of $\mathcal{O}_u(b)$ (from the family of insertions
$\{\ins a\cdot\mid a\; \mbox{\rm is an \Np-word}\}$ we keep only
$\ins e\cdot$, which amounts to $\ins$). It is also a substructure
of $\mathcal{O}_u(b\cdot)$, and is a kind of common denominator of
$\mathcal{O}_u(b)$ and $\mathcal{O}_u(b\cdot)$. We can prove the
following.

\begin{prop}\label{prop6.13}
For all \Np-words $a$ and $b$ the
structures $\mathcal{O}_u^\ins(a)$ and $\mathcal{O}_u^\ins(ba)$
are isomorphic.
\end{prop}

\begin{proof}
We restrict the maps $b\cdot$ and $b\bsl$, which we used for
the proof of the preceding proposition.
\end{proof}

\noindent As a corollary of this proposition we obtain that for
every \Np-word $b$ the structures $\mathcal{O}_u^\ins(e)$ and
$\mathcal{O}_u^\ins(b)$ are isomorphic.

Everything we said in this section about the relationship between
\Oe\ and \Ou\ can be restricted to the non-unitary structures
\Oem\ and \Oum.

\vspace{4ex}

{\samepage\noindent{\large {\sc Part II}}

\section{The category \WOum}}\label{sec7}
We introduce a category that will be a weakened version of \Ou,
with axiomatic equations replaced by isomorphisms. We deal first
in this section with the non-unitary category, and add what is
required for the unitary category in the next section. We deal in
these two sections with the $u$ versions of the weakened notions,
leaving for sections 9-10 the more complicated $e$ versions.

The object of \WOum\ are the terms of \Oum\ (not the operadic
operations of \Oum; see Section~5). We define inductively the
\emph{arrow terms} of \WOum. Every arrow term has a \emph{type},
which is a pair of objects $(f,g)$; as usual, we write $u\!:f\str
g$ to indicate that the arrow term $u$ is of that type. The object
$f$ is the \emph{source}, and $g$ the \emph{target}, of $u$. We
specify first the \emph{basic} arrow terms:
\begin{tabbing}
\hspace{4em}\=$\mj_f\!:f\str f$,
\\[1ex]
\>$\beta_{h,g,f}\!:(h\ins g)\ins f\str h\ins(g\ins f)$,
\hspace{1em} $\beta^{-1}_{h,g,f}\!:h\ins(g\ins f)\str(h\ins g)\ins
f$,
\\[1ex]
\>$\theta_{h,g,f}\!:(h\ins g)\ins f\str(h\ins f)\ins g$.
\end{tabbing}
We can make for these arrow terms comments on the legitimacy of
expressions for objects in their types exactly analogous to those
made for the equations (\emph{assoc}~1$_u$) and
(\emph{assoc}~2$_u$) in Section~5, from which these arrow terms
are derived.

Next we have the following two partial operations on arrow terms:
\begin{itemize}
\item[]{if $u\!:f\str f'$ and $v\!:g\str g'$ are arrow terms, then
$v\cirk u\!:f\str g'$ is an arrow term when $f'$ is $g$, and
$v\ins u\!:g\ins f\str g'\ins f'$ is an arrow term when $g\ins f$
and $g'\ins f'$ are legitimate.}
\end{itemize}
This concludes the definition of the arrow terms of \WOum. We use
$u,v,w,\ldots,$ perhaps with indices, for arrow terms.

Note that $\ins$ occurs now on three levels: first, at the level
of nominal arities in $Y\ins_a X$, which underly the objects,
secondly, at the level of objects in $g\ins f$, which underly the
arrow terms, and thirdly, at the level of arrow terms in $v\ins
u$.

Since the arrow terms of \WOum\ are derived from the equations of
\Oum, we obtain immediately from what we established by induction
on the length of derivation for \Oum\ in Section~5 that for every
arrow term $u\!:f\str g$ of \WOum\ we have $s(f)=s(g)$ and
$t(f)=t(g)$.

The equations of \WOum\ between arrow terms are given through an
axiomatic equational system, and the arrows of \WOum\ will be
formally equivalence classes of arrow terms such that these
equations are satisfied. Besides $u=u$ and the \emph{categorial
equations} $u\cirk \mj_f=u=\mj_g\cirk u$, for $u\!:f\str g$, and
$(w\cirk v)\cirk u=w\cirk(v\cirk u)$, the \emph{axiomatic
equations} of \WOum\ are the following:
\begin{tabbing}
\hspace{0em}\=(\emph{ins}~1)\hspace{1.2em}\=$\mj_g\ins\mj_f=\mj_{g\ins
f}$,\\*[1ex] \>(\emph{ins}~2)\>$(v_2\cirk v_1)\ins(u_2\cirk
u_1)=(v_2\ins
u_2)\cirk(v_1\ins u_1)$,\\[1.5ex]
\>($\beta$~\emph{nat})\>$\beta_{h_2,g_2,f_2}\cirk((w\ins v)\ins
u)\;$\=$=(w\ins(v\ins u))\cirk\beta_{h_1,g_1,f_1}$,\\*[1ex]
\>($\theta$~\emph{nat})\>$\theta_{h_2,g_2,f_2}\cirk((w\ins v)\ins
u)$\>$=((w\ins u)\ins v)\cirk\theta_{h_1,g_1,f_1}$,\\[1.5ex]
\>($\beta\beta$)\>$\beta^{-1}_{h,g,f}\cirk\beta_{h,g,f}\;$\=$=\mj_{(h\ins
g)\ins f}$, \hspace{1em}
$\beta_{h,g,f}\cirk\beta^{-1}_{h,g,f}=\mj_{h\ins (g\ins
f)}$,\\*[1ex]
\>($\theta\theta$)\>$\theta_{h,f,g}\cirk\theta_{h,g,f}$\>$=\mj_{(h\ins
g)\ins f}$,\\[1.5ex]
\>($\beta$~\emph{pent})\>$(\mj_j\ins\beta_{h,g,f})\cirk\beta_{j,h\ins
g,f}\cirk(\beta_{j,h,g}\ins\mj_f)=\beta_{j,h,g\ins
f}\cirk\beta_{j\ins h,g,f}$,\\[1ex]
\>($\theta$~YB)\>$\theta_{j\ins
f,h,g}\cirk(\theta_{j,h,f}\ins\mj_g)\cirk\theta_{j\ins
h,g,f}=(\theta_{j,g,f}\ins\mj_h)\cirk\theta_{j\ins
g,h,f}\cirk(\theta_{j,h,g}\ins\mj_f)$,\\[1.5ex]
\>($\beta\theta 1$)\>$(\mj_j\ins\theta_{h,g,f})\cirk\beta_{j,h\ins
g,f}\cirk(\beta_{j,h,g}\ins\mj_f)=\beta_{j,h\ins
f,g}\cirk(\beta_{j,h,f}\ins\mj_g)\cirk\theta_{j\ins
h,g,f}$,\\*[1ex] \>($\beta\theta 2$)\>$\theta_{j,h\ins
g,f}\cirk(\beta_{j,h,g}\ins\mj_f)=\beta_{j\ins
f,h,g}\cirk(\theta_{j,h,f}\ins\mj_g)\cirk\theta_{j\ins h,g,f}$.
\end{tabbing}
The name of ($\beta$~\emph{pent}) comes from Mac Lane's pentagon
of monoidal categories (see \cite{ML63} and \cite{ML98}, Section VII.1), while
($\theta$~YB) is related to the equation (YB) of Section 16 (and
YB comes from Yang-Baxter).

As rules we have symmetry and transitivity of $=$, and for $\xi$
being $\cirk$ and $\ins$ the congruence rules:
\begin{tabbing}
\hspace{4em}from $u_1=u_2$ and $v_1=v_2$ derive $v_1\,\xi\,
u_1=v_2\,\xi\, u_2$.
\end{tabbing}
This concludes our definition of the equations of \WOum, and of
the category \WOum.

An equation between the arrow terms of a category is legitimate
when both sides are of the same type, and one can easily check by
induction on the length of derivation that the equations of our
axiomatic system for \WOum\ satisfy this requirement for
legitimacy. The same holds for the equational axiomatic system of
all the categories introduced later, and we will not mention this
matter any more.

\section{The category \WOu}\label{sec8}
We add now to \WOum\ what is needed to obtain the unitary category
\WOu.

The objects of \WOu\ are the terms of \Ou\ (see Section~5). The
arrow terms of \WOu\ are defined like those of \WOum\ in the
preceding section with the following additional basic arrow terms,
derived from the equations (\emph{unit}$_u$) of Section~5:
\begin{tabbing}
\hspace{4em}\=$\mu_{f,a}\!:f\ins a\cdot\Iota\str f$,
\hspace{1em} \=$\mu^{-1}_{f,a}\!:f\str f\ins a\cdot\Iota$,
\\[1ex]
\>$\lambda_f\!:t(f)\cdot\Iota\ins f\str f$,\>
$\lambda^{-1}_f\!:f\str t(f)\cdot\Iota\ins f$.
\end{tabbing}

The equation of \WOu\ between arrow terms are defined like those
of \WOum\ in the preceding section with the following additional
axiomatic equations:
\begin{tabbing}
\hspace{0em}\=($\mu$~\emph{nat})\hspace{1.2em}\=
$\mu_{f_2,a}\cirk(u\ins\mj_{a\cdot\Iota})=u\cirk\mu_{f_1,a}$,\\*[1ex]
\>($\lambda$~\emph{nat})\>$\lambda_{f_2}\cirk(\mj_{t(f_1)\cdot\Iota}\ins
u)=
u\cirk\lambda_{f_1}$,\\[1.5ex]
\>($\mu\mu$)\>$\mu^{-1}_{f,a}\cirk\mu_{f,a}\;$\=$=\mj_{f\ins
a\cdot\Iota}$, \hspace{1em}
$\mu_{f,a}\cirk\mu^{-1}_{f,a}\;$\=$=\mj_f$,\\*[1ex]
\>($\lambda\lambda$)\>$\lambda^{-1}_f\cirk\lambda_f$\>$=\mj_{t(f)\cdot\Iota\ins
f}$, \hspace{1em}
$\lambda_f\cirk\lambda^{-1}_f\;$\>$=\mj_f$,\\[1.5ex]
\>($\beta\mu\lambda$)\>$\beta_{h,t(f)\cdot\Iota,f}=
(\mj_h\ins\lambda^{-1}_f)\cirk(\mu_{h,t(f)}\ins\mj_f)$,\\[1ex]
\>($\theta\mu$)\>$\theta_{h,b\cdot\Iota,f}=\mu^{-1}_{h\ins
f,b}\cirk(\mu_{h,b}\ins\mj_f)$.
\end{tabbing}

\section{The category \WOem}\label{sec9}
We introduce in this section the $e$ analogue of the non-unitary
category \WOum\ of Section~7.

The object of \WOem\ are the terms of \Oem\ (see Section~3). To
define inductively the \emph{arrow terms} of \WOem, we specify
first the \emph{basic} arrow terms:
\begin{tabbing}
\hspace{4em}\=$\mj_f$\=$:f\str f$,
\\[1ex]
\>$\beta_{h,(b,g),(a,f)}$\=$:(h\ins_b g)\ins_{ba} f\str
h\ins_b(g\ins_a f)$,\\[1ex]
\>$\beta^{-1}_{h,(b,g),(a,f)}$\>$:h\ins_b(g\ins_a f)\str(h\ins_b
g)\ins_{ba} f$,
\\[1ex]
\>$\theta_{h,(b,g),(a,f)}\!:(h\ins_b g)\ins_a f\str(h\ins_a
f)\ins_b g$.
\end{tabbing}
We can make for these arrow terms comments on the legitimacy of
expressions for objects in their types exactly analogous to those
made for the equations (\emph{assoc}~1$_e$) and
(\emph{assoc}~2$_e$) in Section~3, from which these arrow terms
are derived.

The operations under which the arrow terms are closed are
composition $\cirk$ and the operations $\ins_a$ for which we have
the following clause:
\begin{itemize}
\item[]{if $u\!:f\str f'$ and $v\!:g\str g'$ are arrow terms, then
$v\ins_a u\!:g\ins_a f\str g'\ins_a f'$ is an arrow term, when
$g\ins_a f$ and $g'\ins_a f'$ are legitimate.}
\end{itemize}
This concludes the definition of the arrow terms of \WOem. We use
still $u,v,w,\ldots,$ perhaps with indices, for these newly
introduced arrow terms. It follows immediately from what we
established for \Oem\ in Section~3 that for every arrow term
$u\!:f\str g$ of \WOem\ we have $s(f)=s(g)$.

The equations of \WOem\ between arrow terms are given through an
axiomatic equational system, which besides $u=u$ and the
categorial equations (as those given in Section~7 for \WOum) has
the following \emph{axiomatic equations}:
\begin{tabbing}
\hspace{0em}\=(\emph{ins}~1$_e$)\hspace{1.2em}\=$\mj_g\ins_a\mj_f=\mj_{g\ins_a
f}$,
\\*[1ex]
\>(\emph{ins}~2$_e$)\>$(v_2\cirk v_1)\ins_a(u_2\cirk
u_1)=(v_2\ins_a u_2)\cirk(v_1\ins_a u_1)$,
\\[1.5ex]
\>($\beta$~\emph{nat}$_e$)\>$\beta_{h_2,(b,g_2),(a,f_2)}\cirk((w\ins_b
v)\ins_{ba} u)\;$\=$=(w\ins_b(v\ins_a
u))\cirk\beta_{h_1,(b,g_1),(a,f_1)}$,
\\*[1ex]
\>($\theta$~\emph{nat}$_e$)\>$\theta_{h_2,(b,g_2),(a,f_2)}\cirk((w\ins_b
v)\ins_a u)$\>$=((w\ins_a u)\ins_b
v)\cirk\theta_{h_1,(b,g_1),(a,f_1)}$,
\\[1.5ex]
\>($\beta\beta_e$)\>$\beta^{-1}_{h,(b,g),(a,f)}\cirk\beta_{h,(b,g),(a,f)}\;$\=$=\mj_{(h\ins_b
g)\ins_{ba} f}$,\\*[1ex]
\>\>$\beta_{h,(b,g),(a,f)}\cirk\beta^{-1}_{h,(b,g),(a,f)}$\>$=\mj_{h\ins_b
(g\ins_a f)}$,
\\[1ex]
\>($\theta\theta_e$)\>$\theta_{h,(a,f),(b,g)}\cirk\theta_{h,(b,g),(a,f)}$\>$=\mj_{(h\ins_b
g)\ins_a f}$,
\\[1.5ex]
\>($\beta$~\emph{pent}$_e$)\>$(\mj_j\ins_c\beta_{h,(b,g),(a,f)})\cirk\beta_{j,(c,h\ins_b
g),(a,f)}\cirk(\beta_{j,(c,h)(b,g)}\ins_a\mj_f)=$\\[.5ex]
\`$\beta_{j,(c,h),(b,g\ins_a f)}\cirk\beta_{j\ins_c
h,(b,g),(a,f)}$,
\\[1ex]
\>($\theta$~YB$_e$)\>$\theta_{j\ins_a
f,(c,h),(b,g)}\cirk(\theta_{j,(c,h),(a,f)}\ins_b\mj_g)\cirk\theta_{j\ins_c
h,(b,g),(a,f)}=$\\[.5ex]
\`$(\theta_{j,(b,g),(a,f)}\ins_c\mj_h)\cirk\theta_{j\ins_b
g,(c,h),(a,f)}\cirk(\theta_{j,(c,h),(b,g)}\ins_a\mj_f)$,
\\[1.5ex]
\>($\beta\theta
1_e$)\>$(\mj_j\ins_c\theta_{h,(b,g),(a,f)})\cirk\beta_{j,(c,h\ins_b
g),(a,f)}\cirk(\beta_{j,(c,h),(b,g)}\ins_a\mj_f)=$\\[.5ex]
\`$\beta_{j,(c,h\ins_a
f),(b,g)}\cirk(\beta_{j,(c,h),(a,f)}\ins_b\mj_g)\cirk\theta_{j\ins_c
h,(b,g),(a,f)}$,
\\[1ex]
\>($\beta\theta 2_e$)\>$\theta_{j,(c,h\ins_b
g),(a,f)}\cirk(\beta_{j,(c,h),(b,g)}\ins_a\mj_f)=$\\[.5ex]
\`$\beta_{j\ins_a
f,(c,h),(b,g)}\cirk(\theta_{j,(c,h),(a,f)}\ins_b\mj_g)\cirk\theta_{j\ins_c
h,(b,g),(a,f)}$.
\end{tabbing}

As rules we have symmetry and transitivity of $=$, and the
congruence rules for $\cirk$ and $\ins_a$ (just put $\cirk$ and
$\ins_a$ for $\xi$ in the schema at the end of Section~7). This
concludes our definition of the equations of \WOem, and of the
category \WOem.

\section{The category \WOe}\label{sec10}
We add now to \WOem\ what is needed to obtain the unitary category
\WOe.

The objects of \WOe\ are the terms of \Oe\ (see Section~3). The
arrow terms of \WOe\ are defined like those of \WOem\ in the
preceding section with the following additional basic arrow terms,
derived from the equations (\emph{unit}$_e$) of Section~3:
\begin{tabbing}
\hspace{4em}\=$\mu_{f,a}\!:f\ins_a\Iota\str f$, \hspace{1em}
\=$\mu^{-1}_{f,a}\!:f\str f\ins_a\Iota$,
\\[1ex]
\>$\lambda_f\!:\Iota\ins_e f\str f$,\> $\lambda^{-1}_f\!:f\str
\Iota\ins_e f$.
\end{tabbing}

The equation of \WOe\ between arrow terms are defined like those
of \WOem\ in the preceding section with the following additional
axiomatic equations:
\begin{tabbing}
\hspace{0em}\=($\mu$~\emph{nat}$_e$)\hspace{1.2em}\=
$\mu_{f_2,a}\cirk(u\ins_a\mj_{\Iota})=u\cirk\mu_{f_1,a}$,\\*[1ex]
\>($\lambda$~\emph{nat}$_e$)\>$\lambda_{f_2}\cirk(\mj_{\Iota}\ins_e
u)=
u\cirk\lambda_{f_1}$,\\[1.5ex]
\>($\mu\mu_e$)\>$\mu^{-1}_{f,a}\cirk\mu_{f,a}\;$\=$=\mj_{f\ins_a\Iota}$,
\hspace{1em} $\mu_{f,a}\cirk\mu^{-1}_{f,a}\;$\=$=\mj_f$,\\*[1ex]
\>($\lambda\lambda_e$)\>$\lambda^{-1}_f\cirk\lambda_f$\>$=\mj_{\Iota\ins_e
f}$, \hspace{1em}
$\lambda_f\cirk\lambda^{-1}_f\;$\>$=\mj_f$,\\[1.5ex]
\>($\beta\mu\lambda_e$)\>$\beta_{h,(b,\Iota),(e,f)}\;$\=$=
(\mj_h\ins_b\lambda^{-1}_f)\cirk(\mu_{h,b}\ins_b\mj_f)$,\\[1ex]
\>($\theta\mu_e$)\>$\theta_{h,(b,\Iota),(a,f)}$\>$=\mu^{-1}_{h\ins_a
f,b}\cirk(\mu_{h,b}\ins_a\mj_f)$.
\end{tabbing}

\section{\WOe\ and \WOu}\label{sec11}
In this section we establish the relationship between \WOe\ and
\WOu. We show that \WOu\ is the disjoint union of isomorphic
copies of \WOe.

For every \Np-word $a$ we define inductively a map $a\cdot$ from
the arrow terms of \WOu\ to the arrow terms of \WOu:
\begin{tabbing}
\hspace{4em}\=$a\cdot\nu_f=\nu_{a\cdot f}$, \hspace{.2em} where
$\nu$ is $\mj$, $\lambda$ and $\lambda^{-1}$,\\[1ex]
\>$a\cdot\nu_{f,b}=\nu_{a\cdot f,ab}$, \hspace{.2em} where $\nu$
is $\mu$ and $\mu^{-1}$,\\[1ex]
\>$a\cdot\zeta_{h,g,f}=\zeta_{a\cdot h,a\cdot g,a\cdot f}$,
\hspace{.2em} where $\zeta$ is $\beta$, $\beta^{-1}$ and
$\theta$,\\[1ex]
\>$a\cdot(v\,\xi\,u)=(a\cdot v)\,\xi\,(a\cdot u)$, \hspace{.2em}
where $\xi$ is $\cirk$ and $\ins$.
\end{tabbing}
The definition of $a\cdot$ on the objects of \WOu, i.e.\ the terms
of \Ou, which is mentioned in the indices above, is given in
Section~6. It is clear that $a\cdot$ induces an endofunctor of
\WOu, since $u=v$ in \WOu\ clearly implies $a\cdot u=a\cdot v$ in
\WOu.

We need $a\cdot$ to define inductively a map $U$ from the arrow
terms of \WOe\ to the arrow terms of \WOu:
\begin{tabbing}
\hspace{4em}\=$U(\nu_f)=\nu_{U(f)}$, \hspace{.2em} where
$\nu$ is $\mj$, $\lambda$ and $\lambda^{-1}$,
\\[1ex]
\>$U(\nu_{f,a})=\nu_{U(f),a}$, \hspace{.2em} where $\nu$
is $\mu$ and $\mu^{-1}$,
\\[1ex]
\>$U(\zeta_{h,(b,g),(a,f)})\;$\=$=\zeta_{U(h),b\cdot U(g),ba\cdot
U(f)}$, \hspace{.2em} where $\zeta$ is $\beta$ and $\beta^{-1}$,
\\[1ex]
\>$U(\theta_{h,(b,g),(a,f)})$\>$=\theta_{U(h),b\cdot U(g),a\cdot
U(f)}$,
\\[1ex]
\>$U(v\cirk u)=U(v)\cirk U(u)$,
\\[1ex]
\>$U(v\ins_a u)=U(v)\ins a\cdot U(u)$.
\end{tabbing}
The map $U$ mentioned in the indices of this definition is the map
$U$ defined in Section~6.

Next we define inductively a map $E$ from the arrow terms of \WOu\
to the arrow terms of \WOe:
\begin{tabbing}
\hspace{4em}\=$E(\nu_f)=\nu_{E(f)}$, \hspace{.2em} where $\nu$
is $\mj$, $\lambda$ and $\lambda^{-1}$,
\\[1ex]
\>$E(\nu_{f,a})=\nu_{E(f),a}$, \hspace{.2em} where $\nu$ is $\mu$
and $\mu^{-1}$,
\\[1ex]
\>$E(\zeta_{h,g,f})\;$\=$=\zeta_{E(h),(t(h)\bsl
t(g),E(g)),(t(g)\bsl t(f),E(f))}$, \hspace{.2em} where $\zeta$ is
$\beta$ and $\beta^{-1}$,
\\[1ex]
\>$E(\theta_{h,g,f})$\>$=\theta_{E(h),(t(h)\bsl
t(g),E(g)),(t(h)\bsl t(f),E(f))}$,
\\[1ex]
\>$E(v\cirk u)=E(v)\cirk E(u)$,
\\[1ex]
\>$E(v\ins u)=E(v)\ins_{t(g)\bsl t(f)} E(u)$, \hspace{.2em} for
$u\!:f\str f'$ and $v\!:g\str g'$.
\end{tabbing}
The map $E$ mentioned in the indices of this definition is the map
$E$ defined in Section~6.

With the help of ($st\;a\cdot$) and Lemma \ref{lem6.4} of
Section~6, it is easy to establish the following by induction on
the complexity of~$u$.

\begin{lem}\label{lem11.1}
For every arrow term $u$ of \WOu\ we have that $E(a\cdot u)$ is
$E(u)$.
\end{lem}

\noindent This lemma serves for the first of the following two
propositions, which are proved by inductions on the complexity
of~$u$.

\begin{prop}\label{prop11.2}
For every arrow term $u$ of \WOe\ we have that $E(U(u))$ is~$u$.
\end{prop}

\begin{prop}\label{prop11.3}
For every arrow term $u\!:f\str f'$ of
\WOu\ we have that $t(f)\cdot U(E(u))$ is~$u$.
\end{prop}

\noindent In the proofs of these propositions we rely on
Propositions \ref{prop6.5} and \ref{prop6.6}. As an example, which
is perhaps not quite straightforward, we give the following case
in the proof of Proposition~\ref{prop11.3}:
\begin{tabbing}
\hspace{4em}$t(h)\cdot U(E(\beta_{h,g,f}))=t(h)\cdot
U(\beta_{E(h),(t(h)\bsl t(g),E(g)),(t(g)\bsl t(f),E(f))})$\\*[1ex]
\hspace{8em}\=$=\beta_{t(h)\cdot U(E(h)),t(h)(t(h)\bsl t(g))\cdot
U(E(g)),
t(h)(t(h)\bsl t(g))(t(g)\bsl t(f))\cdot U(E(f))}$\\[1ex]
\>$=\beta_{t(h)\cdot U(E(h)),t(g)\cdot U(E(g)),t(f)\cdot
U(E(f))}=\beta_{h,g,f}$, \hspace{.2em} by
Proposition~\ref{prop6.6}.
\end{tabbing}\bigskip

Next we establish the following two propositions by inductions on
the length of derivation, where the main burden is in the bases of
the inductions, with axiomatic equations.

\begin{prop}\label{prop11.4}
If $u=v$ in \WOe, then $U(u)=U(v)$ in
\WOu.
\end{prop}

\begin{prop}\label{prop11.5}
If $u=v$ in \WOu, then $E(u)=E(v)$ in
\WOe.
\end{prop}

\noindent The inductive proofs of these two propositions, which
are lengthy but straightforward, yield more than what is stated in
the propositions. Every derivation in the equational system of
\WOe\ is translated into a derivation in the equational system of
\WOu\ in a ``homomorphic'' manner, and vice versa. This means, for
example, that ($\beta\theta 1_e$) goes into ($\beta\theta 1$),
etc. These two propositions yield that $U$ induces a functor from
\WOe\ to \WOu, and $E$ a functor in the opposite direction.

Let $\mathcal{WO}_e(X)$ be the full subcategory of \WOe\ whose
objects are all the terms $f$ of \Oe\ such that $s(f)=X$. The
category \WOe\ is the disjoint union of its subcategories
$\mathcal{WO}_e(X)$ for all the nominal arities $X$.

Let $\mathcal{WO}_u(X,b)$ be the full subcategory of \WOu\ whose
objects are all the terms $f$ of \Ou\ such that $s(f)=X$ and
$t(f)=b$. Let \WOub\ be the union of the categories
$\mathcal{WO}_u(X,b)$ for all the nominal arities $X$ (this union
is disjoint). The category \WOu\ is the disjoint union of its
subcategories \WOub\ for all \Np-words $b$. We can prove the
following.

\begin{prop}\label{prop11.6}
The categories \WOe\ and
$\mathcal{WO}_u(e)$ are isomorphic.
\end{prop}

\begin{proof}
Propositions \ref{prop11.2}-\ref{prop11.5} show that, on the
one hand, the functor $U$ and, on the other hand, the functor $E$
restricted to $\mathcal{WO}_u(e)$ are inverse to each other.
\end{proof}

\noindent More precisely, we have that $\mathcal{WO}_e(X)$ and
$\mathcal{WO}_u(X,e)$ are isomorphic categories. For that we just
restrict further the functors $U$ and $E$.

Let \WOubp\ be the union of the categories $\mathcal{WO}_u(ba)$
for all \Np-words $a$. We define a functor $b\bsl$ from \WOubp\ to
\WOu\ by stipulating that $b\bsl u$ is $u$ with every index $ba$
replaced by $a$; moreover, $b\bsl(v\,\xi\,u)$ is $(b\bsl v)\,\xi\,
(b\bsl u)$ for $\xi$ being $\cirk$ and $\ins$ (cf. the definition
of the functor $a\cdot$ at the beginning of the section). We can
prove the following.

\begin{prop}\label{prop11.7}
For every \Np-word $b$ the categories
\WOu\ and \WOubp\ are isomorphic.
\end{prop}

\begin{proof}
From \WOu\ and \WOubp\ we have the functor $b\cdot$, and it
is easy to see that the functor $b\backslash$ is its inverse.
\end{proof}

We can also prove the following.

\begin{prop}\label{prop11.8}
For all \Np-words $a$ and $b$ the
categories $\mathcal{WO}_u(a)$ and $\mathcal{WO}_u(ba)$ are
isomorphic.
\end{prop}

\begin{proof}
We restrict the functors $b\cdot$ and $b\bsl$ to the
subcategory $\mathcal{WO}_u(a)$ of \WOu\ and the subcategory
$\mathcal{WO}_u(ba)$ of \WOubp.
\end{proof}

As a corollary of this proposition, we obtain that for every
\Np-word $b$ the categories $\mathcal{WO}_u(e)$ and \WOub\ are
isomorphic, and hence, with Proposition~\ref{prop11.6}, we have
that \WOe\ and \WOub\ are isomorphic. Since the category \WOu\ is
the disjoint union of the categories \WOub\ for all \Np-words $b$,
we may conclude that \WOu\ is the disjoint union of isomorphic
copies of \WOe.

Propositions \ref{prop6.11}, \ref{prop6.12} and \ref{prop6.13} are
parallel to Propositions \ref{prop11.6}, \ref{prop11.7} and
\ref{prop11.8}. The propositions of Section~6 deal with algebraic
structures with partial operations satisfying various equations
related to (\emph{unit}), (\emph{assoc}~1) and (\emph{assoc}~2),
while the categories \WOe\ and \WOu\ do not satisfy these
equations, but have arrows that are isomorphisms instead of them.
Nevertheless, the former propositions indicate how the functors of
this section do not preserve only $\mj$ and $\cirk$, but also
insertion, as a partial operation both on objects and on arrows.
What is preserved is either the whole family of insertion
operations indexed with \Np-words, or just the single partial
operation $\ins$ (which corresponds to $\ins_e$). What is
preserved is also $\Iota$, or objects derived from it.

Everything we said in this section about the relationship between
\WOe\ and \WOu\ can be restricted to the non-unitary categories
\WOem\ and \WOum.

\section{Operads, Cat-operads and weak Cat-operads}\label{sec12}
{\sc Definition of operad.} The standard general notion of operad
with insertions primitive would be based on the operad \O\ of
Section~2. It defines a class of partial algebras in which \O\ is
freely generated by $G$ and $\alpha_G$. These algebras have a
carrier $C$ made of elements called operadic operations, such that
for every operadic operation $\varphi$ of $C$ we have an arity in
$\mathbf N$, and they have the family of insertions $\{\ins_k\mid
k\in{\mathbf N}^+\}$, which are partial operations on operadic
operations, and $\Iota$ with arity 1, for which equations like
those of \O\ hold. The notion of multicategory is a generalization
of this notion where an arity $n$ is replaced by a sequence of $n$
occurrences of some objects, and moreover we have an object as a
target (see \cite{L89}). An operad is a one-object multicategory.

\vspace{2ex}

\noindent {\sc First definition of Cat-operad.} A Cat-operad is an
operad that in addition has arrows between operadic operations of
the same arity. We have an identity arrow for every operadic
operation, and the arrows are closed under the partial operations
of composition and insertions, which are now not only partial
operations on operadic operations, but also partial operations on
arrows. This structure is a category, i.e.\ identity arrows and
composition satisfy the categorial equations (see Section~7).
Since the enrichment is over Cat, the category of all \emph{small}
categories, it is assumed that the hom-categories of this category
are small, but this is not an essential matter. We have moreover
in this category the equations (\emph{ins}~1$_e$) and
(\emph{ins}~2$_e$) of Section~9 with $a$ replaced by $n$, and the
analogues of the equations (\emph{unit}), (\emph{assoc}~1) and
(\emph{assoc}~2) for arrows (see Section~2). In these analogues
$\Iota$ is replaced by $\mj_\Iota$ and $\varphi$, $\gamma$ and
$\chi$ are replaced by variables for arrows. If $u$ is such a
variable, then $\alpha(u)$ is the arity of the source or target of
$u$, which must have the same arity.

\vspace{2ex}

\noindent {\sc Second definition of Cat-operad.} In other words, a
Cat-operad can be defined as the disjoint union of categories
$\mathcal{C}_k$; in $\mathcal{C}_k$ all the objects, called
operadic operations, have arity $k$. In $\mathcal{C}_1$ we have a
special object $\Iota$. We have moreover the bifunctors
$\ins_n\!:\mathcal{C}_k\times\mathcal{C}_l \str
\mathcal{C}_{k-1+l}$, for ${1\leq n\leq k}$, which satisfy the
equations (\emph{unit}), (\emph{assoc}~1) and (\emph{assoc}~2),
and analogous equations for arrows.


\vspace{2ex}

We base our general notion of weak Cat-operad on \WOe. Let
$\langle G,\alpha_G,\mathcal{G}\rangle$ be a triple (not in the
sense of \emph{monad}), where $G$ is a set, $\alpha_G$  is a
function from $G$ to $\textbf{N}$ and $\mathcal{G}$ is a directed
graph $G\stackrel{dom}\longleftarrow
A\stackrel{cod}\longrightarrow G$ (in the sense of \cite{ML98},
Section I.2) such that for every $u$ in $A$ we have
\[
\alpha_G(dom(u))=\alpha_G(cod(u)).
\]
Let $X$ be a category whose objects are such triples and whose
arrows from $\langle G,\alpha_G,\mathcal{G}\rangle$ to $\langle
G',\alpha_{G'},\mathcal{G}'\rangle$ are graph morphisms
$\varphi\!:\mathcal{G}\str \mathcal{G}'$ such that for every $x$
in $G$ we have that $\alpha_G(x)=\alpha_{G'}(\varphi(x))$.

We define the category $\mbox{\WOe}\langle
G,\alpha_G,\mathcal{G}\rangle$ as we have defined \WOe\ in
Sections 9-10 based on $G$ and $\alpha_G$ save that in the
definition of arrow terms we add the clause
\begin{itemize}
\item[]{if $u\!:x\str y$ is in $A$, then it is an arrow term.}
\end{itemize}
We say that a functor $F\!:\mbox{\WOe}\langle
G,\alpha_G,\mathcal{G}\rangle\str \mbox{\WOe}\langle
G',\alpha_{G'},\mathcal{G}'\rangle$ preserves the \WOe\ structure
when
\begin{tabbing}
for objects:
\\[1ex]
\hspace{4em}\=$F(\Iota)=\Iota$;
\\[1ex]
\>$F(g\ins_a f)=Fg\ins_{a'} Ff$, where
$K_{s(g)}(a)=K_{s(Fg)}(a')$;
\\[2ex]
and for arrows:
\\[1ex]
\>$F(\beta_{h,(b,g),(a,f)})=\beta_{Fh,(b',Fg),(a',Ff)}$, where
$K_{s(g)}(a)=K_{s(Fg)}(a')$ and
\\
\` $K_{s(h)}(b)=K_{s(Fh)}(b')$;
\\[1ex]
\>similarly for $\beta^{-1}$, $\theta$, $\mu$, $\mu^{-1}$, $\lambda$
and $\lambda^{-1}$;
\\[1ex]
\>$F(v\ins_a u)=Fv\ins_{a'} Fu$, where $u\!:f\str f'$, $v\!:g\str
g'$ and
\\
\` $K_{s(g)}(a)=K_{s(Fg)}(a')$
\end{tabbing}
(see Section~4 for the definition of $K$).

Let $Y$ be the category whose objects are all the categories
$\mbox{\WOe}\langle G,\alpha_G,\mathcal{G}\rangle$ and whose
arrows are all the functors that preserve the \WOe\ structure.
Every arrow $\varphi\!:\langle G,\alpha_G,\mathcal{G}\rangle \str
\langle G',\alpha_{G'},\mathcal{G}'\rangle$ of $X$ induces a
function $\varphi_0$ from the objects of $\mbox{\WOe}\langle
G,\alpha_G,\mathcal{G}\rangle$ to the objects of
$\mbox{\WOe}\langle G',\alpha_{G'},\mathcal{G}'\rangle$ and a
function $\varphi_1$ from the arrow terms of $\mbox{\WOe}\langle
G,\alpha_G,\mathcal{G}\rangle$ to the arrow terms of
$\mbox{\WOe}\langle G',\alpha_{G'},\mathcal{G}'\rangle$ in a
natural way. It is straightforward to verify that for the arrow
terms $u$ and $v$, if $u=v$ in $\mbox{\WOe}\langle
G,\alpha_G,\mathcal{G}\rangle$, then $\varphi_1(u)=\varphi_1(v)$
in $\mbox{\WOe}\langle G',\alpha_{G'},\mathcal{G}'\rangle$. Hence
the functions $\varphi_0$ and $\varphi_1$ underlie a functor
$L\varphi$ from $\mbox{\WOe}\langle G,\alpha_G,\mathcal{G}\rangle$
to $\mbox{\WOe}\langle G',\alpha_{G'},\mathcal{G}'\rangle$. It is
not difficult to check that $L\varphi$ preserves the \WOe\
structure.

Hence we have a functor $L\!:X\str Y$ which maps $\langle
G,\alpha_G,\mathcal{G}\rangle$ to $\mbox{\WOe}\langle
G,\alpha_G,\mathcal{G}\rangle$ and $\varphi$ to $L\varphi$. On the
other hand, we can define the forgetful functor $R\!:Y\str X$ so
that for $C$ being $\mbox{\WOe}\langle
G,\alpha_G,\mathcal{G}\rangle$ we have that $R(C)=\langle
H,\alpha_H,\mathcal{H}\rangle$, where
\begin{tabbing}
\hspace{4em}\=$H=Ob(C)$, i.e.\ the set of terms of \Oe\ based on
$G$ and $\alpha_G$,
\\[1ex]
\>$\alpha_H(f)=|s(f)|$, and
\\[1ex]
\>$\mathcal{H}$ is the directed graph underlying the category $C$.
\end{tabbing}

\begin{prop}\label{prop12.1}
The functor $L$ is the left adjoint of
$R$.
\end{prop}

\begin{proof}
The proof that $Y(\mbox{\WOe}\langle
G,\alpha_G,\mathcal{G}\rangle, \mbox{\WOe}\langle
G',\alpha_{G'},\mathcal{G}'\rangle)$ is naturally isomorphic to
\linebreak $X(\langle G,\alpha_G,\mathcal{G}\rangle,
R(\mbox{\WOe}\langle G',\alpha_{G'},\mathcal{G}'\rangle))$ is
similar to the proof that $L\varphi$ is a functor that preserves
the \WOe\ structure.
\end{proof}

\noindent {\sc Definition of weak Cat-operad.} A weak Cat-operad
is an algebra of the monad in $X$ defined by the above adjunction.


\vspace{2ex}

Alternatively, a weak Cat-operad can be defined in the style of the second
definition of Cat-operad above by copying this definition until we
reach the equations (\emph{unit}), (\emph{assoc}~1) and
(\emph{assoc}~2). Instead of these equations, we have natural
isomorphisms corresponding to $\mu$, $\lambda$, $\beta$ and
$\theta$ and equations analogous to the equations
($\beta\mu\lambda_e$), ($\theta\mu_e$), ($\beta$~\emph{pent}$_e$),
($\theta$~YB$_e$), ($\beta\theta 1_e$) and ($\beta\theta 2_e$) of
Sections 10 and~9. For example, the equation analogous to
($\beta\theta 1_e$) would be the following:

\begin{tabbing}
\hspace{4em}$(\mj_j\ins_m\theta_{h,(l,g),(k-1+\alpha(g),f)})\cirk\beta_{j,(m,h\ins_l
g),(k',f)}\cirk(\beta_{j,(m,h),(l',g)}\ins_{k'}\mj_f)=$
\\[.75ex]
\`$\beta_{j,(m,h\ins_k
f),(l',g)}\cirk(\beta_{j,(m,h),(k+m-1,f)}\ins_{l'}\mj_g)\cirk\theta_{j\ins_m
h,(l',g),(k',f)}$,
\end{tabbing}\medskip
where $k'=k\mn1\pl\alpha(g)\pl m\mn 1$ and $l'=l\pl m\mn 1$, with
$k$, $l$ and $m$ standing respectively for $K_{s(h)}(a)$,
$K_{s(h)}(b)$ and $K_{s(h)}(c)$ (see Section~4 for the definition
of $K$).

This equation is considerably more complicated than our equation
($\beta\theta 1_e$), and to eschew such complications is the main
reason for introducing our nominal arities. To infer all the
equations needed for a full definition of weak Cat-operads, as we
inferred the equation above from ($\beta\theta 1_e$), is a
lengthy, but straightforward, matter, into which we will not go
further in this paper.

\section{\WOem\ and hemiassociahedra}\label{sec13}
For particular choices of the nominal arity $X$ the categories
\WOemX\ defined in Section 11 have an interesting shape. They
become representable by polyhedra of a kind analogous to
associahedra and permutohedra (see \cite{S97}, \cite{S97a},
\cite{Z95}, Lecture~0, Example 0.10, and \cite{GR63} for
historical references concerning associahedra and permutohedra).
In some cases they are exactly associahedra that involve only
$\beta$ arrows, and permutohedra that involve only $\theta$
arrows. These associahedra and permutohedra do not differ from
those already considered in the literature in connection with
associativity and commutativity isomorphisms of monoidal and
symmetric strictly monoidal categories.

It is interesting however to describe associahedra and
permutohedra with $\beta$ and $\theta$ arrows mixed. With these
two kinds of arrows mixed we obtain also in three dimensions
another kind of polyhedron, which was called
\emph{hemiassociahedron} in \cite{DP10} (Example 5.14). That paper
investigates in particular the relationship of the
hemiassociahedron to the permutohedron (from which, together with
\cite{T97}, one may gather that the hemiassociahedron, conceived
as an abstract polytope, can be realized; we will however not go
here into this problem, for which one should also consult \cite{DP11}).

In this section we describe nine categories \WOemX. The first four
may be represented by hemiassociahedra; next we have two that may
be represented by three-dimensional associahedra with $\beta$ and
$\theta$ mixed, and one that may be represented by a
three-dimensional permutohedron with $\beta$ and $\theta$ arrows
mixed. Finally we have a category represented by a purely $\beta$
associahedron, and a category represented by a purely $\theta$
permutohedron.

In all our examples, instead of dealing with \WOemX\ we may deal
with the isomorphic category $\mathcal{WO}_u^-(X,e)$. In our first
example, which is given with more details, we concentrate on this
category, simpler to deal with. In other examples, with less
details, we would proceed analogously.

\begin{exa}\label{exa13.1}
Let $G=\{x\}$ and let $\alpha_G(x)=2$. Let $X=\{111, 112, 121,
122, 21, 22\}$. Every object of $\mathcal{WO}_u^-(X,e)$, which is
a term $f$ of \Oum\ with $s(f)=X$ and $t(f)=e$, records a
destruction of the following binary tree, whose leaves make $X$:
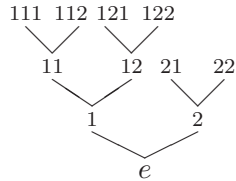
\medskip
\begin{center}
\begin{picture}(90,70)
\put(45,10){\line(-2,1){20}} \put(45,10){\line(2,1){20}}

\put(25,30){\line(-3,2){15}} \put(25,30){\line(3,2){15}}

\put(65,30){\line(-1,1){10}} \put(65,30){\line(1,1){10}}

\put(10,50){\line(-1,1){10}} \put(10,50){\line(1,1){10}}
\put(40,50){\line(-1,1){10}} \put(40,50){\line(1,1){10}}

\put(45,5){\makebox(0,0){$e$}}
\put(25,25){\makebox(0,0){\scriptsize $1$}}
\put(65,25){\makebox(0,0){\scriptsize $2$}}
\put(10,45){\makebox(0,0){\scriptsize $11$}}
\put(40,45){\makebox(0,0){\scriptsize $12$}}
\put(55,45){\makebox(0,0){\scriptsize $21$}}
\put(75,45){\makebox(0,0){\scriptsize $22$}}
\put(0,65){\makebox(0,0){\scriptsize $111$}}
\put(17,65){\makebox(0,0){\scriptsize $112$}}
\put(33,65){\makebox(0,0){\scriptsize $121$}}
\put(50,65){\makebox(0,0){\scriptsize $122$}}

\end{picture}
\end{center}
For example, the object $f$, which is
\begin{tabbing}
\hspace{4em}$(e\cdot x\ins 2\cdot x)\ins((1\cdot x\ins 11\cdot
x)\ins 12\cdot x)$,
\end{tabbing}
records the following destruction. The main insertion $\ins$ of
our object is first removed. This leaves us with the \Oum\ terms
$f_1$ and $f_2$, which are respectively
\begin{tabbing}
\hspace{4em}$(1\cdot x\ins 11\cdot x)\ins 12\cdot x$\quad
and\quad $e\cdot x\ins 2\cdot x$,
\end{tabbing}
and which record destructions of the following two trees:\medskip
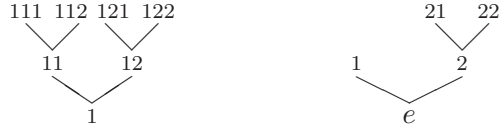
\begin{center}
\begin{picture}(180,50)
\put(145,10){\line(-2,1){20}} \put(145,10){\line(2,1){20}}

\put(25,10){\line(-3,2){15}} \put(25,10){\line(3,2){15}}

\put(165,30){\line(-1,1){10}} \put(165,30){\line(1,1){10}}

\put(10,30){\line(-1,1){10}} \put(10,30){\line(1,1){10}}
\put(40,30){\line(-1,1){10}} \put(40,30){\line(1,1){10}}

\put(145,5){\makebox(0,0){$e$}}
\put(25,5){\makebox(0,0){\scriptsize $1$}}
\put(125,25){\makebox(0,0){\scriptsize $1$}}
\put(165,25){\makebox(0,0){\scriptsize $2$}}
\put(10,25){\makebox(0,0){\scriptsize $11$}}
\put(40,25){\makebox(0,0){\scriptsize $12$}}
\put(155,45){\makebox(0,0){\scriptsize $21$}}
\put(175,45){\makebox(0,0){\scriptsize $22$}}
\put(0,45){\makebox(0,0){\scriptsize $111$}}
\put(17,45){\makebox(0,0){\scriptsize $112$}}
\put(33,45){\makebox(0,0){\scriptsize $121$}}
\put(50,45){\makebox(0,0){\scriptsize $122$}}

\end{picture}
\end{center}
We have $t(f_1)=1\in s(f_2)$.

Next we concentrate on the destruction recorded by $f_1$. We
remove its main insertion, and next the insertion remaining on its
left, which leaves three residual trees:\medskip
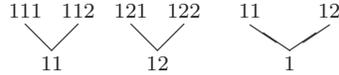
\begin{center}
\begin{picture}(120,30)
\put(100,10){\line(-3,2){15}} \put(100,10){\line(3,2){15}}

\put(10,10){\line(-1,1){10}} \put(10,10){\line(1,1){10}}

\put(50,10){\line(-1,1){10}} \put(50,10){\line(1,1){10}}

\put(100,5){\makebox(0,0){\scriptsize $1$}}
\put(10,5){\makebox(0,0){\scriptsize $11$}}
\put(50,5){\makebox(0,0){\scriptsize $12$}}
\put(0,25){\makebox(0,0){\scriptsize $111$}}
\put(20,25){\makebox(0,0){\scriptsize $112$}}
\put(40,25){\makebox(0,0){\scriptsize $121$}}
\put(60,25){\makebox(0,0){\scriptsize $122$}}
\put(85,25){\makebox(0,0){\scriptsize $11$}}
\put(115,25){\makebox(0,0){\scriptsize $12$}}

\end{picture}
\end{center}
This destruction of $f_1$ may be recorded also by $(12)\cdot(11)$,
which indicate that we broke the tree of $f_1$ first at the vertex
12, and next broke one of the remaining trees at the vertex 11.

The destruction of $f_2$ proceeds by breaking its tree at the
vertex~2, which leaves the two trees
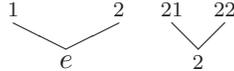
\begin{center}
\begin{picture}(80,30)
\put(20,10){\line(-2,1){20}} \put(20,10){\line(2,1){20}}

\put(70,10){\line(-1,1){10}} \put(70,10){\line(1,1){10}}

\put(20,5){\makebox(0,0){$e$}}
\put(0,25){\makebox(0,0){\scriptsize $1$}}
\put(40,25){\makebox(0,0){\scriptsize $2$}}
\put(70,5){\makebox(0,0){\scriptsize $2$}}
\put(60,25){\makebox(0,0){\scriptsize $21$}}
\put(80,25){\makebox(0,0){\scriptsize $22$}}

\end{picture}
\end{center}
This destruction may also be recorded by 2, the vertex where the
tree broke.

The destruction of $f_1$ and $f_2$ may have gone simultaneously,
and we indicate that by writing $((12)\cdot(11))\pl 2$, which is
equal to $2\pl((12)\cdot(11))$. On the other hand, $\cdot$ is not
commutative, and indicates successive steps of destruction. Both
$\cdot$ and $+$ are on the other hand associative. Our complete
destruction of $f$ may then be recorded by
$1\cdot(((12)\cdot(11))\pl 2)$. (Every destruction is, of course,
a construction in reverse order.)
\end{exa}

Consider terms built out as in this example out of the names of
the inner vertices of a tree with the help of the operations
$\cdot$ and $+$. Some of these terms are in one-to-one
correspondence with the terms of \Oum\ in our examples, so that
the associative and non-commutative operation $\cdot$ corresponds
to successive steps of destruction, while the associative and
commutative operation $+$ corresponds to simultaneous steps of
destruction. These terms are shorter than the terms of \Oum.

This matter is treated more formally in \cite{DP10}, where such
terms with $\cdot$ and $+$ are called \textbf{S}-trees. (The
one-to-one correspondence mentioned in the preceding paragraph
exists because the value of $\alpha_G$ is greater than or equal to
$2$, and $\alpha_G$ happens to be one-one.) To shorten these terms
further we will write $a$, $b$, $c$ and $d$ for respectively 11,
12, 1 and 2. The \Np-words $a$, $b$, $c$ and $d$ stand for the
inner vertices of the tree of $f$, i.e.\ vertices that are neither
leaves nor roots. These vertices in this tree make the following
graph:\medskip
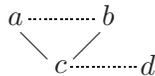
\begin{center}
\begin{picture}(50,20)(0,5)
\put(10,10){\line(-1,1){10}} \put(20,10){\line(1,1){10}}
\multiput(3,25)(2,0){13}{\line(1,0){.5}}
\multiput(19,7)(2,0){13}{\line(1,0){.5}}

\put(15,7){\makebox(0,0){$c$}} \put(-2,25){\makebox(0,0){$a$}}
\put(33,26.5){\makebox(0,0){$b$}} \put(48,8){\makebox(0,0){$d$}}
\end{picture}
\end{center}
where the solid edge $\{c,a\}$ indicates that the vertices $c$ and
$a$ are immediately one above the other in the tree; and
analogously for the solid edge $\{c,b\}$; the dotted edge
$\{a,b\}$ indicates that $a$ and $b$ are vertices growing out of
the same predecessor, and analogously with $\{c,d\}$.

The destruction of the tree of $f$ may be understood as the
destruction of the graph we have just given, because in
destructing the tree we break it at inner vertices. The
destruction of the graph is based on \emph{vertex removal} (which
one finds in Ulam's conjecture; see \cite{Ha69}, Chapter~2). More
details and a general theory concerning the destruction of graphs
may be found in \cite{DP10}.

The \textbf{S}-tree $c\cdot((b\cdot a)+d)$, which we had above,
and which corresponds to our object $f$ of
$\mathcal{WO}_u^-(X,e)$, and all the other \textbf{S}-trees, which
correspond bijectively to all the other objects of
$\mathcal{WO}_u^-(X,e)$, make the vertices of the following
hemiassociahedron:
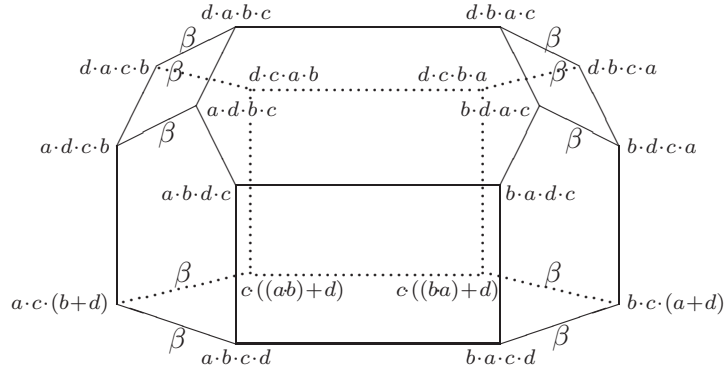
\begin{center}
\begin{picture}(200,139)(0,60)

\put(50,70){\line(-3,1){45}} \put(150,70){\line(3,1){45}}

\put(50,70){\line(1,0){100}} \put(5,85){\line(0,1){60}}
\put(195,85){\line(0,1){60}} \put(50,70){\line(0,1){60}}
\put(150,70){\line(0,1){60}} \put(50,130){\line(1,0){100}}
\put(50,130){\line(-1,2){15}} \put(150,130){\line(1,2){15}}

\put(5,145){\line(1,2){15}} \put(195,145){\line(-1,2){15}}
\put(5,145){\line(2,1){30}} \put(195,145){\line(-2,1){30}}
\put(35,160){\line(1,2){15}} \put(165,160){\line(-1,2){15}}
\put(20,175){\line(2,1){30}} \put(180,175){\line(-2,1){30}}
\put(50,190){\line(1,0){100}}

\multiput(8,85.7)(3,.7){17}{\makebox(0,0){\circle*{.5}}}
\multiput(192,85.7)(-3,.7){17}{\makebox(0,0){\circle*{.5}}}

\multiput(21.5,174.5)(3,-.7){13}{\makebox(0,0){\circle*{.5}}}
\multiput(178.5,174.5)(-3,-.7){13}{\makebox(0,0){\circle*{.5}}}
\multiput(58.3,96)(3,0){29}{\makebox(0,0){\circle*{.5}}}
\multiput(61,166)(3,0){27}{\makebox(0,0){\circle*{.5}}}
\multiput(56,99)(0,3){19}{\makebox(0,0){\circle*{.5}}}
\multiput(56,163)(0,3){2}{\makebox(0,0){\circle*{.5}}}
\multiput(144,99)(0,3){19}{\makebox(0,0){\circle*{.5}}}
\multiput(144,163)(0,3){2}{\makebox(0,0){\circle*{.5}}}

\put(198,85){\makebox(0,0)[l]{\scriptsize $b\!\cdot\! c\!\cdot\!
(a\pl d)$}}

\put(2,85){\makebox(0,0)[r]{\scriptsize $a\!\cdot\! c\!\cdot\!
(b\pl d)$}}

\put(150,68){\makebox(0,0)[t]{\scriptsize $b\!\cdot\! a\!\cdot\!
c\!\cdot\! d$}}

\put(50,68){\makebox(0,0)[t]{\scriptsize $a\!\cdot\! b\!\cdot\!
c\!\cdot\! d$}}

\put(150,94){\makebox(0,0)[tr]{\scriptsize $c\!\!\cdot\!
((b\!\!\cdot\!\! a)\pl d)$}}

\put(52,94){\makebox(0,0)[tl]{\scriptsize $c\!\!\cdot\!
((a\!\!\cdot\!\! b)\pl d)$}}

\put(145,175){\makebox(0,0)[tr]{\scriptsize $d\!\cdot\! c\!\cdot\!
b\!\cdot\! a$}}

\put(55,175){\makebox(0,0)[tl]{\scriptsize $d\!\cdot\! c\!\cdot\!
a\!\cdot\! b$}}

\put(152,130){\makebox(0,0)[tl]{\scriptsize $b\!\cdot\! a\!\cdot\!
d\!\cdot\! c$}}

\put(48,130){\makebox(0,0)[tr]{\scriptsize $a\!\cdot\! b\!\cdot\!
d\!\cdot\! c$}}

\put(198,145){\makebox(0,0)[l]{\scriptsize $b\!\cdot\! d\!\cdot\!
c\!\cdot\! a$}}

\put(2,145){\makebox(0,0)[r]{\scriptsize $a\!\cdot\! d\!\cdot\!
c\!\cdot\! b$}}

\put(183,175){\makebox(0,0)[l]{\scriptsize $d\!\cdot\! b\!\cdot\!
c\!\cdot\! a$}}

\put(17,175){\makebox(0,0)[r]{\scriptsize $d\!\cdot\! a\!\cdot\!
c\!\cdot\! b$}}

\put(50,193){\makebox(0,0)[b]{\scriptsize $d\!\cdot\! a\!\cdot\!
b\!\cdot\! c$}}

\put(150,193){\makebox(0,0)[b]{\scriptsize $d\!\cdot\! b\!\cdot\!
a\!\cdot\! c$}}

\put(39,159){\makebox(0,0)[l]{\scriptsize $a\!\cdot\! d\!\cdot\!
b\!\cdot\! c$}}

\put(161,159){\makebox(0,0)[r]{\scriptsize $b\!\cdot\! d\!\cdot\!
a\!\cdot\! c$}}

\put(28,71){\makebox(0,0){$\beta$}}
\put(30,96){\makebox(0,0){$\beta$}}
\put(24,149){\makebox(0,0){$\beta$}}
\put(32,185){\makebox(0,0){$\beta$}}
\put(27,172){\makebox(0,0){$\beta$}}

\put(173,72){\makebox(0,0){$\beta$}}
\put(170,96){\makebox(0,0){$\beta$}}
\put(178,147){\makebox(0,0){$\beta$}}
\put(170,185){\makebox(0,0){$\beta$}}
\put(172,172){\makebox(0,0){$\beta$}}
\end{picture}
\end{center}
The edges with the label $\beta$ stand for $\beta$ \emph{arrows},
i.e.\ arrows built with $\mj$, $\ins$ and one occurrence of
$\beta$ (or $\beta^{-1}$), while the remaining edges stand for
$\theta$ \emph{arrows}, i.e.\ arrows built with $\mj$, $\ins$ and
one occurrence of $\theta$.

In this hemiassociahedron, the two hexagonal faces with both
$\beta$ and $\theta$ arrows (actually, each with four $\beta$
arrows and two $\theta$ arrows) stand for the commuting diagram
corresponding to the equation ($\beta\theta 1$), or ($\beta\theta
1_e$) if we are in \WOemX. The remaining hexagonal face, which has
only $\theta$ arrows, is analogously related to the equation
($\theta$~YB), or ($\theta$~YB$_e$). The four pentagonal faces are
all of the $\beta$ and $\theta$ mixed type (each with two $\beta$
arrows and three $\theta$ arrows), and they are related to the
equation ($\beta\theta 2$), or ($\beta\theta 2_e$). The four
square faces are related to the equation ($\theta$~\emph{nat}), or
($\theta$~\emph{nat}$_e$).

This hemiassociahedron shows that the particular instance of
($\theta$~YB) that corresponds to one of the hexagonal faces is
derivable from the equations corresponding to the other faces.
However, not all instances of ($\theta$~YB) are derivable in this
manner (cf.\ Examples \ref{exa13.4} and \ref{exa13.9}). The same
applies to the equation ($\theta$~YB$_e$).

In the remaining examples we will not make comments so detailed as
here. It is however easy to recognize in the pictures of other
hemiassociahedra, and of associahedra and permutohedra, the
equations related to the two-dimensional faces. Besides equations
we had in this example, we will encounter also
($\beta$~\emph{pent}) and ($\beta$~\emph{nat}), which appear
already in the next example.

\begin{exa}\label{exa13.2}
Let $G$ and $\alpha_G$ be as in the preceding example. Let
$X=\{1111,1112,112,\linebreak 121,122,2\}$. Let $a$, $b$, $c$ and
$d$ stand respectively for 12, 1, 11 and 111, which are the inner
vertices of the following binary tree, whose leaves make $X$:
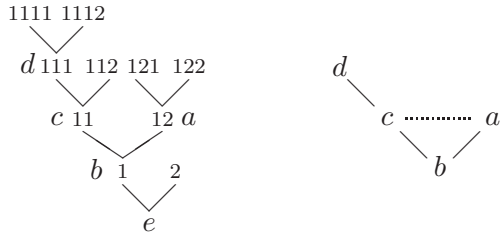
\medskip
\begin{center}
\begin{picture}(170,85)(0,3)
\put(45,10){\line(-1,1){10}} \put(45,10){\line(1,1){10}}

\put(35,30){\line(-3,2){15}} \put(35,30){\line(3,2){15}}

\put(10,70){\line(-1,1){10}} \put(10,70){\line(1,1){10}}

\put(20,50){\line(-1,1){10}} \put(20,50){\line(1,1){10}}
\put(50,50){\line(-1,1){10}} \put(50,50){\line(1,1){10}}

\put(45,5){\makebox(0,0){$e$}}
\put(35,25){\makebox(0,0){\scriptsize $1$}}
\put(25,26){\makebox(0,0){$b$}}
\put(55,25){\makebox(0,0){\scriptsize $2$}}
\put(20,45){\makebox(0,0){\scriptsize $11$}}
\put(10,45){\makebox(0,0){$c$}}
\put(50,45){\makebox(0,0){\scriptsize $12$}}
\put(60,45){\makebox(0,0){$a$}}
\put(0,85){\makebox(0,0){\scriptsize $1111$}}
\put(20,85){\makebox(0,0){\scriptsize $1112$}}
\put(10,65){\makebox(0,0){\scriptsize $111$}}
\put(-1,66){\makebox(0,0){$d$}}
\put(27,65){\makebox(0,0){\scriptsize $112$}}
\put(43,65){\makebox(0,0){\scriptsize $121$}}
\put(60,65){\makebox(0,0){\scriptsize $122$}}

\put(150,30){\line(-1,1){10}} \put(160,30){\line(1,1){10}}
\put(130,50){\line(-1,1){10}}
\multiput(142,45)(2,0){13}{\line(1,0){.5}}

\put(155,27){\makebox(0,0){$b$}} \put(135,45){\makebox(0,0){$c$}}
\put(175,45){\makebox(0,0){$a$}} \put(117,65){\makebox(0,0){$d$}}
\end{picture}
\end{center}
On the right of our tree one finds the graph of the inner
vertices.

The objects of \WOemX\ correspond bijectively to the vertices of
the following hemiassociahedron:
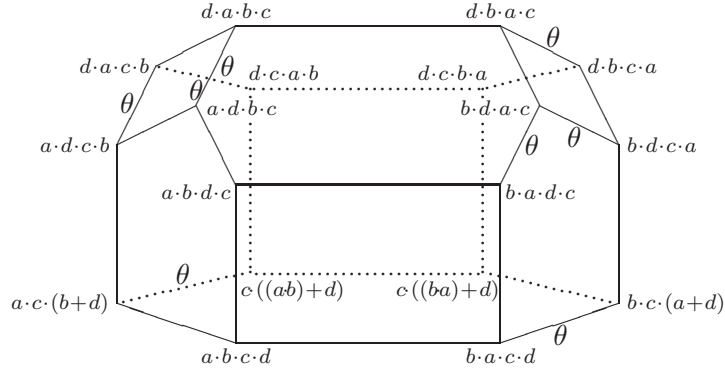
\medskip
\begin{center}
\begin{picture}(200,137)(0,60)

\put(50,70){\line(-3,1){45}} \put(150,70){\line(3,1){45}}

\put(50,70){\line(1,0){100}} \put(5,85){\line(0,1){60}}
\put(195,85){\line(0,1){60}} \put(50,70){\line(0,1){60}}
\put(150,70){\line(0,1){60}} \put(50,130){\line(1,0){100}}
\put(50,130){\line(-1,2){15}} \put(150,130){\line(1,2){15}}

\put(5,145){\line(1,2){15}} \put(195,145){\line(-1,2){15}}
\put(5,145){\line(2,1){30}} \put(195,145){\line(-2,1){30}}
\put(35,160){\line(1,2){15}} \put(165,160){\line(-1,2){15}}
\put(20,175){\line(2,1){30}} \put(180,175){\line(-2,1){30}}
\put(50,190){\line(1,0){100}}

\multiput(8,85.7)(3,.7){17}{\makebox(0,0){\circle*{.5}}}
\multiput(192,85.7)(-3,.7){17}{\makebox(0,0){\circle*{.5}}}

\multiput(21.5,174.5)(3,-.7){13}{\makebox(0,0){\circle*{.5}}}
\multiput(178.5,174.5)(-3,-.7){13}{\makebox(0,0){\circle*{.5}}}
\multiput(58.3,96)(3,0){29}{\makebox(0,0){\circle*{.5}}}
\multiput(61,166)(3,0){27}{\makebox(0,0){\circle*{.5}}}
\multiput(56,99)(0,3){19}{\makebox(0,0){\circle*{.5}}}
\multiput(56,163)(0,3){2}{\makebox(0,0){\circle*{.5}}}
\multiput(144,99)(0,3){19}{\makebox(0,0){\circle*{.5}}}
\multiput(144,163)(0,3){2}{\makebox(0,0){\circle*{.5}}}

\put(198,85){\makebox(0,0)[l]{\scriptsize $b\!\cdot\! c\!\cdot\!
(a\pl d)$}}

\put(2,85){\makebox(0,0)[r]{\scriptsize $a\!\cdot\! c\!\cdot\!
(b\pl d)$}}

\put(150,68){\makebox(0,0)[t]{\scriptsize $b\!\cdot\! a\!\cdot\!
c\!\cdot\! d$}}

\put(50,68){\makebox(0,0)[t]{\scriptsize $a\!\cdot\! b\!\cdot\!
c\!\cdot\! d$}}

\put(150,94){\makebox(0,0)[tr]{\scriptsize $c\!\!\cdot\!
((b\!\!\cdot\!\! a)\pl d)$}}

\put(52,94){\makebox(0,0)[tl]{\scriptsize $c\!\!\cdot\!
((a\!\!\cdot\!\! b)\pl d)$}}

\put(145,175){\makebox(0,0)[tr]{\scriptsize $d\!\cdot\! c\!\cdot\!
b\!\cdot\! a$}}

\put(55,175){\makebox(0,0)[tl]{\scriptsize $d\!\cdot\! c\!\cdot\!
a\!\cdot\! b$}}

\put(152,130){\makebox(0,0)[tl]{\scriptsize $b\!\cdot\! a\!\cdot\!
d\!\cdot\! c$}}

\put(48,130){\makebox(0,0)[tr]{\scriptsize $a\!\cdot\! b\!\cdot\!
d\!\cdot\! c$}}

\put(198,145){\makebox(0,0)[l]{\scriptsize $b\!\cdot\! d\!\cdot\!
c\!\cdot\! a$}}

\put(2,145){\makebox(0,0)[r]{\scriptsize $a\!\cdot\! d\!\cdot\!
c\!\cdot\! b$}}

\put(183,175){\makebox(0,0)[l]{\scriptsize $d\!\cdot\! b\!\cdot\!
c\!\cdot\! a$}}

\put(17,175){\makebox(0,0)[r]{\scriptsize $d\!\cdot\! a\!\cdot\!
c\!\cdot\! b$}}

\put(50,193){\makebox(0,0)[b]{\scriptsize $d\!\cdot\! a\!\cdot\!
b\!\cdot\! c$}}

\put(150,193){\makebox(0,0)[b]{\scriptsize $d\!\cdot\! b\!\cdot\!
a\!\cdot\! c$}}

\put(39,159){\makebox(0,0)[l]{\scriptsize $a\!\cdot\! d\!\cdot\!
b\!\cdot\! c$}}

\put(161,159){\makebox(0,0)[r]{\scriptsize $b\!\cdot\! d\!\cdot\!
a\!\cdot\! c$}}

\put(30,96){\makebox(0,0){$\theta$}}
\put(9,161){\makebox(0,0){$\theta$}}
\put(47,175){\makebox(0,0){$\theta$}}
\put(35,166){\makebox(0,0){$\theta$}}

\put(173,73){\makebox(0,0){$\theta$}}
\put(178,148){\makebox(0,0){$\theta$}}
\put(170,186){\makebox(0,0){$\theta$}}
\put(162,145){\makebox(0,0){$\theta$}}
\end{picture}
\end{center}
whose edges with $\theta$ stand for $\theta$ arrows, while the
remaining edges stand for $\beta$ arrows.
\end{exa}

\begin{exa}\label{exa13.3}
Let $G$ and $\alpha_G$ be as in the preceding two examples. Let
$X=\{1111,1112,\linebreak 1121,1122,12,2\}$. Let $a$, $b$, $c$ and
$d$ stand respectively for 111, 112, 11 and 1, which are the inner
vertices of the following binary tree, whose leaves make~$X$:
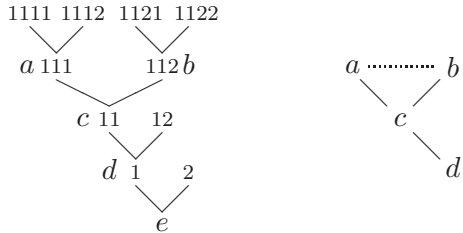
\medskip
\begin{center}
\begin{picture}(170,90)
\put(45,10){\line(-1,1){10}} \put(45,10){\line(1,1){10}}

\put(25,50){\line(-2,1){20}} \put(25,50){\line(2,1){20}}

\put(5,70){\line(-1,1){10}} \put(5,70){\line(1,1){10}}

\put(35,30){\line(-1,1){10}} \put(35,30){\line(1,1){10}}
\put(45,70){\line(-1,1){10}} \put(45,70){\line(1,1){10}}

\put(45,5){\makebox(0,0){$e$}}
\put(35,25){\makebox(0,0){\scriptsize $1$}}
\put(25,26){\makebox(0,0){$d$}}
\put(55,25){\makebox(0,0){\scriptsize $2$}}
\put(25,45){\makebox(0,0){\scriptsize $11$}}
\put(15,45){\makebox(0,0){$c$}}
\put(45,45){\makebox(0,0){\scriptsize $12$}}
\put(-5,85){\makebox(0,0){\scriptsize $1111$}}
\put(15,85){\makebox(0,0){\scriptsize $1112$}}
\put(5,65){\makebox(0,0){\scriptsize $111$}}
\put(-6,65){\makebox(0,0){$a$}}
\put(45,65){\makebox(0,0){\scriptsize $112$}}
\put(55,66){\makebox(0,0){$b$}}
\put(37,85){\makebox(0,0){\scriptsize $1121$}}
\put(58,85){\makebox(0,0){\scriptsize $1122$}}

\put(150,30){\line(-1,1){10}} \put(140,50){\line(1,1){10}}
\put(130,50){\line(-1,1){10}}
\multiput(123,65)(2,0){13}{\line(1,0){.5}}

\put(155,27){\makebox(0,0){$d$}} \put(135,45){\makebox(0,0){$c$}}
\put(155,65){\makebox(0,0){$b$}} \put(117,65){\makebox(0,0){$a$}}
\end{picture}
\end{center}
On the right of our tree one finds the graph of the inner
vertices.

The objects of \WOemX\ correspond bijectively to the vertices of
the following hemiassociahedron:
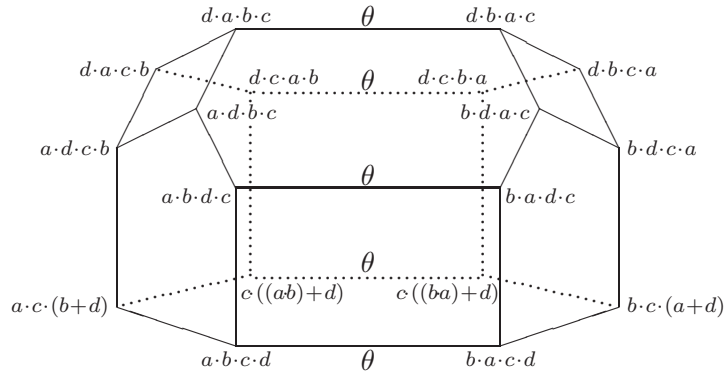
\medskip
\begin{center}
\begin{picture}(200,137)(0,60)

\put(50,70){\line(-3,1){45}} \put(150,70){\line(3,1){45}}

\put(50,70){\line(1,0){100}} \put(5,85){\line(0,1){60}}
\put(195,85){\line(0,1){60}} \put(50,70){\line(0,1){60}}
\put(150,70){\line(0,1){60}} \put(50,130){\line(1,0){100}}
\put(50,130){\line(-1,2){15}} \put(150,130){\line(1,2){15}}

\put(5,145){\line(1,2){15}} \put(195,145){\line(-1,2){15}}
\put(5,145){\line(2,1){30}} \put(195,145){\line(-2,1){30}}
\put(35,160){\line(1,2){15}} \put(165,160){\line(-1,2){15}}
\put(20,175){\line(2,1){30}} \put(180,175){\line(-2,1){30}}
\put(50,190){\line(1,0){100}}

\multiput(8,85.7)(3,.7){17}{\makebox(0,0){\circle*{.5}}}
\multiput(192,85.7)(-3,.7){17}{\makebox(0,0){\circle*{.5}}}

\multiput(21.5,174.5)(3,-.7){13}{\makebox(0,0){\circle*{.5}}}
\multiput(178.5,174.5)(-3,-.7){13}{\makebox(0,0){\circle*{.5}}}
\multiput(58.3,96)(3,0){29}{\makebox(0,0){\circle*{.5}}}
\multiput(61,166)(3,0){27}{\makebox(0,0){\circle*{.5}}}
\multiput(56,99)(0,3){19}{\makebox(0,0){\circle*{.5}}}
\multiput(56,163)(0,3){2}{\makebox(0,0){\circle*{.5}}}
\multiput(144,99)(0,3){19}{\makebox(0,0){\circle*{.5}}}
\multiput(144,163)(0,3){2}{\makebox(0,0){\circle*{.5}}}

\put(198,85){\makebox(0,0)[l]{\scriptsize $b\!\cdot\! c\!\cdot\!
(a\pl d)$}}

\put(2,85){\makebox(0,0)[r]{\scriptsize $a\!\cdot\! c\!\cdot\!
(b\pl d)$}}

\put(150,68){\makebox(0,0)[t]{\scriptsize $b\!\cdot\! a\!\cdot\!
c\!\cdot\! d$}}

\put(50,68){\makebox(0,0)[t]{\scriptsize $a\!\cdot\! b\!\cdot\!
c\!\cdot\! d$}}

\put(150,94){\makebox(0,0)[tr]{\scriptsize $c\!\!\cdot\!
((b\!\!\cdot\!\! a)\pl d)$}}

\put(52,94){\makebox(0,0)[tl]{\scriptsize $c\!\!\cdot\!
((a\!\!\cdot\!\! b)\pl d)$}}

\put(145,175){\makebox(0,0)[tr]{\scriptsize $d\!\cdot\! c\!\cdot\!
b\!\cdot\! a$}}

\put(55,175){\makebox(0,0)[tl]{\scriptsize $d\!\cdot\! c\!\cdot\!
a\!\cdot\! b$}}

\put(152,130){\makebox(0,0)[tl]{\scriptsize $b\!\cdot\! a\!\cdot\!
d\!\cdot\! c$}}

\put(48,130){\makebox(0,0)[tr]{\scriptsize $a\!\cdot\! b\!\cdot\!
d\!\cdot\! c$}}

\put(198,145){\makebox(0,0)[l]{\scriptsize $b\!\cdot\! d\!\cdot\!
c\!\cdot\! a$}}

\put(2,145){\makebox(0,0)[r]{\scriptsize $a\!\cdot\! d\!\cdot\!
c\!\cdot\! b$}}

\put(183,175){\makebox(0,0)[l]{\scriptsize $d\!\cdot\! b\!\cdot\!
c\!\cdot\! a$}}

\put(17,175){\makebox(0,0)[r]{\scriptsize $d\!\cdot\! a\!\cdot\!
c\!\cdot\! b$}}

\put(50,193){\makebox(0,0)[b]{\scriptsize $d\!\cdot\! a\!\cdot\!
b\!\cdot\! c$}}

\put(150,193){\makebox(0,0)[b]{\scriptsize $d\!\cdot\! b\!\cdot\!
a\!\cdot\! c$}}

\put(39,159){\makebox(0,0)[l]{\scriptsize $a\!\cdot\! d\!\cdot\!
b\!\cdot\! c$}}

\put(161,159){\makebox(0,0)[r]{\scriptsize $b\!\cdot\! d\!\cdot\!
a\!\cdot\! c$}}

\put(100,65){\makebox(0,0){$\theta$}}
\put(100,102){\makebox(0,0){$\theta$}}
\put(100,195){\makebox(0,0){$\theta$}}
\put(100,171){\makebox(0,0){$\theta$}}
\put(100,135){\makebox(0,0){$\theta$}}

\end{picture}
\end{center}
whose edges with $\theta$ stand for $\theta$ arrows, while the
remaining edges stand for $\beta$ arrows.
\end{exa}

\begin{exa}\label{exa13.4}
Let $G=\{x,y\}$, and let $\alpha_G(x)=2$ and $\alpha_G(y)=3$. Let
$X=\{111,112,12,\linebreak 21,22,31,32\}$. Let $a$, $b$, $c$ and
$d$ stand respectively for 3, 2, 1 and 11, which are the inner
vertices of the following tree, whose leaves make~$X$:
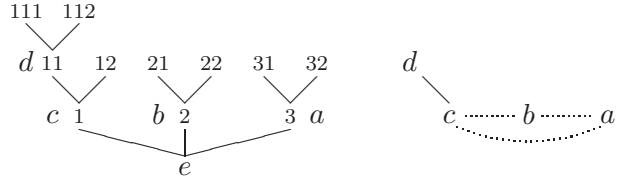
\medskip
\begin{center}
\begin{picture}(220,70)
\put(60,10){\line(-4,1){40}} \put(60,10){\line(0,1){10}}
\put(60,10){\line(4,1){40}}

\put(20,30){\line(-1,1){10}} \put(20,30){\line(1,1){10}}
\put(60,30){\line(-1,1){10}} \put(60,30){\line(1,1){10}}
\put(100,30){\line(-1,1){10}} \put(100,30){\line(1,1){10}}

\put(10,50){\line(-1,1){10}} \put(10,50){\line(1,1){10}}

\put(60,5){\makebox(0,0){$e$}}
\put(20,25){\makebox(0,0){\scriptsize $1$}}
\put(10,25){\makebox(0,0){$c$}}
\put(60,25){\makebox(0,0){\scriptsize $2$}}
\put(50,26){\makebox(0,0){$b$}}
\put(100,25){\makebox(0,0){\scriptsize $3$}}
\put(110,25){\makebox(0,0){$a$}}
\put(10,45){\makebox(0,0){\scriptsize $11$}}
\put(0,46){\makebox(0,0){$d$}}
\put(30,45){\makebox(0,0){\scriptsize $12$}}
\put(50,45){\makebox(0,0){\scriptsize $21$}}
\put(70,45){\makebox(0,0){\scriptsize $22$}}
\put(90,45){\makebox(0,0){\scriptsize $31$}}
\put(110,45){\makebox(0,0){\scriptsize $32$}}
\put(0,65){\makebox(0,0){\scriptsize $111$}}
\put(20,65){\makebox(0,0){\scriptsize $112$}}

\put(160,30){\line(-1,1){10}}
\multiput(166,25)(2,0){10}{\line(1,0){.5}}
\multiput(195,25)(2,0){10}{\line(1,0){.5}}

\put(160,25){\makebox(0,0){$c$}} \put(190,26){\makebox(0,0){$b$}}
\put(220,25){\makebox(0,0){$a$}} \put(145,46){\makebox(0,0){$d$}}

\qbezier[25](163,21)(190,10)(217,21)
\end{picture}
\end{center}
On the right of our tree one finds the graph of the inner
vertices.

The objects of \WOemX\ correspond bijectively to the vertices of
the following hemiassociahedron:
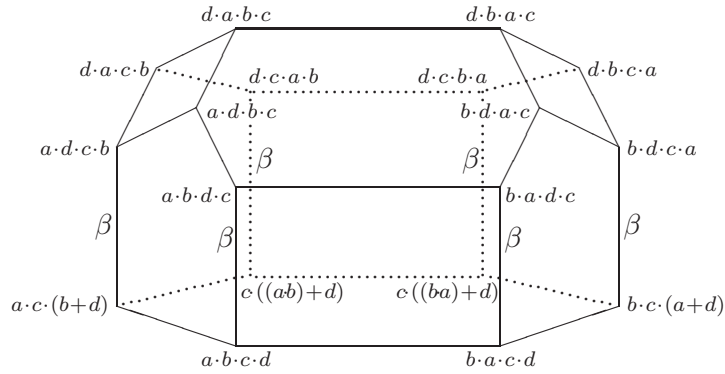
\medskip
\begin{center}
\begin{picture}(200,137)(0,60)

\put(50,70){\line(-3,1){45}} \put(150,70){\line(3,1){45}}

\put(50,70){\line(1,0){100}} \put(5,85){\line(0,1){60}}
\put(195,85){\line(0,1){60}} \put(50,70){\line(0,1){60}}
\put(150,70){\line(0,1){60}} \put(50,130){\line(1,0){100}}
\put(50,130){\line(-1,2){15}} \put(150,130){\line(1,2){15}}

\put(5,145){\line(1,2){15}} \put(195,145){\line(-1,2){15}}
\put(5,145){\line(2,1){30}} \put(195,145){\line(-2,1){30}}
\put(35,160){\line(1,2){15}} \put(165,160){\line(-1,2){15}}
\put(20,175){\line(2,1){30}} \put(180,175){\line(-2,1){30}}
\put(50,190){\line(1,0){100}}

\multiput(8,85.7)(3,.7){17}{\makebox(0,0){\circle*{.5}}}
\multiput(192,85.7)(-3,.7){17}{\makebox(0,0){\circle*{.5}}}

\multiput(21.5,174.5)(3,-.7){13}{\makebox(0,0){\circle*{.5}}}
\multiput(178.5,174.5)(-3,-.7){13}{\makebox(0,0){\circle*{.5}}}
\multiput(58.3,96)(3,0){29}{\makebox(0,0){\circle*{.5}}}
\multiput(61,166)(3,0){27}{\makebox(0,0){\circle*{.5}}}
\multiput(56,99)(0,3){19}{\makebox(0,0){\circle*{.5}}}
\multiput(56,163)(0,3){2}{\makebox(0,0){\circle*{.5}}}
\multiput(144,99)(0,3){19}{\makebox(0,0){\circle*{.5}}}
\multiput(144,163)(0,3){2}{\makebox(0,0){\circle*{.5}}}

\put(198,85){\makebox(0,0)[l]{\scriptsize $b\!\cdot\! c\!\cdot\!
(a\pl d)$}}

\put(2,85){\makebox(0,0)[r]{\scriptsize $a\!\cdot\! c\!\cdot\!
(b\pl d)$}}

\put(150,68){\makebox(0,0)[t]{\scriptsize $b\!\cdot\! a\!\cdot\!
c\!\cdot\! d$}}

\put(50,68){\makebox(0,0)[t]{\scriptsize $a\!\cdot\! b\!\cdot\!
c\!\cdot\! d$}}

\put(150,94){\makebox(0,0)[tr]{\scriptsize $c\!\!\cdot\!
((b\!\!\cdot\!\! a)\pl d)$}}

\put(52,94){\makebox(0,0)[tl]{\scriptsize $c\!\!\cdot\!
((a\!\!\cdot\!\! b)\pl d)$}}

\put(145,175){\makebox(0,0)[tr]{\scriptsize $d\!\cdot\! c\!\cdot\!
b\!\cdot\! a$}}

\put(55,175){\makebox(0,0)[tl]{\scriptsize $d\!\cdot\! c\!\cdot\!
a\!\cdot\! b$}}

\put(152,130){\makebox(0,0)[tl]{\scriptsize $b\!\cdot\! a\!\cdot\!
d\!\cdot\! c$}}

\put(48,130){\makebox(0,0)[tr]{\scriptsize $a\!\cdot\! b\!\cdot\!
d\!\cdot\! c$}}

\put(198,145){\makebox(0,0)[l]{\scriptsize $b\!\cdot\! d\!\cdot\!
c\!\cdot\! a$}}

\put(2,145){\makebox(0,0)[r]{\scriptsize $a\!\cdot\! d\!\cdot\!
c\!\cdot\! b$}}

\put(183,175){\makebox(0,0)[l]{\scriptsize $d\!\cdot\! b\!\cdot\!
c\!\cdot\! a$}}

\put(17,175){\makebox(0,0)[r]{\scriptsize $d\!\cdot\! a\!\cdot\!
c\!\cdot\! b$}}

\put(50,193){\makebox(0,0)[b]{\scriptsize $d\!\cdot\! a\!\cdot\!
b\!\cdot\! c$}}

\put(150,193){\makebox(0,0)[b]{\scriptsize $d\!\cdot\! b\!\cdot\!
a\!\cdot\! c$}}

\put(39,159){\makebox(0,0)[l]{\scriptsize $a\!\cdot\! d\!\cdot\!
b\!\cdot\! c$}}

\put(161,159){\makebox(0,0)[r]{\scriptsize $b\!\cdot\! d\!\cdot\!
a\!\cdot\! c$}}

\put(0,115){\makebox(0,0){$\beta$}}
\put(46,110){\makebox(0,0){$\beta$}}
\put(61,140){\makebox(0,0){$\beta$}}
\put(139,140){\makebox(0,0){$\beta$}}
\put(155,110){\makebox(0,0){$\beta$}}
\put(200,115){\makebox(0,0){$\beta$}}

\end{picture}
\end{center}
whose edges with $\beta$ stand for $\beta$ arrows, while the
remaining edges stand for $\theta$ arrows.
\end{exa}

\begin{exa}\label{exa13.5}
Let $G$ and $\alpha_G$ be as in Examples
\ref{exa13.1}-\ref{exa13.3}. Let $X=\{1111, 1112, 112,
12,\linebreak 21,22\}$. Let $a$, $b$, $c$ and $d$ stand
respectively for 2, 1, 11 and 111, which are the inner vertices of
the following binary tree, whose leaves make~$X$:
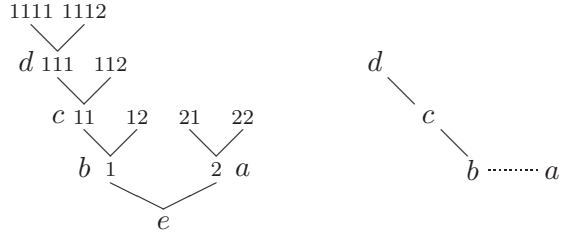
\medskip
\begin{center}
\begin{picture}(200,90)
\put(45,10){\line(-2,1){20}} \put(45,10){\line(2,1){20}}

\put(25,30){\line(-1,1){10}} \put(25,30){\line(1,1){10}}

\put(15,50){\line(-1,1){10}} \put(15,50){\line(1,1){10}}

\put(5,70){\line(-1,1){10}} \put(5,70){\line(1,1){10}}
\put(65,30){\line(-1,1){10}} \put(65,30){\line(1,1){10}}

\put(45,5){\makebox(0,0){$e$}}
\put(25,25){\makebox(0,0){\scriptsize $1$}}
\put(15,26){\makebox(0,0){$b$}}
\put(65,25){\makebox(0,0){\scriptsize $2$}}
\put(75,25){\makebox(0,0){$a$}}
\put(15,45){\makebox(0,0){\scriptsize $11$}}
\put(5,45){\makebox(0,0){$c$}}
\put(35,45){\makebox(0,0){\scriptsize $12$}}
\put(55,45){\makebox(0,0){\scriptsize $21$}}
\put(75,45){\makebox(0,0){\scriptsize $22$}}
\put(5,65){\makebox(0,0){\scriptsize $111$}}
\put(-7,66){\makebox(0,0){$d$}}
\put(25,65){\makebox(0,0){\scriptsize $112$}}
\put(-5,85){\makebox(0,0){\scriptsize $1111$}}
\put(15,85){\makebox(0,0){\scriptsize $1112$}}

\put(160,30){\line(-1,1){10}} \put(140,50){\line(-1,1){10}}
\multiput(168,25)(2,0){10}{\line(1,0){.5}}

\put(162,25){\makebox(0,0){$b$}} \put(192,24){\makebox(0,0){$a$}}
\put(145,45){\makebox(0,0){$c$}} \put(125,66){\makebox(0,0){$d$}}

\end{picture}
\end{center}
On the right of our tree one finds the graph of the inner
vertices.

The objects of \WOemX\ correspond bijectively to the vertices of
the following three-dimensional associahedron:\medskip
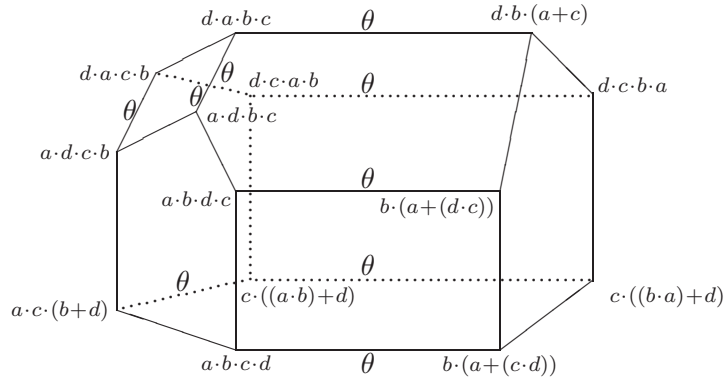
\begin{center}
\begin{picture}(200,140)(0,60)

\put(50,70){\line(-3,1){45}}

\put(50,70){\line(1,0){100}} \put(5,85){\line(0,1){60}}

\put(50,70){\line(0,1){60}} \put(150,70){\line(0,1){60}}
\put(50,130){\line(1,0){100}} \put(50,130){\line(-1,2){15}}
\put(5,145){\line(1,2){15}} \put(5,145){\line(2,1){30}}
\put(35,160){\line(1,2){15}} \put(20,175){\line(2,1){30}}
\put(50,190){\line(1,0){112.2}} \put(150,70){\line(4,3){35}}
\put(185,96){\line(0,1){71.3}}

\put(150,130){\line(1,5){11.9}} \put(162,190){\line(1,-1){23}}

\multiput(8,85.7)(3,.7){17}{\makebox(0,0){\circle*{.5}}}

\multiput(21.5,174.5)(3,-.7){13}{\makebox(0,0){\circle*{.5}}}

\multiput(58.3,96.5)(3,0){43}{\makebox(0,0){\circle*{.5}}}
\multiput(61,166)(3,0){42}{\makebox(0,0){\circle*{.5}}}
\multiput(56,99)(0,3){19}{\makebox(0,0){\circle*{.5}}}
\multiput(56,163)(0,3){2}{\makebox(0,0){\circle*{.5}}}

\put(2,85){\makebox(0,0)[r]{\scriptsize $a\!\cdot\! c\!\cdot\!
(b\pl d)$}}

\put(150,68){\makebox(0,0)[t]{\scriptsize $b\!\cdot\! (a\pl
(c\!\cdot\! d))$}}

\put(50,68){\makebox(0,0)[t]{\scriptsize $a\!\cdot\! b\!\cdot\!
c\!\cdot\! d$}}

\put(235,94){\makebox(0,0)[tr]{\scriptsize $c\!\cdot\!
((b\!\cdot\! a)\pl d)$}}

\put(52,94){\makebox(0,0)[tl]{\scriptsize $c\!\cdot\! ((a\!\cdot\!
b)\pl d)$}}

\put(187,173){\makebox(0,0)[tl]{\scriptsize $d\!\cdot\! c\!\cdot\!
b\!\cdot\! a$}}

\put(55,175){\makebox(0,0)[tl]{\scriptsize $d\!\cdot\! c\!\cdot\!
a\!\cdot\! b$}}

\put(148,128){\makebox(0,0)[tr]{\scriptsize $b\!\cdot\!
(a\pl(d\!\cdot\! c))$}}

\put(48,130){\makebox(0,0)[tr]{\scriptsize $a\!\cdot\! b\!\cdot\!
d\!\cdot\! c$}}

\put(2,145){\makebox(0,0)[r]{\scriptsize $a\!\cdot\! d\!\cdot\!
c\!\cdot\! b$}}

\put(17,175){\makebox(0,0)[r]{\scriptsize $d\!\cdot\! a\!\cdot\!
c\!\cdot\! b$}}

\put(50,193){\makebox(0,0)[b]{\scriptsize $d\!\cdot\! a\!\cdot\!
b\!\cdot\! c$}}

\put(165,193){\makebox(0,0)[b]{\scriptsize $d\!\cdot\! b\!\cdot\!
(a\pl c)$}}

\put(39,159){\makebox(0,0)[l]{\scriptsize $a\!\cdot\! d\!\cdot\!
b\!\cdot\! c$}}

\put(30,96){\makebox(0,0){$\theta$}}
\put(47,175){\makebox(0,0){$\theta$}}
\put(35,166){\makebox(0,0){$\theta$}}
\put(10,161){\makebox(0,0){$\theta$}}

\put(100,65){\makebox(0,0){$\theta$}}
\put(100,102){\makebox(0,0){$\theta$}}
\put(100,195){\makebox(0,0){$\theta$}}
\put(100,171){\makebox(0,0){$\theta$}}
\put(100,135){\makebox(0,0){$\theta$}}
\end{picture}
\end{center}
whose edges with $\theta$ stand for $\theta$ arrows, while the
remaining edges stand for $\beta$ arrows.
\end{exa}

\begin{exa}\label{exa13.6}
Let $G$ and
$\alpha_G$ be as in Examples \ref{exa13.1}-\ref{exa13.3} and the
preceding example. Let $X=\{111,112,12,21,221,222\}$. Let $a$,
$b$, $c$ and $d$ stand respectively for 11, 1, 2 and 22, which are
the inner vertices of the following binary tree, whose leaves
make~$X$:
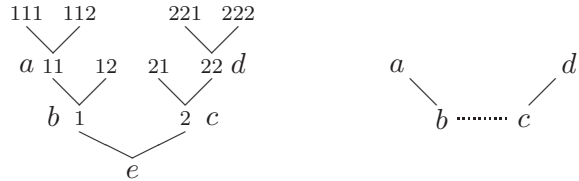
\medskip
\begin{center}
\begin{picture}(200,70)
\put(45,10){\line(-2,1){20}} \put(45,10){\line(2,1){20}}

\put(25,30){\line(-1,1){10}} \put(25,30){\line(1,1){10}}

\put(15,50){\line(-1,1){10}} \put(15,50){\line(1,1){10}}

\put(75,50){\line(-1,1){10}} \put(75,50){\line(1,1){10}}
\put(65,30){\line(-1,1){10}} \put(65,30){\line(1,1){10}}

\put(45,5){\makebox(0,0){$e$}}
\put(25,25){\makebox(0,0){\scriptsize $1$}}
\put(15,26){\makebox(0,0){$b$}}
\put(65,25){\makebox(0,0){\scriptsize $2$}}
\put(75,25){\makebox(0,0){$c$}}
\put(15,45){\makebox(0,0){\scriptsize $11$}}
\put(5,45){\makebox(0,0){$a$}}
\put(35,45){\makebox(0,0){\scriptsize $12$}}
\put(55,45){\makebox(0,0){\scriptsize $21$}}
\put(75,45){\makebox(0,0){\scriptsize $22$}}
\put(85,46){\makebox(0,0){$d$}}
\put(5,65){\makebox(0,0){\scriptsize $111$}}
\put(25,65){\makebox(0,0){\scriptsize $112$}}
\put(65,65){\makebox(0,0){\scriptsize $221$}}
\put(85,65){\makebox(0,0){\scriptsize $222$}}

\put(160,30){\line(-1,1){10}} \put(195,30){\line(1,1){10}}
\multiput(168,25)(2,0){10}{\line(1,0){.5}}

\put(162,25){\makebox(0,0){$b$}} \put(193,24){\makebox(0,0){$c$}}
\put(145,45){\makebox(0,0){$a$}} \put(210,46){\makebox(0,0){$d$}}

\end{picture}
\end{center}
On the right of our tree one finds the graph of the inner
vertices.

The objects of \WOemX\ correspond bijectively to the vertices of
the following three-dimensional associahedron:\medskip
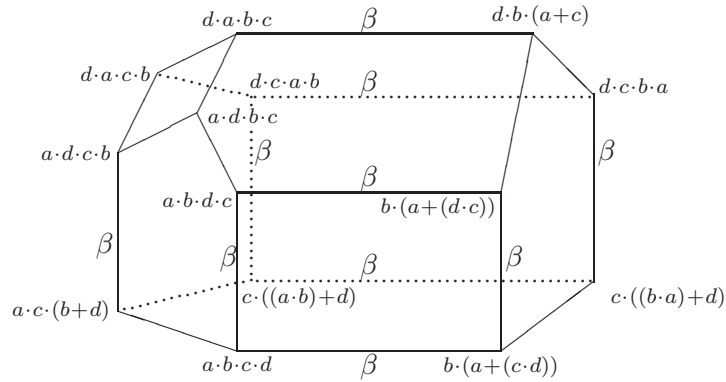
\begin{center}
\begin{picture}(200,140)(0,60)

\put(50,70){\line(-3,1){45}}

\put(50,70){\line(1,0){100}} \put(5,85){\line(0,1){60}}

\put(50,70){\line(0,1){60}} \put(150,70){\line(0,1){60}}
\put(50,130){\line(1,0){100}} \put(50,130){\line(-1,2){15}}
\put(5,145){\line(1,2){15}} \put(5,145){\line(2,1){30}}
\put(35,160){\line(1,2){15}} \put(20,175){\line(2,1){30}}
\put(50,190){\line(1,0){112.2}} \put(150,70){\line(4,3){35}}
\put(185,96){\line(0,1){71.3}}

\put(150,130){\line(1,5){11.9}} \put(162,190){\line(1,-1){23}}

\multiput(8,85.7)(3,.7){17}{\makebox(0,0){\circle*{.5}}}

\multiput(21.5,174.5)(3,-.7){13}{\makebox(0,0){\circle*{.5}}}

\multiput(58.3,96.5)(3,0){43}{\makebox(0,0){\circle*{.5}}}
\multiput(61,166)(3,0){42}{\makebox(0,0){\circle*{.5}}}
\multiput(56,99)(0,3){19}{\makebox(0,0){\circle*{.5}}}
\multiput(56,163)(0,3){2}{\makebox(0,0){\circle*{.5}}}

\put(2,85){\makebox(0,0)[r]{\scriptsize $a\!\cdot\! c\!\cdot\!
(b\pl d)$}}

\put(150,68){\makebox(0,0)[t]{\scriptsize $b\!\cdot\! (a\pl
(c\!\cdot\! d))$}}

\put(50,68){\makebox(0,0)[t]{\scriptsize $a\!\cdot\! b\!\cdot\!
c\!\cdot\! d$}}

\put(235,94){\makebox(0,0)[tr]{\scriptsize $c\!\cdot\!
((b\!\cdot\! a)\pl d)$}}

\put(52,94){\makebox(0,0)[tl]{\scriptsize $c\!\cdot\! ((a\!\cdot\!
b)\pl d)$}}

\put(187,173){\makebox(0,0)[tl]{\scriptsize $d\!\cdot\! c\!\cdot\!
b\!\cdot\! a$}}

\put(55,175){\makebox(0,0)[tl]{\scriptsize $d\!\cdot\! c\!\cdot\!
a\!\cdot\! b$}}

\put(148,128){\makebox(0,0)[tr]{\scriptsize $b\!\cdot\!
(a\pl(d\!\cdot\! c))$}}

\put(48,130){\makebox(0,0)[tr]{\scriptsize $a\!\cdot\! b\!\cdot\!
d\!\cdot\! c$}}

\put(2,145){\makebox(0,0)[r]{\scriptsize $a\!\cdot\! d\!\cdot\!
c\!\cdot\! b$}}

\put(17,175){\makebox(0,0)[r]{\scriptsize $d\!\cdot\! a\!\cdot\!
c\!\cdot\! b$}}

\put(50,193){\makebox(0,0)[b]{\scriptsize $d\!\cdot\! a\!\cdot\!
b\!\cdot\! c$}}

\put(165,193){\makebox(0,0)[b]{\scriptsize $d\!\cdot\! b\!\cdot\!
(a\pl c)$}}

\put(39,159){\makebox(0,0)[l]{\scriptsize $a\!\cdot\! d\!\cdot\!
b\!\cdot\! c$}}

\put(0,110){\makebox(0,0){$\beta$}}
\put(47,105){\makebox(0,0){$\beta$}}
\put(60,145){\makebox(0,0){$\beta$}}
\put(155,105){\makebox(0,0){$\beta$}}
\put(190,145){\makebox(0,0){$\beta$}}

\put(100,64){\makebox(0,0){$\beta$}}
\put(100,102){\makebox(0,0){$\beta$}}
\put(100,195){\makebox(0,0){$\beta$}}
\put(100,171){\makebox(0,0){$\beta$}}
\put(100,135){\makebox(0,0){$\beta$}}
\end{picture}
\end{center}
whose edges with $\beta$ stand for $\beta$ arrows, while the
remaining edges stand for $\theta$ arrows.
\end{exa}

\begin{exa}\label{exa13.7}
Let $G=\{x,y\}$, and let $\alpha_G(x)=2$ and $\alpha_G(y)=3$. Let
$X=\{111,112,121,\linebreak 122,131,132,2\}$. Let $a$, $b$, $c$
and $d$ stand respectively for 1, 11, 12 and 13, which are the
inner vertices of the following tree, whose leaves make~$X$:
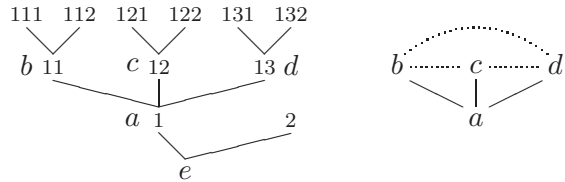
\medskip
\begin{center}
\begin{picture}(200,70)
\put(60,10){\line(-1,1){10}} \put(60,10){\line(4,1){40}}

\put(50,30){\line(-4,1){40}} \put(50,30){\line(0,1){10}}
\put(50,30){\line(4,1){40}}

\put(10,50){\line(-1,1){10}} \put(10,50){\line(1,1){10}}

\put(50,50){\line(-1,1){10}} \put(50,50){\line(1,1){10}}
\put(90,50){\line(-1,1){10}} \put(90,50){\line(1,1){10}}

\put(60,5){\makebox(0,0){$e$}}
\put(50,25){\makebox(0,0){\scriptsize $1$}}
\put(40,25){\makebox(0,0){$a$}}
\put(100,25){\makebox(0,0){\scriptsize $2$}}
\put(10,45){\makebox(0,0){\scriptsize $11$}}
\put(0,46){\makebox(0,0){$b$}}
\put(50,45){\makebox(0,0){\scriptsize $12$}}
\put(40,46){\makebox(0,0){$c$}}
\put(90,45){\makebox(0,0){\scriptsize $13$}}
\put(100,46){\makebox(0,0){$d$}}
\put(0,65){\makebox(0,0){\scriptsize $111$}}
\put(20,65){\makebox(0,0){\scriptsize $112$}}
\put(40,65){\makebox(0,0){\scriptsize $121$}}
\put(60,65){\makebox(0,0){\scriptsize $122$}}
\put(80,65){\makebox(0,0){\scriptsize $131$}}
\put(100,65){\makebox(0,0){\scriptsize $132$}}

\put(165,30){\line(-2,1){20}} \put(175,30){\line(2,1){20}}
\put(170,30){\line(0,1){10}}
\multiput(145,45)(2,0){10}{\line(1,0){.5}}
\multiput(175,45)(2,0){10}{\line(1,0){.5}}
\qbezier[25](143,50)(170,70)(197,50)

\put(170,25){\makebox(0,0){$a$}} \put(140,46){\makebox(0,0){$b$}}
\put(170,45){\makebox(0,0){$c$}} \put(200,46){\makebox(0,0){$d$}}

\end{picture}
\end{center}
On the right of our tree one finds the graph of the inner
vertices.

The objects of \WOemX\ correspond bijectively to the vertices of
the following three-dimensional permutohedron:\medskip
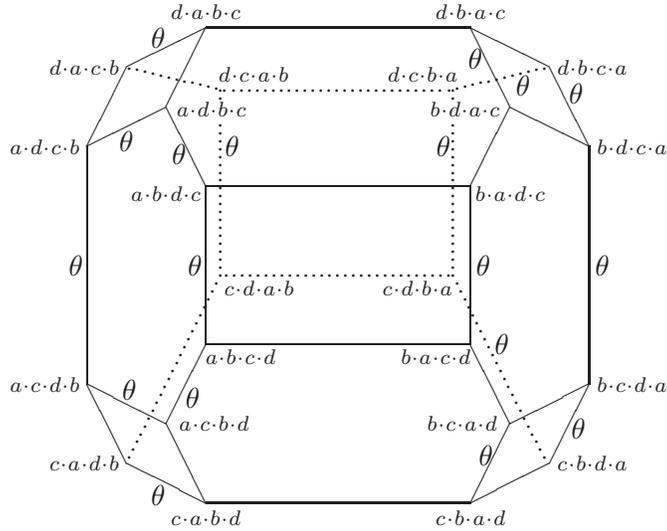
\begin{center}
\begin{picture}(200,220)(0,-10)

\put(50,10){\line(1,0){100}} \put(50,10){\line(-2,1){30}}
\put(150,10){\line(2,1){30}} \put(50,10){\line(-1,2){15}}
\put(150,10){\line(1,2){15}} \put(20,25){\line(-1,2){15}}
\put(180,25){\line(1,2){15}} \put(35,40){\line(-2,1){30}}
\put(165,40){\line(2,1){30}}

\put(35,40){\line(1,2){15}} \put(165,40){\line(-1,2){15}}
\put(50,70){\line(1,0){100}} \put(5,55){\line(0,1){90}}
\put(195,55){\line(0,1){90}} \put(50,70){\line(0,1){60}}
\put(150,70){\line(0,1){60}} \put(50,130){\line(1,0){100}}
\put(50,130){\line(-1,2){15}} \put(150,130){\line(1,2){15}}

\put(5,145){\line(1,2){15}} \put(195,145){\line(-1,2){15}}
\put(5,145){\line(2,1){30}} \put(195,145){\line(-2,1){30}}
\put(35,160){\line(1,2){15}} \put(165,160){\line(-1,2){15}}
\put(20,175){\line(2,1){30}} \put(180,175){\line(-2,1){30}}
\put(50,190){\line(1,0){100}}

\multiput(21.5,27)(1.5,3){24}{\makebox(0,0){\circle*{.5}}}
\multiput(178.5,27)(-1.5,3){24}{\makebox(0,0){\circle*{.5}}}

\multiput(21.5,174.5)(3,-.7){13}{\makebox(0,0){\circle*{.5}}}
\multiput(178.5,174.5)(-3,-.7){13}{\makebox(0,0){\circle*{.5}}}
\multiput(58.3,96)(3,0){29}{\makebox(0,0){\circle*{.5}}}
\multiput(61,166)(3,0){27}{\makebox(0,0){\circle*{.5}}}
\multiput(56,99)(0,3){19}{\makebox(0,0){\circle*{.5}}}
\multiput(56,163)(0,3){2}{\makebox(0,0){\circle*{.5}}}
\multiput(144,99)(0,3){19}{\makebox(0,0){\circle*{.5}}}
\multiput(144,163)(0,3){2}{\makebox(0,0){\circle*{.5}}}

\put(50,7){\makebox(0,0)[t]{\scriptsize $c\!\cdot\! a\!\cdot\!
b\!\cdot\! d$}}

\put(150,7){\makebox(0,0)[t]{\scriptsize $c\!\cdot\! b\!\cdot\!
a\!\cdot\! d$}}

\put(183,25){\makebox(0,0)[l]{\scriptsize $c\!\cdot\! b\!\cdot\!
d\!\cdot\! a$}}

\put(17,25){\makebox(0,0)[r]{\scriptsize $c\!\cdot\! a\!\cdot\!
d\!\cdot\! b$}}

\put(198,55){\makebox(0,0)[l]{\scriptsize $b\!\cdot\! c\!\cdot\!
d\!\cdot\! a$}}

\put(2,55){\makebox(0,0)[r]{\scriptsize $a\!\cdot\! c\!\cdot\!
d\!\cdot\! b$}}

\put(160,40){\makebox(0,0)[r]{\scriptsize $b\!\cdot\! c\!\cdot\!
a\!\cdot\! d$}}

\put(40,40){\makebox(0,0)[l]{\scriptsize $a\!\cdot\! c\!\cdot\!
b\!\cdot\! d$}}

\put(150,68){\makebox(0,0)[tr]{\scriptsize $b\!\cdot\! a\!\cdot\!
c\!\cdot\! d$}}

\put(50,68){\makebox(0,0)[tl]{\scriptsize $a\!\cdot\! b\!\cdot\!
c\!\cdot\! d$}}

\put(143,94){\makebox(0,0)[tr]{\scriptsize $c\!\cdot\! d\!\cdot\!
b\!\cdot\! a$}}

\put(57,94){\makebox(0,0)[tl]{\scriptsize $c\!\cdot\! d\!\cdot\!
a\!\cdot\! b$}}

\put(145,175){\makebox(0,0)[tr]{\scriptsize $d\!\cdot\! c\!\cdot\!
b\!\cdot\! a$}}

\put(55,175){\makebox(0,0)[tl]{\scriptsize $d\!\cdot\! c\!\cdot\!
a\!\cdot\! b$}}

\put(152,130){\makebox(0,0)[tl]{\scriptsize $b\!\cdot\! a\!\cdot\!
d\!\cdot\! c$}}

\put(48,130){\makebox(0,0)[tr]{\scriptsize $a\!\cdot\! b\!\cdot\!
d\!\cdot\! c$}}

\put(198,145){\makebox(0,0)[l]{\scriptsize $b\!\cdot\! d\!\cdot\!
c\!\cdot\! a$}}

\put(2,145){\makebox(0,0)[r]{\scriptsize $a\!\cdot\! d\!\cdot\!
c\!\cdot\! b$}}

\put(183,175){\makebox(0,0)[l]{\scriptsize $d\!\cdot\! b\!\cdot\!
c\!\cdot\! a$}}

\put(17,175){\makebox(0,0)[r]{\scriptsize $d\!\cdot\! a\!\cdot\!
c\!\cdot\! b$}}

\put(50,193){\makebox(0,0)[b]{\scriptsize $d\!\cdot\! a\!\cdot\!
b\!\cdot\! c$}}

\put(150,193){\makebox(0,0)[b]{\scriptsize $d\!\cdot\! b\!\cdot\!
a\!\cdot\! c$}}

\put(39,159){\makebox(0,0)[l]{\scriptsize $a\!\cdot\! d\!\cdot\!
b\!\cdot\! c$}}

\put(161,159){\makebox(0,0)[r]{\scriptsize $b\!\cdot\! d\!\cdot\!
a\!\cdot\! c$}}

\put(32,13){\makebox(0,0){$\theta$}}
\put(21,53){\makebox(0,0){$\theta$}}
\put(1,100){\makebox(0,0){$\theta$}}
\put(46,100){\makebox(0,0){$\theta$}}
\put(45,50){\makebox(0,0){$\theta$}}
\put(40,142){\makebox(0,0){$\theta$}}
\put(20,147){\makebox(0,0){$\theta$}}
\put(32,186){\makebox(0,0){$\theta$}}
\put(60,145){\makebox(0,0){$\theta$}}
\put(140,145){\makebox(0,0){$\theta$}}
\put(152,178){\makebox(0,0){$\theta$}}
\put(190,165){\makebox(0,0){$\theta$}}
\put(170,168){\makebox(0,0){$\theta$}}
\put(200,100){\makebox(0,0){$\theta$}}
\put(155,100){\makebox(0,0){$\theta$}}
\put(156,28){\makebox(0,0){$\theta$}}
\put(191,38){\makebox(0,0){$\theta$}}
\put(162,70){\makebox(0,0){$\theta$}}
\end{picture}
\end{center}
whose edges with $\theta$ stand for $\theta$ arrows, while the
remaining edges stand for $\beta$ arrows.
\end{exa}

\begin{exa}\label{exa13.8}Let $G$ and
$\alpha_G$ be as in Examples \ref{exa13.1}-\ref{exa13.3} and
\ref{exa13.5}-\ref{exa13.6}. Let
$X=\{11111,11112,1112,112,12,2\}$. Let $a$, $b$, $c$ and $d$ stand
respectively for 1, 11, 111 and 1111, which are the inner vertices
of the following binary tree, whose leaves make~$X$:
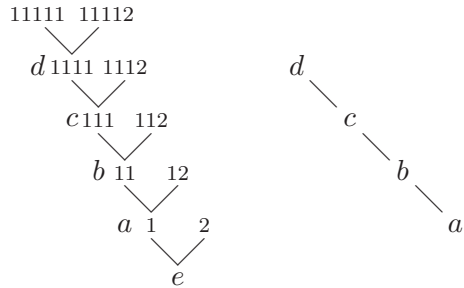
\medskip
\begin{center}
\begin{picture}(170,110)
\put(60,10){\line(-1,1){10}} \put(60,10){\line(1,1){10}}

\put(50,30){\line(-1,1){10}} \put(50,30){\line(1,1){10}}

\put(40,50){\line(-1,1){10}} \put(40,50){\line(1,1){10}}

\put(30,70){\line(-1,1){10}} \put(30,70){\line(1,1){10}}
\put(20,90){\line(-1,1){10}} \put(20,90){\line(1,1){10}}

\put(60,5){\makebox(0,0){$e$}}
\put(50,25){\makebox(0,0){\scriptsize $1$}}
\put(40,25){\makebox(0,0){$a$}}
\put(70,25){\makebox(0,0){\scriptsize $2$}}
\put(40,45){\makebox(0,0){\scriptsize $11$}}
\put(30,46){\makebox(0,0){$b$}}
\put(60,45){\makebox(0,0){\scriptsize $12$}}
\put(30,65){\makebox(0,0){\scriptsize $111$}}
\put(20,65){\makebox(0,0){$c$}}
\put(50,65){\makebox(0,0){\scriptsize $112$}}
\put(20,85){\makebox(0,0){\scriptsize $1111$}}
\put(7,86){\makebox(0,0){$d$}}
\put(40,85){\makebox(0,0){\scriptsize $1112$}}
\put(7,105){\makebox(0,0){\scriptsize $11111$}}
\put(33,105){\makebox(0,0){\scriptsize $11112$}}

\put(160,30){\line(-1,1){10}} \put(140,50){\line(-1,1){10}}
\put(120,70){\line(-1,1){10}}

\put(165,25){\makebox(0,0){$a$}} \put(145,46){\makebox(0,0){$b$}}
\put(125,65){\makebox(0,0){$c$}} \put(105,86){\makebox(0,0){$d$}}

\end{picture}
\end{center}
On the right of our tree one finds the graph of the inner
vertices.

The objects of \WOemX\ correspond bijectively to the vertices of
the three-dimensional associahedron whose edges all stand for
$\beta$ arrows (its picture is like that for Examples
\ref{exa13.5} and \ref{exa13.6} without the labels $\theta$ and
$\beta$).
\end{exa}

\begin{exa}\label{exa13.9}
Let $G=\{x,y\}$, and let $\alpha_G(x)=2$ and $\alpha_G(y)=4$. Let
$X=\{11,12,21,22,\linebreak 31,32,41,42\}$. Let $a$, $b$, $c$ and
$d$ stand respectively for 1, 2, 3 and 4, which are the inner
vertices of the following tree, whose leaves make~$X$:
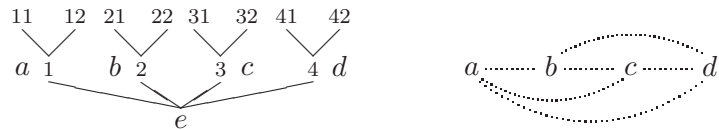
\medskip
\begin{center}
\begin{picture}(260,50)
\put(60,10){\line(-5,1){50}} \put(60,10){\line(5,1){50}}

\put(60,10){\line(-3,2){15}} \put(60,10){\line(3,2){15}}

\put(10,30){\line(-1,1){10}} \put(10,30){\line(1,1){10}}
\put(45,30){\line(-1,1){10}} \put(45,30){\line(1,1){10}}
\put(75,30){\line(-1,1){10}} \put(75,30){\line(1,1){10}}
\put(110,30){\line(-1,1){10}} \put(110,30){\line(1,1){10}}

\put(60,5){\makebox(0,0){$e$}}
\put(10,25){\makebox(0,0){\scriptsize $1$}}
\put(0,25){\makebox(0,0){$a$}}
\put(45,25){\makebox(0,0){\scriptsize $2$}}
\put(35,26){\makebox(0,0){$b$}}
\put(75,25){\makebox(0,0){\scriptsize $3$}}
\put(85,25){\makebox(0,0){$c$}}
\put(110,25){\makebox(0,0){\scriptsize $4$}}
\put(120,26){\makebox(0,0){$d$}}
\put(0,45){\makebox(0,0){\scriptsize $11$}}
\put(20,45){\makebox(0,0){\scriptsize $12$}}
\put(35,45){\makebox(0,0){\scriptsize $21$}}
\put(53,45){\makebox(0,0){\scriptsize $22$}}
\put(67,45){\makebox(0,0){\scriptsize $31$}}
\put(85,45){\makebox(0,0){\scriptsize $32$}}
\put(100,45){\makebox(0,0){\scriptsize $41$}}
\put(120,45){\makebox(0,0){\scriptsize $42$}}

\multiput(175,25)(2,0){10}{\line(1,0){.5}}
\multiput(205,25)(2,0){10}{\line(1,0){.5}}
\multiput(235,25)(2,0){10}{\line(1,0){.5}}
\qbezier[40](173,20)(215,-10)(257,20)
\qbezier[30](174,21)(200,5)(227,21)
\qbezier[30](204,31)(230,45)(257,31)

\put(170,25){\makebox(0,0){$a$}} \put(200,26){\makebox(0,0){$b$}}
\put(230,25){\makebox(0,0){$c$}} \put(260,26){\makebox(0,0){$d$}}

\end{picture}
\end{center}
On the right of our tree one finds the graph of the inner
vertices.

The objects of \WOemX\ correspond bijectively to the vertices of
the three-dimensional permutohedron whose edges all stand for
$\theta$ arrows (its picture is like that for
Example~\ref{exa13.7} without labels for edges).
\end{exa}

\vspace{4ex}

{\samepage \noindent{\large {\sc Part III}}

\section{Coherence of \Monu}}\label{sec14}
In \cite{MLP85} (Section~2) it was established that the notion of
bicategory is coherent, in the sense that all diagrams of
canonical arrows commute. The proof of this coherence result is
obtained by imitating the proof of monoidal coherence (see
\cite{ML63} \cite{ML98}, Section VII.2).

Our purpose in this, concluding part, of the paper is to prove an
analogous coherence result for our notion of weak Cat-operad,
which is analogous to the notion of bicategory. This amounts to
showing that the category \WOe\ is a preorder (i.e.\ that all
diagrams commute in this category, or that for a given source and
target there is not more than one arrow). This coherence result
for \WOe\ does not rely only on monoidal coherence, as coherence
for bicategories does. Besides relying on monoidal coherence, it
relies also on a generalization of coherence for symmetric
monoidal categories, which is related to the presentation of
structures related to symmetric groups.

Coherence for \WOe\ is a justification of our definition of this
category, and of our notion of weak Cat-operad. We will establish
the coherence of \WOe\ by establishing the coherence of \WOu. We
first introduce a category \Monu\ derived from \WOu, which is
analogous to a monoidal category.

The category \Monu\ is defined like \WOu\ save that we omit the
basic arrow terms $\theta_{h,g,f}$ and the axiomatic equations
that involve $\theta$ explicitly (these are ($\theta$~\emph{nat}),
($\theta\theta$), ($\theta$~YB), ($\beta\theta 1$) and
($\beta\theta 2$)). The remaining axiomatic equations are
analogous to Mac Lane's postulates for monoidal categories (see \cite{ML63} and
\cite{ML98}, Section VII.1); the difference is that $\ins$ is a
partial operation on objects and on arrows, and is hence not a
real biendofunctor, though it is analogous to such a functor.

Let an arrow term of \Monu\ be called \emph{directed} when
$\beta$, $\mu^{-1}$ and $\lambda^{-1}$ do not occur in it (but
$\beta^{-1}$, $\mu$ and $\lambda$ may occur). Let us call an
object $f$ of \Monu, i.e.\ term of \Ou, \emph{normal} when either
all parentheses in $f$ are associated to the left and $\Iota$ does
not occur in $f$, or $f$ is of the form $a\cdot\Iota$. One can
then prove the following.

\begin{lem}\label{lem14.1}
If $u$ and $v$ are directed arrow terms of \Monu\ of the same type
with a normal target, then ${u=v}$ in \Monu.
\end{lem}

\begin{proof}
The proof of this lemma is obtained by imitating a part of
the proof of monoidal coherence in \cite{ML98} (Section VII.2,
Theorem 1; see also \cite{DP04}, Sections 4.2, 4.3 and 4.6,
Directedness Lemmata). There is nothing essentially new in this
inductive proof, which consists in showing a kind of confluence
property, related to what one has in term-rewriting. For example,
if $u$ and $v$ are respectively of the form
$u'\cirk(\mj_j\ins\beta^{-1}_{h,g,f})$ and
$v'\cirk\beta^{-1}_{j,h,g\ins f}$, then, since the target of $u$
and $v$ is normal, there is a $w$ such that by the induction
hypothesis
\begin{tabbing}
\hspace{4em}\=$u'=w\cirk(\beta^{-1}_{j,h,g}\ins\mj_f)\cirk\beta^{-1}_{j,h\ins
g,f}$,
\\[.5ex]
\>$v'=w\cirk\beta^{-1}_{j,h,g\ins f}$;
\end{tabbing}
and then by using the equation obtained from ($\beta$~\emph{pent})
and ($\beta\beta$), analogous to Mac Lane's pentagonal diagram, we
obtain ${u=v}$.
\end{proof}

We can then establish the following.

\begin{prop}\label{prop14.2}
\Monu\ is a preorder.
\end{prop}

\noindent To prove this proposition we may proceed as for
Associative Coherence in \cite{DP04} (Section~4.3).

\section{The category \WOuth}\label{sec15}
Next we define a category \WOuth, which is \WOu\ strictified in
the monoidal structure. This means that in \WOuth\ the arrows
$\beta$, $\mu$ and $\lambda$, and their inverses, become all
identity arrows, and the monoidal structure of \WOuth\ is trivial.
As in the case of monoidal categories, and as in \cite{DP04}
(Section 3.2), by relying on Proposition~\ref{prop14.2} one can
show that \WOu\ and \WOuth\ are equivalent categories.

Now we will define \WOuth\ syntactically. Its objects are the
normal terms of \Ou\ as defined in the preceding section. We may
identify these terms by terms where $\ins$ and parentheses are
deleted, and in which $\Iota$ does not occur, except if the term
is of the form $a\cdot\Iota$. We use this abbreviated notation
below.

As basic arrow terms we have the following:
\begin{tabbing}
\hspace{4em}\=$\mj_f\!:f\str f$,\quad for every object $f$,
\\[.5ex]
\>$\theta_{h,g,f}\!:hgf\str hfg$,\quad provided $t(f)\in s(h)$ and
$t(g)\in s(h)$.
\end{tabbing}
Next we have the operations on arrow terms as for \WOu\ (see
Section~7) save that $v\ins u\!:g\ins f\str g'\ins f'$ is written
${vu\!: gf\str g'f'}$. (We write $v\cirk u$ as before.)

Besides $u=u$ and the categorial equations, and the equations
(\emph{ins}~1), (\emph{ins}~2), ($\theta$~\emph{nat}),
($\theta\theta$) and ($\theta$~YB) (all written without $\ins$),
the \emph{axiomatic equations} of \WOuth\ are the following
strictified versions of ($\beta\theta 1$) and ($\beta\theta 2$):
\begin{tabbing}
\hspace{4em}\=($\theta
1$)\hspace{2em}\=$\mj_j\theta_{h,g,f}=\theta_{jh,g,f}$,
\\[1ex]
\>($\theta
2$)\>$\theta_{j,hg,f}=(\theta_{j,h,f}\mj_g)\cirk\theta_{jh,g,f}$.
\end{tabbing}
(Note that if the left-hand sides of these equations are
legitimate, then the right-hand sides are legitimate too, but not
conversely.)

The rules of the equational axiomatic system are the same as for
\WOu. This concludes the definition of the category \WOuth.

\section{\Cg\ and \BCg}\label{sec16}
We introduce a family of categories we call \Cg\ and show in the
next section that every category in the family is coherent in a
sense to be made precise. A particular category in the family, for
a particular choice of $\Gamma$, will be shown isomorphic to the
category \WOuth\ of the preceding section, and this will establish
that \WOuth\ is coherent, which in this particular case implies
that \WOuth\ is a preorder. Symmetric groups arise as particular
members of this family, and the proof that \Cg\ is coherent will
proceed as a proof that would show the completeness of a standard
presentation of symmetric groups (see the next section).

Now we introduce \Cg. Let $A$, $B$, $P$, $Q$, $R$, $S$,
$U,\ldots,$ perhaps with indices, stand for finite (possibly
empty) sequences, i.e.\ for words, in an alphabet whose members we
call \emph{atoms}; we use $p$, $q$, $r,\ldots,$ perhaps with
indices, for atoms. We use $e$, as before, for the empty word. The
set of objects of \Cg\ is some set of these words, not necessarily
all. (So the objects of \Cg\ make a language.)

The \emph{basic} arrow terms of \Cg\ make a set $\Gamma$, which
satisfies the following. For every object $A$, the arrow term
$A\!:A\str A$ is in $\Gamma$ (we abbreviate $\mj_A$ by writing
just $A$). We have in $\Gamma$ also some arrow terms of the form
\begin{tabbing}
\hspace{4em}$A[p,q]B\!:ApqB\str AqpB$,
\end{tabbing}
provided both the source and the target are objects. All the arrow
terms in $\Gamma$ are of these two kinds. Finally, $\Gamma$ must
satisfy the following condition:
\begin{itemize}[label=($\Gamma$)]
\item if the basic arrow terms on one side of the
equations ($C1$) and ($C2$) below are in $\Gamma$, then all the
basic arrow terms on the other side of these equations are in
$\Gamma$ too.
\end{itemize}
(It is natural to call the arrow terms $A\!:A\str A$ basic, though
they could have been left out from $\Gamma$, and introduced a bit
later, to produce \Cg.)

All the arrow terms of \Cg\ are defined by starting from $\Gamma$,
and closing under composition: if ${u\!:A\str B}$ and ${v\!:B\str
C}$ are arrow terms, then ${v\cirk u\!:A\str C}$ is an arrow term.
We use $u$, $v$, $w,\ldots,$ perhaps with indices, for arrow
terms, and we use the abbreviation given by the following
inductive clause:
\begin{tabbing}
\hspace{4em}$A(v\cirk u)B=_{df}AvB\cirk AuB$.
\end{tabbing}

Besides ${u=u}$ and the categorial equations $u\cirk A=u=B\cirk
u$, for ${u\!:A\str B}$, and $(w\cirk v)\cirk u=w\cirk(v\cirk u)$,
the \emph{axiomatic equations} of \Cg\ are the following:
\begin{tabbing}
\hspace{0em}\=($C
1$)\hspace{2em}\=$A(srU[p,q]\cirk[r,s]Upq)B=A([r,s]Uqp\cirk
rsU[p,q])B$,
\\[.5ex]
\>($C2$)\>$A(s[p,r]\cirk[p,s]r\cirk p[r,s])B=A([r,s]p\cirk
r[p,s]\cirk[p,r]s)B$,
\\[.5ex]
\>($C3$)\>$A([p,r]\cirk[r,p])B=ArpB$.
\end{tabbing}
As rules we have symmetry and transitivity of $=$ and congruence
for $\cirk$ (see Section~7). This concludes our definition of the
equations of \Cg, and of the category \Cg.

The axiomatic equations of \Cg\ are analogous to the equations of
the standard presentation of the symmetric group $S_n$, for $n\geq
1$, with the generators $\sigma_i$, for $1\leq i<n$, being the
transpositions of $i$ and $i\pl 1$ (see \cite{CM57}, Section 6.2).
The equation ($C1$) corresponds to the permutability of $\sigma_i$
and $\sigma_j$ when ${i\mn j}$ is at least 2. The equation ($C2$)
corresponds to the equation
\begin{tabbing}
\hspace{0em}\=($C1$)\hspace{2em}\=\kill
\>(YB)\>$\sigma_{i+1}\cirk\sigma_i\cirk\sigma_{i+1}=\sigma_i\cirk\sigma_{i+1}\cirk\sigma_i$
\end{tabbing}
(YB comes from Yang-Baxter), and ($C3$) corresponds to the
$\sigma_i$'s being self-inverse.

The symmetric group $S_n$ is \Cg\ that has a unique object $p^n$,
which is a sequence of $n$ occurrences of $p$, and the set
$\Gamma$ is made of the arrow terms $p^{i-1}[p,p]p^k$, where
$i\geq 1$, $k\geq 0$ and ${i\pl 1\pl k=n}$, which correspond to
$\sigma_i$. Note that in general \Cg\ is not a group. It need not
even be a groupoid (in the categorial sense, a Brandt groupoid;
see \cite{ML98}, Section I.5); we may have $A[p,q]B$ in $\Gamma$
without having its inverse $A[q,p]B$.

To reduce the arrow terms of \Cg\ to normal form we introduce the
category \BCg, a variant of \Cg, which we will show isomorphic to
\Cg. The objects of \BCg\ are those of \Cg. The arrow terms of
\BCg\ are defined starting from the same basic arrow terms
$\Gamma$, and closing under composition and under the following:
\begin{tabbing}
\hspace{0em}\=($C1$)\hspace{1em}\=\kill
\>($\dagger$)\>if for every $s$ in the word $S$, such that $S$ is $S'sS''$, the arrow term\\[.5ex]
\hspace{6em}$AS'[r,s]S''B\!:AS'rsS''B\str AS'srS''B$\\[.5ex]
\>\>is in $\Gamma$, then\\[.5ex]
\hspace{6em}$A[r,S]B\!:ArSB\str ASrB$\\[.5ex]
\>\>is an arrow term.
\end{tabbing}
Note that according to this clause $A[r,e]B\!:ArB\str ArB$ is
always an arrow term. Note also that according to this definition we also have the implication converse to ($\dagger$).

The equations of \BCg\ are defined like those of \Cg\ save that
the axiomatic equations ($C1$), ($C2$) and ($C3$) are replaced by
the following axiomatic equations:
\begin{tabbing}
\hspace{0em}\=($C1$)\hspace{2em}\=\kill
\>($BC1$)\>$A(SrU[p,Q]\cirk[r,S]UpQ)B\;$\=$=A([r,S]UQp\cirk
rSU[p,Q])B$,
\\*[.5ex]
\` provided neither $Q$ nor $S$ is $e$,
\\[1ex]
\>($BC2$)\>$A([p,QSrU]\cirk
pQ[r,S]U)B$\>$=A(Q[r,S]Up\cirk[p,QrSU])B$,
\hspace{.5em} provided $S$ is not $e$,
\\[1ex]
\>($BC3$)\>$A(S[p,UrQ]\cirk[r,SpU]Q)B$\>$=A([r,SU]Qp\cirk
rS[p,UQ])B$,
\\[1ex]
\>($BC4$)\>$A(S[r,Q]\cirk[r,S]Q)B=A[r,SQ]B$,\hspace{.5em} provided
neither $Q$ nor $S$ is $e$,
\\[1ex]
\>($BC5$)\>$A[r,e]B=ArB$.
\end{tabbing}
This concludes the definition of \BCg.

The axiomatic equations of \BCg\ are analogous to the equations
that may be found in \cite{DP04} (Section 5.2). These equations
are such as to enable us to reach quickly a normal form for arrow
terms, with which we will deal in the next section. In the
remainder of this section we will establish that \Cg\ and \BCg\
are isomorphic.

We show first that we have in \Cg\ the structure of \BCg. We
define $A[r,S]B$ in \Cg\ by the following inductive clauses:
\begin{tabbing}
\hspace{4em}\=$A[r,e]B=_{df}ArB$,
\\[.5ex]
\>$A[r,sQ]B=_{df}A(s[r,Q]\cirk[r,s]Q)B$.
\end{tabbing}
Then it remains to derive the equations ($BC1$)-($BC5$) in \Cg.

We derive first ($BC1$) by induction on the sum $n$ of the lengths
of $Q$ and $S$. In the basis when $n$ is 2, we use ($C1$), and in
the induction step we just use the induction hypothesis.

Next we derive by induction on the length $n$ of $S$ the equation
($C2\;S$), which is ($C2$) with $s$ replaced by $S$. The basis,
when $n$ is 0, is trivial. In the induction step, to derive our
equation (read from left to right) we use ($BC1$) (read from left
to right), the induction hypothesis and ($C2$)(read from left to
right).

To derive ($BC2$) (read from left to right) we use ($BC1$) (read
in both directions) and ($C2\;S$) (read from left to right).

To derive ($BC3$) (read from left to right) we use ($C2\;S$),
($C3$) and ($BC1$) (all read from left to right).

The equations ($BC4$) and ($BC5$) hold in \Cg\ by definition. This
establishes that we have the structure of \BCg\ in \Cg.

We have noted above in parentheses when we needed the equations
($C2$) and ($C3$) only from left to right (while ($C1$), via
($BC1$), is needed in both directions). This may be interesting
when our procedure is connected with the reduction procedure of
\cite{La03} (Section 2.1 and Appendix~A).

To establish the converse --- namely, that we have the structure
of \Cg\ in \BCg\ --- is an easy matter. The equations ($C1$),
($C2$) and ($C3$) amount to particular cases of ($BC1$), ($BC2$)
and ($BC3$), with the help of ($BC4$) and ($BC5$). The definitions
of $A[r,e]B$ and $A[r,sQ]B$, which we introduced in \Cg, clearly
hold in \BCg\ by ($BC5$) and ($BC4$). So we may conclude that \Cg\
and \BCg\ are isomorphic.

As a consequence of the isomorphism of \Cg\ and \BCg\ we obtain
that if on one side of the equations ($BC1$)-($BC5$) we have arrow
terms of \BCg, then on the other side we have such arrow terms
too. The straightforward proof of that is based essentially on
($\dagger$), the implication converse to ($\dagger$), and the fact that if on one side of an equation of \Cg\ we have arrow
terms of \Cg, then on the other side we have such arrow terms
too.

\section{Coherence of \Cg}\label{sec17}
We say that an arrow term of \BCg\ is in \emph{normal form} when
it is of the form
\begin{tabbing}
\hspace{4em}$A_1Q_1p_1B_1\cirk A_1[p_1,Q_1]B_1\cirk\ldots\cirk
A_n[p_n,Q_n]B_n$,
\end{tabbing}
for $n\geq 0$, for the words $Q_1,\ldots,Q_n$ nonempty and for
\begin{tabbing}
\hspace{4em}$|B_1|>|B_2|>\ldots>|B_n|$,
\end{tabbing}
where $|B_i|$ is the length of $B_i$. This, or an analogous normal
form, for symmetric groups is implicit in \cite{M96} and
\cite{B11} (Note~C), and occurs explicitly in \cite{La95} (Section
3.2), \cite{La03} (Section 2.1) and \cite{DP04} (Section 5.2).

A survey of all possible cases shows that if an arrow term of
\BCg\ is not in normal form, then it has, after perhaps applying
categorial equations, a subterm of the form of the left-hand side
of one of the equations ($BC1$)-($BC5$), and hence one of these
equations may be applied. We can then establish the following.

\begin{lem}\label{lem17.1}
Every arrow term of \BCg\ is equal in \BCg\ to an arrow term in
normal form.
\end{lem}

\begin{proof}
For every arrow term of the form
\begin{tabbing}
\hspace{4em}$C_m[r_m,S_m]D_m\cirk\ldots\cirk C_1[r_1,S_1]D_1$,
\end{tabbing}
where $m\geq 2$, consider the following measure of this arrow
term:
\begin{tabbing}
\hspace{4em}$\displaystyle \sum_{i=1}^m(|C_i|\pl 1\pl|S_i|)\cdot
i$.
\end{tabbing}

Then it can be checked that with each application of
($BC1$)-($BC4$) from left to right the measure decreases. The
equation ($BC5$) from left to right works together with the
categorial equations $u\cirk A=u=B\cirk u$ to reduce our measure.
\mbox{\hspace{1em}}
\end{proof}

The \emph{graph} of an arrow term of \BCg\ is derived from a
bijection between finite ordinals, defined as for symmetric groups
or as for symmetric monoidal categories (see \cite{DP04}). This
bijection induces a graph with edges connecting an occurrence of
an atom in the source to an occurrence of the same atom in the
target. The graphs of $A[p,q]B$, $A[p,Q]B$ and of the identity
arrow $A\!:A\str A$ are given by
\begin{center}
\begin{picture}(265,40)
\put(0,10){\line(0,1){20}} \put(15,10){\line(0,1){20}}
\put(25,10){\line(1,2){10}} \put(35,10){\line(-1,2){10}}
\put(50,10){\line(0,1){20}} \put(65,10){\line(0,1){20}}

\put(7,5){\makebox(0,0){$A$}} \put(25,5){\makebox(0,0){$q$}}
\put(35,5){\makebox(0,0){$p$}} \put(57,5){\makebox(0,0){$B$}}
\put(7,35){\makebox(0,0){$A$}} \put(25,35){\makebox(0,0){$p$}}
\put(35,35){\makebox(0,0){$q$}} \put(57,35){\makebox(0,0){$B$}}

\put(7.5,20){\makebox(0,0){\scriptsize $\cdots$}}
\put(57.5,20){\makebox(0,0){\scriptsize $\cdots$}}

\put(120,10){\line(0,1){20}} \put(135,10){\line(0,1){20}}
\put(143,10){\line(1,2){10}} \put(155,10){\line(1,2){10}}
\put(165,10){\line(-1,1){20}} \put(180,10){\line(0,1){20}}
\put(195,10){\line(0,1){20}}

\put(127,5){\makebox(0,0){$A$}} \put(149,5){\makebox(0,0){$Q$}}
\put(165,5){\makebox(0,0){$p$}} \put(187,5){\makebox(0,0){$B$}}
\put(127,35){\makebox(0,0){$A$}} \put(145,35){\makebox(0,0){$p$}}
\put(161,35){\makebox(0,0){$Q$}} \put(187,35){\makebox(0,0){$B$}}

\put(127.5,20){\makebox(0,0){\scriptsize $\cdots$}}
\put(187.5,20){\makebox(0,0){\scriptsize $\cdots$}}
\put(157.5,27){\makebox(0,0){\tiny $\cdots$}}

\put(250,10){\line(0,1){20}} \put(265,10){\line(0,1){20}}
\put(257,5){\makebox(0,0){$A$}} \put(257,35){\makebox(0,0){$A$}}
\put(257.5,20){\makebox(0,0){\scriptsize $\cdots$}}

\end{picture}
\end{center}
The graphs for arrow terms with $\cirk$ are obtained by composing
the underlying bijections. We can then prove the following.

\begin{lem}\label{lem17.2}
If the arrow terms $u,v\!:A\str B$ of \BCg\ are
in normal form and their graphs are the same, then $u$ and $v$ are
the same arrow term.
\end{lem}

\begin{proof}
The proof of this lemma is analogous to the proof of the
Uniqueness Lemma in \cite{DP04} (Section 5.2), to which we refer
for details. Here is a sketch of the proof.

Let $u$ and $v$ be respectively the arrow terms
\begin{tabbing}
\hspace{4em}\=$A\cirk A_1[p_1,Q_1]B_1\cirk\ldots\cirk
A_n[p_n,Q_n]B_n$,
\\[.5ex]
\>$A\cirk C_1[r_1,S_1]D_1\cirk\ldots\cirk C_m[r_m,S_m]D_m$
\end{tabbing}
in normal form. We proceed by induction on $n$. If $n=0$, then we
show that $m$ must be 0 too; otherwise the graphs would differ.

If $n>0$, then we must have $m>0$ too, as we have just shown, and
$A_n[p_n,Q_n]B_n$ must be equal to $C_m[r_m,S_m]D_m$; otherwise
the graphs would differ in the edges of $p_n$ and $r_m$:
\begin{center}
\begin{picture}(140,40)
\put(20,10){\line(-1,1){20}} \put(30,10){\line(0,1){20}}
\put(45,10){\line(0,1){20}}

\put(21,5){\makebox(0,0){$p_n$}} \put(39,5){\makebox(0,0){$B_n$}}
\put(0,35){\makebox(0,0){$p_n$}} \put(39,35){\makebox(0,0){$B_n$}}

\put(37.5,20){\makebox(0,0){\scriptsize $\cdots$}}
\put(0,5){\makebox(0,0){\scriptsize $\cdots$}}
\put(20,35){\makebox(0,0){\scriptsize $\cdots$}}

\put(120,10){\line(-1,1){20}} \put(130,10){\line(0,1){20}}
\put(145,10){\line(0,1){20}}

\put(121,5){\makebox(0,0){$r_m$}}
\put(139,5){\makebox(0,0){$D_m$}}
\put(100,35){\makebox(0,0){$r_m$}}
\put(139,35){\makebox(0,0){$D_m$}}

\put(137.5,20){\makebox(0,0){\scriptsize $\cdots$}}
\put(100,5){\makebox(0,0){\scriptsize $\cdots$}}
\put(120,35){\makebox(0,0){\scriptsize $\cdots$}}

\end{picture}
\end{center}
(see \cite{DP04}, Section 5.2, for details). We then conclude that
\begin{tabbing}
\hspace{4em}\=$A\cirk A_1[p_1,Q_1]B_1\cirk\ldots\cirk
A_{n-1}[p_{n-1},Q_{n-1}]B_{n-1}$,
\\[.5ex]
\>$A\cirk C_1[r_1,S_1]D_1\cirk\ldots\cirk
C_{m-1}[r_{m-1},S_{m-1}]D_{m-1}$
\end{tabbing}
must have the same graph, and we apply the induction hypothesis to
them.
\end{proof}

Then we obtain the following

\begin{prop}\label{prop17.3}
For $u$ and $v$ arrow terms of \BCg\ of the same type we have
$u=v$ in \BCg\ iff $u$ and $v$ have the same graph.
\end{prop}

\begin{proof}
It is easy to establish the implication from left to right by
induction on the length of derivation. For the converse
implication, we have by Lemma~\ref{lem17.1} that $u=u'$ and $v=v'$
in \BCg\ for $u'$ and $v'$ in normal form. From the assumption
that $u$ and $v$ have the same graph we conclude, by the
implication from left to right, that $u'$ and $v'$ have the same
graph. But then by Lemma~\ref{lem17.2} we have that $u'$ is $v'$,
and hence $u=v$ in \BCg\ by symmetry and transitivity of~$=$.
\end{proof}

From Proposition~\ref{prop17.3} from left to right and
Lemma~\ref{lem17.2} we may infer that for every arrow term $u$ of
\BCg\ there is a unique arrow term $u'$ of \BCg\ in normal form
such that $u=u'$ in \BCg. Note that we did not need this kind of
uniqueness proposition to establish coherence. Note also that we
did not establish this uniqueness proposition by means like
confluence of term rewriting. Instead we have a uniqueness
proposition, our Lemma~\ref{lem17.2}, which involves graphs.
Establishing uniqueness in the latter way may often be more easier
than doing it in the former one. The paper \cite{La03} (Section
2.1 and Appendix~A) considers uniqueness of normal form
established via confluence of term rewriting.

As an alternative to the style of proof of coherence of this paper
there is the style of the original paper of \cite{M96}, which one
finds also in \cite{B11}. This style works for symmetric groups,
and relies on the fact that in that case one can establish that
for the symmetric group $S_n$ there are $n!$ normal forms and $n!$
permutations. Then it is enough to establish that the map from the
syntax for $S_n$ to permutations is \emph{onto}, which means that
every permutation is represented by a term of $S_n$, i.e.\ a
composition of generators. (We must also establish that the map in
question is a homomorphism, which means that every equation of the
syntax holds for permutations.) It follows then that the map is
one-one.

This old style of argument seems however too complicated and
chaotic in the case of \Cg. It may work for a regularly chosen
$\Gamma$, but $\Gamma$ may be irregular, and it is not clear with
what the number $n!$ can be replaced. It would be a number lesser
than or equal to $n!$, but it may change irregularly, and preclude
an inductive argument on $n$.

Since \Cg\ and \BCg\ are isomorphic, we may establish
Proposition~\ref{prop17.3} for \Cg, and we call that proposition
the \emph{Coherence of} \Cg. (The graphs for \Cg\ are obtained in
an obvious manner through the isomorphism with \BCg.)

\section{\Cg\ and \WOuth\ --- Coherence of \WOuth}\label{sec18}
Let us now show that \WOuth\ may be conceived as a category \Cg.
The atoms are all the terms $a\cdot x$ and $a\cdot \Iota$ of \Ou.
The objects of \WOuth\ are the normal terms of \Ou, which may be
identified with some words made of atoms. The set $\Gamma$ is made
first of the arrow terms $\mj_f$ for every object $f$ of \WOuth;
here $\mj_f$ stands for ${f\!:f\str f}$. Next we have in $\Gamma$
all the arrow terms $\theta_{h,a\cdot x,b\cdot y}\mj_j$ and
$\theta_{h,a\cdot x,b\cdot y}$ of \WOuth; they stand for
\begin{tabbing}
\hspace{4em}$h[a\cdot x,b\cdot y]j\!:h(a\cdot x)(b\cdot y)j\str
h(b\cdot y)(a\cdot x)j$
\end{tabbing}
and the same without~$j$.

We use (\emph{ins}~1), (\emph{ins}~2), ($\theta 1$) and ($\theta
2$) to reduce every arrow term $u$ of \WOuth\ to the form
$u_n\cirk\ldots\cirk u_1$ where each $u_i$ is either an identity
${f\!:f\str f}$, or $\theta_{h,a\cdot x,b\cdot y}\mj_j$, or
$\theta_{h,a\cdot x,b\cdot y}$. We check next that we have the
equations of \Cg\ in \WOuth. We have of course $u=u$ and the
categorial equations. For ($C1$) we use essentially
($\theta$~\emph{nat}), for ($C2$) we use essentially
($\theta$~YB), and for ($C3$) we use essentially ($\theta\theta$);
for all that we need also (\emph{ins}~2). To show that,
conversely, all the equations of \WOuth may be derived from the
\Cg\ assumptions is a consequence of the coherence of \Cg.
Finally, we check easily that the set $\Gamma$ of \WOuth\
satisfies condition ($\Gamma$) (see Section 16).

So \WOuth\ is a \Cg\ category, and hence coherence for \Cg\ holds
for it. But in this particular case coherence for \Cg\ becomes the
following.

\begin{prop}\label{prop18.1}
For every arrow terms $u$ and $v$ of \WOuth\ of the same type we
have $u=v$ in \WOuth.
\end{prop}

\noindent In other words, the category \WOuth\ is a preorder. This
proposition is a consequence of the coherence of \Cg\ and of the
fact that in \WOuth\ the type of an arrow term determines uniquely
the graph. The reason for that is that in the set $\Gamma$ of
\WOuth\ we do not have $h[a\cdot x,a\cdot x]j$ and $h[a\cdot
x,a\cdot x] $.

From Proposition~\ref{prop18.1} and the equivalence of \WOu\ with
\WOuth\ we may conclude that \WOu\ is a preorder. From that and
from Proposition~\ref{prop11.6} we conclude that the category
\WOe\ is a preorder. This establishes the coherence of our notion
of weak Cat-operad of Section 12. This notion is coherent in the
same sense in which Mac Lane' s notion of monoidal category and
the notion of bicategory are coherent. All diagrams of canonical
arrows commute in it.

\section*{Acknowledgment}
  \noindent Work on this paper was
supported by the Ministry of Science of Serbia (Grant ON174026).

\end{document}